\documentclass[10pt]{article}
\usepackage{amsmath}
\usepackage{amsfonts}
\usepackage{mathrsfs}
\usepackage{amssymb}
\usepackage{amsthm}
\usepackage{palatino}
\usepackage[margin=2.5cm, vmargin={1.5cm},includefoot]{geometry}

\usepackage{pst-plot}
\usepackage{pst-grad}
\usepackage{pstcol}

\newrgbcolor{LemonChiffon}{1.0 0.98 0.8}
\newrgbcolor{magenta}{1.0 0.3 0.4}
\newrgbcolor{mistyrose}{1.0 0.0 0.6}
\newrgbcolor{deeppink}{1.0 0.08 0.58}
\newrgbcolor{lightsalmon}{1.0 0.63 0.48}
\newrgbcolor{lightred}{0.7 0.0 0.0}
\newrgbcolor{lightblue}{0.68 0.85 0.90}
\newrgbcolor{indianred}{0.80  0.36  0.36}
\newrgbcolor{lightgreen}{0.56 0.93 0.56}
\newrgbcolor{byellow}{0.93 0.76 0.3}

\DeclareMathOperator{\cl}{Cl_2}

\begin{document}

\title{Asymptotics for the partial fractions of the restricted partition generating function II}

\author{Cormac O'Sullivan\footnote{
\newline
Support for this project was provided by a PSC-CUNY Award, jointly funded by The Professional Staff Congress and The City University of New York.}}

\maketitle


\begin{abstract}
The generating function for $p_N(n)$, the number of  partitions of $n$ into at most $N$ parts, may be written as a product of $N$ factors. In part I, \cite{OS1}, we studied the behavior of  coefficients in the partial fraction decomposition of this product as $N \to \infty$ by applying the saddle-point method to get the asymptotics of the main terms. In this second part we bound the error terms. This involves estimating products of sines and further  saddle-point arguments. The saddle-points needed are associated to zeros of the analytically continued dilogarithm.
\end{abstract}

\def\s#1#2{\langle \,#1 , #2 \,\rangle}

\def\H{{\mathbf{H}}}
\def\F{{\frak F}}
\def\C{{\mathbb C}}
\def\R{{\mathbb R}}
\def\Z{{\mathbb Z}}
\def\Q{{\mathbb Q}}
\def\N{{\mathbb N}}
\def\G{{\Gamma}}
\def\GH{{\G \backslash \H}}
\def\g{{\gamma}}
\def\L{{\Lambda}}
\def\ee{{\varepsilon}}
\def\K{{\mathcal K}}
\def\Re{\mathrm{Re}}
\def\Im{\mathrm{Im}}
\def\PSL{\mathrm{PSL}}
\def\SL{\mathrm{SL}}
\def\Vol{\operatorname{Vol}}
\def\lqs{\leqslant}
\def\gqs{\geqslant}
\def\sgn{\operatorname{sgn}}
\def\res{\operatornamewithlimits{Res}}
\def\li{\operatorname{Li_2}}
\def\lis{\operatorname{Li_2^*}}
\def\clp{\operatorname{Cl}'_2}
\def\clpp{\operatorname{Cl}''_2}
\def\farey{\mathscr F}

\newcommand{\stira}[2]{{\left[ #1 \atop #2  \right]}}
\newcommand{\stirb}[2]{{\left\{ #1 \atop #2  \right\}}}
\newcommand{\norm}[1]{\left\lVert #1 \right\rVert}

\newcommand{\spr}[2]{\sideset{}{_{#2}^{-1}}{\textstyle \prod}({#1})}
\newcommand{\spn}[2]{\sideset{}{_{#2}}{\textstyle \prod}({#1})}
\newcommand{\spp}[2]{\sideset{}{_{#2}^{-2}}{\textstyle \prod}({#1})}

\newtheorem{theorem}{Theorem}[section]
\newtheorem{lemma}[theorem]{Lemma}
\newtheorem{prop}[theorem]{Proposition}
\newtheorem{conj}[theorem]{Conjecture}
\newtheorem{cor}[theorem]{Corollary}
\newtheorem{remark}[theorem]{Remark}

\newtheorem*{mainp}{Proposition 1.4}
\newtheorem*{maind}{Theorem 1.6}
\newtheorem*{maine}{Theorem 1.7}

\newcounter{coundef}
\newtheorem{adef}[coundef]{Definition}

\renewcommand{\labelenumi}{(\roman{enumi})}

\numberwithin{equation}{section}

\bibliographystyle{alpha}

\tableofcontents


\section{Introduction}
\subsection{Background}
The generating function for $p_N(n)$, the number of  partitions of $n$ into at most $N$ parts, and its partial fraction decomposition may be written as
\begin{equation}\label{tp}
\sum_{n=0}^\infty p_N(n) q^n =
\prod_{j=1}^N \frac{1}{1-q^j}=\sum_{\substack{0\leqslant h<k \leqslant N \\ (h,k)=1}}
\sum_{\ell=1}^{\lfloor N/k \rfloor} \frac{C_{hk\ell}(N)}{(q-e^{2\pi ih/k})^\ell}
\end{equation}
for coefficients $C_{hk\ell}(N)$ studied by Rademacher in \cite{Ra}. Each $C_{hk\ell}(N)$ is in the field $\Q(e^{2\pi i h/k})$ by \cite[Prop. 3.3]{OS}.  Let   $\li$ denote the dilogarithm.
It is shown in \cite[Sect. 1]{OS3} that
\begin{equation} \label{dilogw0}
\li(w)-2\pi i \log(w) = 0
\end{equation} has a unique solution, $w_0 \approx 0.916198 - 0.182459 i$,  and set $z_0:= 1+\log(1-w_0)/(2\pi i) \approx 1.18147 + 0.255528 i$.
With
\begin{equation*}\label{farey}
    \farey_N:=\Bigl\{ h/k \ : \ 1 \lqs k \lqs N, \ 0\lqs h < k,  \ (h,k)=1\Bigr\}
\end{equation*}
denoting the Farey fractions of order $N$ in $[0,1)$,  the asymptotic result
\begin{equation} \label{maineq1}
    \sum_{h/k \in \farey_{100}} C_{hk1}(N) = \Re\left[(-2  z_0 e^{-\pi i z_0})\frac{w_0^{-N}}{N^2}\right] +O\left( \frac{|w_0|^{-N}}{N^3}\right)
\end{equation}
is given in \cite[Thm. 1.2]{OS1}. This resolves an old conjecture of Rademacher in \cite[p. 302]{Ra} by showing that the limit of $C_{hk\ell}(N)$ as $N \to \infty$ does not exist in general since $|1/w_0|>1$, see \cite[Cor. 1.3]{OS1}.

Equation \eqref{maineq1} is a special case of the  more general theorem, \cite[Thm. 1.4]{OS1}, which we state next. Note that $C_{01\ell}(N)$ is the coefficient of $1/(q-1)^\ell$ in \eqref{tp}.

\begin{theorem}\label{mainthmb} There are explicit coefficients $c_{\ell,0},$ $c_{\ell,1}, \dots $ so that
\begin{multline} \label{maineq2}
   C_{01\ell}(N)+ \sum_{0<h/k \in \farey_{100}} \sum_{j=1}^\ell (e^{2\pi i h/k} -1)^{\ell-j} C_{hkj}(N) \\
   = \Re\left[\frac{w_0^{-N}}{N^{\ell+1}} \left( c_{\ell,0}+\frac{c_{\ell,1}}{N}+ \dots +\frac{c_{\ell,m-1}}{N^{m-1}}\right)\right] + O\left(\frac{|w_0|^{-N}}{N^{\ell+m+1}}\right)
\end{multline}
where $c_{\ell,0}=-2  z_0 e^{-\pi i z_0} (2\pi i z_0)^{\ell -1}$ and the implied constant depends only on  $\ell$ and $m$.
\end{theorem}

The main term of Theorem \ref{mainthmb} is shown in \cite{OS1}. The proof that  the size of the error term above is  $O\left(|w_0|^{-N}/N^{\ell+m+1}\right)$  is  sketched in \cite{OS1}, due to its length, and the detailed proof of this error bound is the main result of this paper.

Rademacher's coefficients $C_{hk\ell}(N)$ are fascinating numbers and their properties have been coming into focus with the recent papers \cite{DG,An,Mu,SZ,DrGe,OS}. Andrews gave the first formulas for them in \cite[Thm. 1]{An}. Further expressions were given in \cite{OS} with, for example, the relatively simple
\begin{equation*}
    C_{01\ell}(N)  = \frac{(-1)^N (\ell-1)!}{ N!} \sum_{j_0+j_1+j_2+ \cdots + j_N = N-\ell}
     \stirb{\ell+j_0}{\ell}\frac{B_{j_1}B_{j_2}  \cdots  B_{j_N}}{(\ell-1+j_0)!}\frac{ 1^{j_1} 2^{j_2} \cdots
N^{j_N}}{j_1 ! j_2 ! \cdots j_N!}
\end{equation*}
where $B_n$ is the $n$th Bernoulli number and $\stirb{n}{m}$ is the Stirling number,  denoting the number of ways to partition a set of size $n$ into $m$  non-empty subsets. Also, with $s_m(N):=1^m+2^m+ \cdots + N^m$,
\begin{multline*}
    C_{01\ell}(N) = \frac{(-1)^N }{ N!} \sum_{j_0+1j_1+2j_2+ \cdots + N j_{N} = N-\ell}
     \frac{1}{j_0! j_1! j_2! \cdots j_N!} \\
    \times
     \left(\frac{B_1}{1 \cdot 1!}\bigl(s_1(N)+1-\ell\bigr) \right)^{j_1}  \cdots \left(\frac{(-1)^{N-1}B_{N}}{N \cdot N!}\bigl(s_{N}(N)+1-\ell \bigr) \right)^{j_{N}}.
\end{multline*}
These results are \cite[Eq. (2.12), Prop 2.4]{OS} and in that paper the close connection is described between Rademacher's coefficients $C_{hk\ell}(N)$ and Sylvester's waves. In forthcoming work we develop this link  and obtain the asymptotics of the individual waves in Sylvester's  decomposition of the (unrestricted) partition function $p(n)$.

It is also shown in \cite[Thm. 7.3]{OS} that, for $r \gqs 1$,
\begin{equation*}
P_{01r}(N) :=(-1)^N N!  \cdot (-4)^r r! \cdot C_{01(N-r)}(N)
\end{equation*}
 is a monic polynomial in $N$ of degree $2r$ with $0$ and $1$ as roots. This proved part of Conjecture 7.1 in \cite{SZ}. In the remaining part, Sills and Zeilberger  conjecture that $P_{01r}(N)$ is convex and has coefficients that alternate in sign.

Rademacher realized, already in the 1937 paper \cite{Ra2}, that his celebrated formula for $p(n)$ leads to a decomposition similar to \eqref{tp}:
\begin{equation}\label{37}
\sum_{n=0}^\infty p(n) q^n =
\prod_{j=1}^\infty \frac{1}{1-q^j}=\sum_{\substack{0\leqslant h<k  \\ (h,k)=1}}
\sum_{\ell=1}^{\infty} \frac{C_{hk\ell}(\infty)}{(q-e^{2\pi ih/k})^\ell} \qquad(|q|<1),
\end{equation}
with numbers $C_{hk\ell}(\infty)$  computed explicitly in \cite[Eq. (130.6)]{Ra}. Using limited numerical evidence he conjectured that $\lim_{N\to \infty}C_{hk\ell}(N) = C_{hk\ell}(\infty)$. Numerical computations were extended in \cite{An,DG,SZ} with the results in \cite{SZ} indicating clearly
that Rademacher's conjecture was almost certainly false. Confirmation of this was given independently in \cite{DrGe} and  \cite{OS1}. The work of Drmota and Gerhold in \cite{DrGe} gives the main term in the asymptotics of $C_{01\ell}(N)$ as $N \to \infty$ using techniques involving the Mellin transform. 
The proof of our Theorem \ref{mainthmb}, in \cite{OS1} and this paper,  is  based on a different, conceptually simple idea that is described in the next subsection. Though certainly very long when all details are included, our proof results in the complete asymptotic expansion of a finite average containing $C_{01\ell}(N)$. With further improvements it should be possible to replace the average on the left side of \eqref{maineq2} with just $C_{01\ell}(N)$, see \cite[Conj. 1.5]{OS1}.

We highlight two further interesting directions for investigation leading from this paper.
\begin{enumerate}
\item It should be possible to obtain the asymptotics for all coefficients $C_{hk\ell}(N)$ with $k$ small. Based on Theorem \ref{c2sed} below, the asymptotic expansion of $C_{121}(N)$ was conjectured in \cite[Conj. 6.3]{OS} and \cite[Conj. 6.3]{OS1}. Elements possibly leading to the asymptotic expansion of $C_{131}(N)+C_{231}(N)$ are given in  \cite[Eq. (6.12)]{OS1}.
\item Rademacher's original conjecture on the relationship between the sequence $C_{hk\ell}(1),C_{hk\ell}(2), \dots$ and $C_{hk\ell}(\infty)$ was too simplistic. However, it seems clear that there is indeed a close relationship between them, as shown in \cite[Sect. 4]{SZ} and \cite[Table 2]{OS}. The precise nature of this link remains to be found.
\end{enumerate}

\subsection{Proof of Theorem \ref{mainthmb}}
We  introduce some notation and results from \cite[Sect. 1.3]{OS1} to describe the proof of Theorem \ref{mainthmb}.
Define the numbers
\begin{equation}\label{qhkl}
    Q_{hk\sigma}(N)   :=    2\pi i \operatornamewithlimits{Res}_{z=h/k} \frac{e^{2\pi i \sigma z}}{(1-e^{2\pi i z})(1-e^{2\pi i2 z}) \cdots (1-e^{2\pi i N z})}.
\end{equation}
The Rademacher coefficients $ C_{hk\ell }(N)$ are related to them by
\begin{equation} \label{exprc}
    C_{hk\ell }(N) =   \sum_{\sigma=1}^\ell \binom{\ell -1}{\sigma-1} (-e^{2\pi i h/k})^{\ell-\sigma} Q_{hk\sigma}(N)
\end{equation}
and for $\sigma$  a positive integer they satisfy
\begin{equation}
\sum_{h/k \in \farey_N} Q_{hk\sigma}(N) =  0  \label{ch}
\end{equation}
for $N(N+1)/2 > \sigma$. Put
\begin{equation} \label{a(n)}
    \mathcal A(N)  := \Bigl\{ h/k \ : \ N/2  <k \lqs N,  \ h=1 \text{ \ or \ } h=k-1 \Bigr\} \subseteq \farey_N
\end{equation}
and
decompose \eqref{ch} into
\begin{equation} \label{qqq0}
    \sum_{h/k \in \farey_{100}} Q_{hk\sigma}(N) + \sum_{h/k \in \farey_N- (\farey_{100} \cup \mathcal A(N))} Q_{hk\sigma}(N) + \sum_{h/k \in \mathcal A(N)} Q_{hk\sigma}(N)=  0.
\end{equation}
Theorem \ref{mainthmb} breaks into two natural parts. The first is proved in \cite{OS1}:
\begin{theorem}\label{maina} With $b_{0}=2  z_0 e^{-\pi i z_0}$ and  explicit  $b_{1}(\sigma),$ $b_{2}(\sigma), \dots $ depending on $\sigma \in \Z$ we have 
\begin{equation*}
   \sum_{h/k \in \mathcal A(N)} Q_{hk\sigma}(N) = \Re\left[\frac{w_0^{-N}}{N^{2}} \left( b_{0}+\frac{b_{1}(\sigma)}{N}+ \dots +\frac{b_{m-1}(\sigma)}{N^{m-1}}\right)\right] + O\left(\frac{|w_0|^{-N}}{N^{m+2}}\right)
\end{equation*}
for an implied constant depending only on  $\sigma$ and $m$.
\end{theorem}
The proof of the second part is  sketched in \cite{OS1}:
\begin{theorem}\label{mainb} There exists  $W<U:=-\log |w_0| \approx 0.068076$ so that
\begin{equation*}
   \sum_{h/k \in \farey_N-(\farey_{100} \cup \mathcal A(N))} Q_{hk\sigma}(N) = O\left(e^{WN}\right)
\end{equation*}
for an implied constant depending only on  $\sigma$. We may take $W=0.055$.
\end{theorem}
Theorem \ref{mainthmb} follows from combining Theorems \ref{maina} and \ref{mainb} with \eqref{qqq0} and \eqref{exprc}. This is done in \cite[Sect. 5.4]{OS1}.

\subsection{Main Results}
In this paper we give the details of the proof of Theorem \ref{mainb}.  This therefore completes the proof of Theorem \ref{mainthmb} and \eqref{maineq1}. The work in this paper and \cite{OS1} will also be useful in describing the asymptotics of Sylvester waves and restricted partitions; this corresponds to estimating $Q_{hk\sigma}(N)$ for $\sigma <0$ as discussed in \cite[Sect. 6.2]{OS1}. Further natural extensions and possible generalizations of our results are  given there as well.

Define the sine product
\begin{equation} \label{sidef}
    \spn{\theta}{m} :=\prod_{j=1}^m 2\sin (\pi j \theta)
\end{equation}
with $\spn{\theta}{0}:=1$. In Section \ref{sec-bndq} we show
\begin{prop} \label{prc}
For $2 \lqs k \lqs N$, $\sigma \in \R$ and $s := \lfloor N/k \rfloor$
\begin{equation*}
    |Q_{hk\sigma}(N)|  \lqs \frac{3}{k^3} \exp\left( N \frac{2  +  \log \left(1+ 3k/4 \right)}{k} +\frac{|\sigma|}{N} \right)\left|\spr{h/k}{N-s k} \right|.
\end{equation*}
\end{prop}
In Section \ref{sec-sinprod} we find  sharp general bounds for $\spr{h/k}{m}$. This requires the interesting sum
\begin{equation} \label{smvx}
    S(m;h,k):=\sum_{(\beta,\gamma)\in Z(h,k)} \frac{\sin(2\pi m \gamma/k)}{|\beta \gamma|}
\end{equation}
for
\begin{equation} \label{Z(h,k)}
Z(h,k):=\Bigl\{ (\beta,\gamma)\in \Z\times \Z \ : \ 1\lqs |\beta| < k, \ 1\lqs \gamma < k,  \ \beta h \equiv \gamma \bmod k  \Bigr\}.
\end{equation}
We will see that $\spr{h/k}{m}$ and $S(m;h,k)$ may be bounded in terms of $1/|\beta_0\gamma_0|$ where $(\beta_0,\gamma_0)$ is a pair in $Z(h,k)$ with $|\beta_0\gamma_0|$ minimal.

Combining a refinement of Proposition \ref{prc} with our bound for $\spr{h/k}{m}$ allows us to prove Theorem \ref{mainb} except for $h/k$ in the following sets
\begin{align}
    \mathcal C(N) & := \Bigl\{ h/k \ : \ \frac{N}{2}  <k \lqs N,  \ k \text{ odd}, \ h=2 \text{ \ or \ } h=k-2 \Bigr\}, \label{cnsub}\\
    \mathcal D(N) & := \Bigl\{ h/k \ : \ \frac{N}{2}  <k \lqs N,  \ k \text{ odd}, \ h=\frac{k-1}2 \text{ \ or \ } h=\frac{k+1}2 \Bigr\}, \label{dnsub}\\
    \mathcal E(N) & := \Bigl\{ h/k \ : \ \frac{N}{3}  <k \lqs \frac{N}{2},  \ h=1 \text{ \ or \ } h=k-1 \Bigr\}. \label{ensub}
\end{align}

For the next results we need a brief description of the zeros of the dilogarithm; see \cite[Sect. 2.3]{OS1} and \cite{OS3} for a fuller discussion. Initially defined as
\begin{equation}\label{def0}
\li(z):=\sum_{n=1}^\infty \frac{z^n}{n^2} \quad \text{ for }|z|\lqs 1,
\end{equation}
the dilogarithm has an analytic continuation given by
$
 -\int_{C(z)} \log(1-u) \frac{du}{u}
$
where the contour of integration $C(z)$ is a path from $0$ to $z \in \C$. This makes the dilogarithm a multi-valued holomorphic function with branch points at $0,$ $1$ and $\infty$.  See for example \cite{max}, \cite{Zag07}.
 We let $\li(z)$ denote the dilogarithm on its principal branch so that $\li(z)$ is a single-valued  holomorphic function on $\C-[1,\infty)$. It can be shown that the value of the analytically continued dilogarithm is always given by
\begin{equation}\label{dilogcont}
\li(z) + 4\pi^2  A +   2\pi i  B  \log \left(z\right)
\end{equation}
for some $A$, $B \in \Z$.

Let $w(A,B)$ be a zero of \eqref{dilogcont}. It is shown in \cite[Thm. 1.1]{OS3} that for $B\neq 0$, a zero $w(A,B)$ exists if and only if $-|B|/2<A\lqs |B|/2$ and is unique in this case. Each zero may be found to arbitrary precision using Newton's method according to \cite[Thm. 1.3]{OS3}. We already met $w_0 = w(0,-1)$ and we also need the two further zeros $w(1,-3) \approx -0.459473 - 0.848535i$, $w(0,-2) \approx 0.968482 - 0.109531i$ and the associated saddle-points
\begin{equation*}
    z_3:=3+\log \bigl(1-w(1,-3)\bigr)/(2\pi i), \qquad z_1:=2+\log \bigl(1-w(0,-2)\bigr)/(2\pi i).
\end{equation*}

\begin{theorem} \label{thmc} With $c_{0}^*=-z_3 e^{-\pi i z_3}/4$ and  explicit  $c_{1}^*(\sigma),$ $c_{2}^*(\sigma), \dots $ depending on $\sigma \in \Z$ we have
\begin{equation} \label{pressx}
   \sum_{h/k \in \mathcal C(N)} Q_{hk\sigma}(N) = \Re\left[\frac{w(1,-3)^{-N}}{N^{2}} \left( c_{0}^*+\frac{c_{1}^*(\sigma)}{N}+ \dots +\frac{c_{m-1}^*(\sigma)}{N^{m-1}}\right)\right] + O\left(\frac{|w(1,-3)|^{-N}}{N^{m+2}}\right)
\end{equation}
for an implied constant depending only on  $\sigma$ and $m$.
\end{theorem}

\begin{theorem} \label{c2sed} Let $\overline{N}$ denote $N \bmod 2$. With
\begin{equation}\label{d0cx}
    d_{0}\bigl(\overline{N}\bigr) = z_0 \sqrt{2  e^{-\pi i z_0}\bigl(e^{-\pi i z_0}+(-1)^N \bigr)}
\end{equation}
 and  explicit  $d_{1}\bigl(\sigma, \overline{N}\bigr),$ $d_{2}\bigl(\sigma, \overline{N}\bigr), \dots $ depending on $\sigma \in \Z$ and $\overline{N}$, we have
\begin{equation} \label{presdx}
   \sum_{h/k \in \mathcal D(N)} Q_{hk\sigma}(N) = \Re\left[\frac{w_0^{-N/2}}{N^{2}} \left( d_{0}\bigl(\overline{N}\bigr) +\frac{d_{1}\bigl(\sigma, \overline{N}\bigr)}{N}+ \dots +\frac{d_{m-1}\bigl(\sigma, \overline{N}\bigr)}{N^{m-1}}\right)\right] + O\left(\frac{|w_0|^{-N/2}}{N^{m+2}}\right)
\end{equation}
for an implied constant depending only on  $\sigma$ and $m$.
\end{theorem}
(By $w_0^{-N/2 }$ we mean $\left(\sqrt{w_0} \right)^{-N}$ where $\sqrt{w_0}$ is chosen as usual with $\Re(\sqrt{w_0})>0$.)

\begin{theorem} \label{thme} With $e_{0}=-3z_1 e^{-\pi i z_1}/2$ and  explicit  $e_{1}(\sigma),$ $e_{2}(\sigma), \dots $ depending on $\sigma \in \Z$ we have
\begin{equation} \label{presex}
   \sum_{h/k \in \mathcal E(N)} Q_{hk\sigma}(N) = \Re\left[\frac{w(0,-2)^{-N}}{N^{2}} \left( e_{0}+\frac{e_{1}(\sigma)}{N}+ \dots +\frac{e_{m-1}(\sigma)}{N^{m-1}}\right)\right] + O\left(\frac{|w(0,-2)|^{-N}}{N^{m+2}}\right)
\end{equation}
for an implied constant depending only on  $\sigma$ and $m$.
\end{theorem}

The above three estimates are the final elements required for Theorem \ref{mainb}, and its proof is given near the end of Section \ref{sec-e1}.
Theorems \ref{thmc},  \ref{c2sed} and \ref{thme} above are proved using the techniques developed in \cite{OS1} for Theorem \ref{maina}, though they each present new challenges. These techniques use the saddle-point method described in the next subsection.

In fact, Theorems \ref{thmc},  \ref{c2sed} and \ref{thme} are more than is needed for Theorem \ref{mainb}, but we included them for two reasons. First, they allow us to check our work, see Tables \ref{c2n1} -- \ref{e1n1}. Secondly, their asymptotic expansions point the way to further results and a better understanding of relations in the left side of the identity \eqref{ch}. Examples of these relations, from \cite[Sect. 6.2]{OS1}, are
\begin{align}
      Q_{011}(N) \quad & \sim \quad  -\sum_{h/k \in \mathcal A(N)} Q_{hk1}(N), \label{simsum1} \\
      Q_{121}(N) \quad & \sim \quad  -\sum_{h/k \in \mathcal D(N)} Q_{hk1}(N) \label{simsum2}
\end{align}
where by \eqref{simsum1} and \eqref{simsum2} (and  \eqref{simsum3}) we mean that, at least numerically, the asymptotic expansions of both sides seem to be  identical.
With Theorems \ref{thmc} and \ref{thme} we discover another  asymptotic relation. To describe it, let $\mathcal C'(N)$ be all $h/k \in \mathcal C(N)$   with $2N/3  <k \lqs N$, so that $\mathcal C'(N)$ is about two thirds of $\mathcal C(N)$ . Then
\begin{equation} \label{simsum3}
    3   \sum_{h/k \in \mathcal C'(N)} Q_{hk\sigma}(N) \quad \sim \quad \sum_{h/k \in \mathcal E(N)} Q_{hk\sigma}(N).
\end{equation}
 See the end of Section \ref{sec-e1} for more about \eqref{simsum3}.

\subsection{The saddle-point method}
The next result was used in \cite[Sect. 5.1]{OS1} and is a simpler version of \cite[Theorem 7.1, p. 127]{Ol}.
\begin{theorem}[Saddle-point method] \label{sdle}
Let $\mathcal P$ be a finite  polygonal path in $\C$ with $p(z)$, $q(z)$  holomorphic functions in a neighborhood of $\mathcal P$. Assume $p$, $q$ and $\mathcal P$ are independent of a parameter $N>0$. Suppose $p'(z)$ has a simple zero at a non-corner point $z_0 \in \mathcal P$ with  $\Re(p(z)-p(z_0))>0$ for $z\in \mathcal P$ except at $z=z_0$. Then there exist explicit numbers $a_{2s}$ depending on $p$, $q$, $z_0$ and $\mathcal P$ so that we have
\begin{equation}\label{sad}
\int_{\mathcal P} e^{-N \cdot p(z)}q(z)\, dz = 2e^{-N \cdot p(z_0)}\left(\sum_{s=0}^{S-1} \G(s+1/2)\frac{a_{2s}}{N^{s+1/2}} + O\left(\frac{1}{N^{S+1/2}}\right)\right)
\end{equation}
for $S$ an arbitrary positive integer and an implied constant independent of $N$.
\end{theorem}

Write the power series for $p$ and $q$ near $z_0$ as
\begin{align}
    p(z) & = p(z_0)+ p_0(z-z_0)^2+p_1(z-z_0)^3+ \cdots, \label{psp}\\
    q(z) & = q_0+q_1(z-z_0)+q_2(z-z_0)^2+ \cdots. \label{psq}
\end{align}
Choose $\omega \in \C$ giving the direction of the path $\mathcal P$ through $z_0$: near $z_0$,  $\mathcal P$  looks like $z=z_0+\omega t$ for small $t\in \R$ increasing.
Wojdylo in \cite[Theorem 1.1]{Woj} found  an explicit formula for the numbers $a_{2s}$:
\begin{equation} \label{a2s}
    a_{2s}= \frac{\omega}{2(\omega^2 p_0)^{1/2}} \sum_{i=0}^{2s} q_{2s-i} \sum_{j=0}^i p_0^{-s-j} \binom{-s-1/2}{j} \hat{B}_{i,j}(p_1, p_2, \dots)
\end{equation}
where we must choose the square root $(\omega^2 p_0)^{1/2}$ in \eqref{a2s} so that $\Re \bigl((\omega^2 p_0)^{1/2}\bigr)>0$ and $\hat{B}_{i,j}$ is the {\em partial ordinary Bell polynomial}.
The first cases are
\begin{equation} \label{a2sb}
    a_0= \frac{\omega}{2(\omega^2 p_0)^{1/2}} q_0, \qquad a_2 = \frac{\omega}{2(\omega^2 p_0)^{1/2}}\left(
    \frac{q_2}{p_0} - \frac{3}{2} \frac{p_1 q_1 + p_2 q_0}{p_0^2} + \frac{15}{8} \frac{p_1^2 q_0}{p_0^3}\right),
\end{equation}
agreeing with \cite[p. 127]{Ol}.

We will be applying Theorem \ref{sdle} to functions $p$ of the form
\begin{equation}\label{pdfn}
    p_d(z):=\frac{ - \li\left(e^{2\pi i z}\right) +\li(1) +4\pi^2 d}{2\pi i z}.
\end{equation}
 Recall that $\li(z)$ is holomorphic on $\C - [1,\infty)$. Hence $p_d(z)$ is
a single-valued holomorphic function away from the vertical branch cuts $(-i\infty,n]$ for $n \in \Z$. (We use $(-i\infty,n]$ to indicate all points in $\C$ with real part $n$ and imaginary part at most $0$.) The next result is shown in \cite[Sect. 2.3]{OS1}. The notation $w(A,B)$ for the dilogarithm zeros is defined after \eqref{dilogcont}.

\begin{theorem} \label{disol}
Fix integers $m$ and $d$ with $-|m|/2<d\lqs |m|/2$. Then there is a unique solution to $p_d'(z)=0$ for $z \in \C$ with $m-1/2<\Re(z)<m+1/2$ and $z \not\in (-i\infty,m]$. Denoting this saddle-point by $z^*$, it is given by
\begin{equation}\label{uniq}
    z^*=m+\frac{\log \bigl(1-w(d,-m)\bigr)}{2\pi i}
\end{equation}
and satisfies
\begin{equation} \label{pzlogw}
    p_d(z^*)=\log \bigl(w(d,-m)\bigr).
\end{equation}
\end{theorem}

\section{The maxima and minima of $\spn{h/k}{m}$} \label{sec-sinprod}
Recall the set $Z(h,k)$ from \eqref{Z(h,k)}.
We will also need Clausen's integral,
\begin{align}\label{simo}
\cl(\theta) & :=-\int_0^\theta \log |2\sin( x/2) | \, dx \qquad (\theta \in \R)\\
& \phantom{:}= \sum_{n=1}^\infty \frac{\sin(n\theta)}{n^2}. \label{clsd}
\end{align}
 The maximum value of $\cl(\theta)$ is
$\cl(\pi/3)  \approx 1.0149416$.

\begin{theorem} \label{mainest}
For all $m$, $h$, $k\in \Z$ with $1 \lqs  h < k$, $(h,k)=1$ and $0 \lqs m < k$ we have
\begin{equation}\label{logpx}
\frac{1}{k}  \log \left| \spr{h/k}{m} \right|  =  \frac{\cl(2\pi m \gamma_0 h/k)}{2\pi |\beta_0 \gamma_0|}  + O\left(\frac{\log k}{\sqrt{k}}\right)
\end{equation}
where $(\beta_0,\gamma_0)$ is a pair in $Z(h,k)$ with $|\beta_0 \gamma_0|$ minimal. The   implied constant in \eqref{logpx} is absolute and in fact this error is  bounded by $(16.05+\sqrt{2}/\pi \log k)/\sqrt{k}$.
\end{theorem}

We prove Theorem \ref{mainest} in the following subsections, assuming throughout  that $m$, $h$, $k$ satisfy its conditions.
Define $D(h,k)$ to be the above minimal value $|\beta_0 \gamma_0|$. For example, it is easy to see that
\begin{equation}\label{egd}
     D(h,k)  = 1 \iff  h \equiv  \pm 1 \bmod k
\end{equation}
and if $D(h,k) \neq 1$ then
\begin{equation}\label{egd2}
     D(h,k)  = 2 \iff   h \text{ or }  h^{-1} \equiv  \pm 2 \bmod k
\end{equation}
with $k$ necessarily odd.
Since $(1,h)\in Z(h,k)$ we have $D(h,k)\lqs h<k$. We will see later in Lemma \ref{l2} that there is a unique $(\beta_0,\gamma_0) \in Z(h,k)$ with $|\beta_0 \gamma_0|$ minimal if $|\beta_0 \gamma_0| < \sqrt{k/2}$.

The corollary we will need, Corollary \ref{needc},  says there exists an absolute constant $\tau$ such that
\begin{equation}\label{logpxx}
\frac{1}{k} \Bigl| \log \bigl| \spn{h/k}{m} \bigr| \Bigr| \lqs  \frac{\cl(\pi/3)}{2\pi D(h,k)}  + \tau \frac{\log k}{\sqrt{k}}.
\end{equation}
For example, Figure \ref{5efig} compares both sides of \eqref{logpxx} with $k=101$, $\tau=0$ and
\begin{equation}\label{phihk}
\Psi(h,k):=\max_{0\lqs m <k} \left\{ \frac{1}{k} \Bigl| \log \bigl| \spn{h/k}{m} \bigr| \Bigr|  \right\}.
\end{equation}



\SpecialCoor
\psset{griddots=5,subgriddiv=0,gridlabels=0pt}
\psset{xunit=0.09cm, yunit=40cm}
\psset{linewidth=1pt}
\psset{dotsize=2pt 0,dotstyle=*}

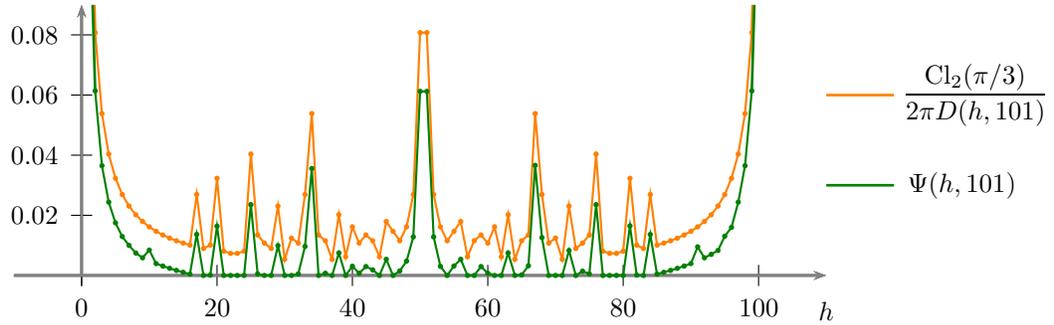
\begin{figure}[h]
\begin{center}
\begin{pspicture}(-15,-0.02)(120,0.09) 

\psaxes[linecolor=gray,Ox=0,Oy=0,Dx=20,dx=20,Dy=0.02,dy=0.02]{->}(0,0)(-10,0)(110,0.09)

\newrgbcolor{darkgrn}{0 0.5 0}
\newrgbcolor{lightgrn}{0 1 0}

\savedata{\mydata}[
{{2., 0.0807665}, {3., 0.0538443}, {4.,
  0.0403832}, {5., 0.0323066}, {6., 0.0269222}, {7., 0.0230761}, {8.,
  0.0201916}, {9., 0.0179481}, {10., 0.0161533}, {11.,
  0.0146848}, {12., 0.0134611}, {13., 0.0124256}, {14.,
  0.0115381}, {15., 0.0107689}, {16., 0.0100958}, {17.,
  0.0269222}, {18., 0.00897405}, {19., 0.0100958}, {20.,
  0.0323066}, {21., 0.00807665}, {22., 0.00734241}, {23.,
  0.00734241}, {24., 0.00807665}, {25., 0.0403832}, {26.,
  0.0134611}, {27., 0.0107689}, {28., 0.00897405}, {29.,
  0.0230761}, {30., 0.00538443}, {31., 0.0124256}, {32.,
  0.0107689}, {33., 0.0269222}, {34., 0.0538443}, {35.,
  0.0134611}, {36., 0.0115381}, {37., 0.00538443}, {38.,
  0.0201916}, {39., 0.00621281}, {40., 0.0161533}, {41.,
  0.0107689}, {42., 0.0134611}, {43., 0.0115381}, {44.,
  0.00621281}, {45., 0.0179481}, {46., 0.0146848}, {47.,
  0.0115381}, {48., 0.0161533}, {49., 0.0269222}, {50.,
  0.0807665}, {51., 0.0807665}, {52., 0.0269222}, {53.,
  0.0161533}, {54., 0.0115381}, {55., 0.0146848}, {56.,
  0.0179481}, {57., 0.00621281}, {58., 0.0115381}, {59.,
  0.0134611}, {60., 0.0107689}, {61., 0.0161533}, {62.,
  0.00621281}, {63., 0.0201916}, {64., 0.00538443}, {65.,
  0.0115381}, {66., 0.0134611}, {67., 0.0538443}, {68.,
  0.0269222}, {69., 0.0107689}, {70., 0.0124256}, {71.,
  0.00538443}, {72., 0.0230761}, {73., 0.00897405}, {74.,
  0.0107689}, {75., 0.0134611}, {76., 0.0403832}, {77.,
  0.00807665}, {78., 0.00734241}, {79., 0.00734241}, {80.,
  0.00807665}, {81., 0.0323066}, {82., 0.0100958}, {83.,
  0.00897405}, {84., 0.0269222}, {85., 0.0100958}, {86.,
  0.0107689}, {87., 0.0115381}, {88., 0.0124256}, {89.,
  0.0134611}, {90., 0.0146848}, {91., 0.0161533}, {92.,
  0.0179481}, {93., 0.0201916}, {94., 0.0230761}, {95.,
  0.0269222}, {96., 0.0323066}, {97., 0.0403832}, {98.,
  0.0538443}, {99., 0.0807665}}
  ]
\dataplot[linecolor=orange,linewidth=0.8pt,plotstyle=line]{\mydata}
\dataplot[linecolor=orange,linewidth=0.8pt,plotstyle=dots]{\mydata}

\savedata{\mydatab}[
{{2., 0.061391}, {3., 0.0364941}, {4.,
  0.0243923}, {5., 0.0175091}, {6., 0.0129402}, {7., 0.00999248}, {8.,
   0.00746955}, {9., 0.00587242}, {10., 0.00843565}, {11.,
  0.00394987}, {12., 0.00312541}, {13., 0.00237323}, {14.,
  0.00168311}, {15., 0.00104695}, {16., 0.000458233}, {17.,
  0.0136426}, {18., 0.}, {19., 0.000103942}, {20., 0.0163961}, {21.,
  0.}, {22., 0.}, {23., 0.}, {24., 0.}, {25., 0.0235572}, {26.,
  0.000509879}, {27., 0.}, {28., 0.}, {29., 0.00991222}, {30.,
  0.}, {31., 0.}, {32., 0.000515158}, {33., 0.00965916}, {34.,
  0.0355716}, {35., 0.}, {36., 0.00073106}, {37., 0.}, {38.,
  0.00748098}, {39., 0.}, {40., 0.00306658}, {41., 0.00072161}, {42.,
  0.00296196}, {43., 0.00184661}, {44., 0.}, {45., 0.00537417}, {46.,
  0.}, {47., 0.00150532}, {48., 0.0047695}, {49., 0.0128317}, {50.,
  0.0612263}, {51., 0.0612263}, {52., 0.0128317}, {53.,
  0.00303832}, {54., 0.}, {55., 0.00314456}, {56., 0.00537417}, {57.,
  0.}, {58., 0.}, {59., 0.00296196}, {60., 0.00072161}, {61.,
  0.}, {62., 0.}, {63., 0.00748098}, {64., 0.}, {65.,
  0.000184613}, {66., 0.00321521}, {67., 0.0366182}, {68.,
  0.0125759}, {69., 0.}, {70., 0.}, {71., 0.}, {72.,
  0.00830841}, {73., 0.}, {74., 0.00143161}, {75., 0.000509879}, {76.,
   0.0235572}, {77., 0.}, {78., 0.}, {79., 0.}, {80., 0.}, {81.,
  0.0164845}, {82., 0.000103942}, {83., 0.}, {84., 0.0136426}, {85.,
  0.000458233}, {86., 0.00104695}, {87., 0.00168311}, {88.,
  0.00237323}, {89., 0.00312541}, {90., 0.00394987}, {91.,
  0.00950783}, {92., 0.00587242}, {93., 0.00701132}, {94.,
  0.00830938}, {95., 0.0130541}, {96., 0.0159929}, {97.,
  0.0243923}, {98., 0.0364941}, {99., 0.061391}}
  ]
\dataplot[linecolor=darkgrn,linewidth=0.8pt,plotstyle=line]{\mydatab}
\dataplot[linecolor=darkgrn,linewidth=0.8pt,plotstyle=dots]{\mydatab}

\psline[linecolor=orange](2., 0.0807665)(1.89,0.09)
\psline[linecolor=darkgrn](2., 0.061391)(1.63,0.09)
\psline[linecolor=orange](99., 0.0807665)(99.11,0.09)
\psline[linecolor=darkgrn](99., 0.061391)(99.37,0.09)

\rput(110,-0.012){$h$}
\rput(132,0.06){$\displaystyle \frac{\cl(\pi/3)}{2\pi D(h,101)}$}
\rput(130,0.03){$\Psi(h,101)$}
\psline[linecolor=orange](110,0.06)(120,0.06)
\psline[linecolor=darkgrn](110,0.03)(120,0.03)


\end{pspicture}
\caption{Bounding $\Psi(h,k)$ for $1 \leqslant h \lqs k-1$ and $k=101$}\label{5efig}
\end{center}
\end{figure}


\subsection{Relating $\spr{h/k}{m}$ to $S(m;h,k)$}

 By \eqref{simo} we have $\clp(\theta)=-\log |2\sin( \theta/2) |$ and
\begin{equation}\label{cata}
\log \left| \spr{h/k}{m} \right| = \sum_{j=1}^m \clp(2\pi j h/k).
\end{equation}
With the sum $S(m;h,k)$ defined in \eqref{smvx}, our first goal is to prove:
\begin{prop} \label{sap}
For $0\lqs m <k$ and an absolute implied constant
\begin{equation*}
\sum_{j=1}^m \clp(2\pi j h/k)= \frac{k}{2\pi} S(m;h,k) +O\left(\log^2 k\right).
\end{equation*}
\end{prop}

Let
\begin{equation}\label{fdef}
f_L(x):=\sum_{n=1}^L \frac{\cos(nx)}{n}
\end{equation}
and define  $\norm{x}$ as the distance from $x\in \R$ to the nearest integer, so that $0\lqs \norm{x} \lqs 1/2$.
\begin{lemma} \label{lema}
For $L\gqs 1$ and $x\in \R$, $x \not\in \Z$ we have
$$
\clp(2\pi x)=f_L(2\pi x)+O\left(\frac 1{L\norm{x}}\right).
$$
\end{lemma}
\begin{proof} We first claim that
\begin{equation}\label{vir}
\left| \sum_{r=L}^M \frac{\cos(2\pi  r x)}{r}\right| \lqs \frac 1{L\norm{x}}
\end{equation}
for $x \not\in \Z$.
Let $A_m(2\pi x):=\sum_{r=1}^m e^{2\pi i r x}$. Then this geometric series evaluates to
$$
A_m(2\pi x)=-\frac{i}{2}\frac{e^{2\pi i(m+1/2)x}-e^{\pi i x}}{\sin \pi x}
$$
and the inequality $|\sin \pi x| \gqs 2\norm{x}$ implies $|A_m(2\pi x)| \lqs 1/(2\norm{x})$.
By partial summation
$$
\sum_{r=L}^M \frac{e^{2\pi i r x}}{r} = \frac{A_M}{M}-\frac{A_{L-1}}{L}+\sum_{d=L}^{M-1} \frac{A_d}{d(d+1)}.
$$
Taking real parts, using the bound for $A_m$ and evaluating the telescoping sum shows \eqref{vir}.

Now $\sum_{n=1}^L \sin(n x) n^{-2}$ as $L\to \infty$ converges uniformly to $\cl(x)$. The derivative of the above partial sum is $f_L(x)$. As
$L\to \infty$,  \eqref{vir} implies that $f_L(2\pi x)$ converges uniformly for $x$ in any closed interval not containing an integer. Hence, with \cite[Thm. 7.17]{Ru},
  $\lim_{L\to\infty}f_L(x) = \clp(2\pi x)$ for $x \not\in \Z$ and  the lemma follows.
\end{proof}

\begin{cor} \label{rome}
We have
$$
\sum_{j=1}^m \clp(2\pi j h/k)=\sum_{j=1}^m f_k(2\pi j h/k)+O\left(\log k\right).
$$
\end{cor}
\begin{proof}
Use
$$
\sum_{j=1}^m \frac{1}{\norm{jh/k}} \lqs \sum_{j=1}^{k-1} \frac{1}{\norm{jh/k}} \lqs 2 \sum_{j=1}^{k/2} \frac{1}{\norm{j/k}}
=2 \sum_{j=1}^{k/2} \frac{k}{j}.
$$
With $\sum_{j=1}^{k} 1/j \lqs 1+\log k$ we get
$$
\sum_{j=1}^m \frac{1}{k\norm{jh/k}} \ll \log k
$$
and the corollary now follows from Lemma \ref{lema}.
\end{proof}

\begin{lemma}\label{er}
For $0 \lqs m <k$ and $L=k^2$,
\begin{equation}\label{dsu}
\sum_{j=1}^m \clp(2\pi j h/k)= \frac{k}{2\pi} \sum_{l=-L}^L\sum_{n=1}^{k-1} \frac{\sin(2\pi m(nh + lk)/k)}{n(nh+lk)}+O\left(\log k\right).
\end{equation}
\end{lemma}
\begin{proof}
Apply  Euler-Maclaurin summation, in the form of \cite[Corollary 4.3]{IwKo}, to find
\begin{multline} \label{mul}
\sum_{j=1}^m f_k(2\pi j h/k) = \sum_{l=-L}^L \int_{0}^m f_k(2\pi x h/k) e^{2\pi i l x}\,dx \\
+ \frac 12 f_k(2\pi m h/k)- \frac 12 f_k(0)
+O\left( \int_0^m \frac{|f_k'(2\pi x h/k) 2\pi h/k|}{1+L\norm{x}}\, dx \right)
\end{multline}
where the implied constant is absolute.
Clearly we see $|f_k(x)|\lqs 1+\log k$ and $|f'_k(x)|\lqs  k$. To bound the error term in \eqref{mul} note that
$$
\int_0^m \frac{dx}{1+L\norm{x}} \lqs \int_0^{k-1} \frac{dx}{1+L\norm{x}} = 2k \int_0^{1/2} \frac{dx}{1+L x} =\frac{2k(1+\log L/2)}{L}.
$$
Hence, on choosing $L=k^2$, \eqref{mul} implies
\begin{equation}\label{fvs}
    \sum_{j=1}^m f_k(2\pi j h/k) = \sum_{l=-L}^L \int_{0}^m f_k(2\pi x h/k) e^{2\pi i l x}\,dx +O(\log k).
\end{equation}
Use $\cos \theta =(e^{i\theta}+e^{-i\theta})/2$ to evaluate the right side of \eqref{fvs} as follows.
\begin{multline*}
    \sum_{l=-L}^L \int_{0}^m f_k(2\pi x h/k) e^{2\pi i l x}\,dx
    = \sum_{l=-L}^L \sum_{n=1}^k \int_{0}^m \frac{\cos(2\pi n x h/k)}{n} e^{2\pi i l x}\,dx\\
    = \frac{1}{4\pi i}\sum_{l=-L}^L \sum_{n=1}^k \left( \frac{e^{2\pi i m(nh/k+l)}-1}{n(nh/k+l)}
    + \frac{e^{2\pi i m(-nh/k+l)}-1}{n(-nh/k+l)}\right)\\
    = \frac{1}{4\pi i}\sum_{l=-L}^L \sum_{n=1}^{k-1}  \frac{e^{2\pi i m(nh/k+l)}-e^{-2\pi i m(nh/k+l)}}{n(nh/k+l)}.
\end{multline*}
Combining this with Corollary  \ref{rome}
completes the proof.
\end{proof}

To simplify the right of \eqref{dsu} set
$$
H(d)=H(d,L;h,k):=\#\Big\{(l,n) \ :\ nh+lk=d, 1\lqs n \lqs k-1, -L \lqs l \lqs L \Big\}.
$$
Then the double sum  equals
\begin{equation}\label{sak}
\sum_{d \in \Z} H(d) \frac{\sin(2\pi m d/k)}{(d h^{-1} \bmod k)d}
\end{equation}
where we exclude $d$s that are multiples of $k$, since $H(d)$ is necessarily $0$ if $k | d$, and we understand here and throughout that $0 \lqs (* \bmod k) \lqs k-1$.

\begin{lemma} \label{aks}
Recall that $L=k^2$. For all $d \in \Z$ we have $H(d)=H(d,L;h,k)$ equalling $0$ or $1$. Also
\begin{eqnarray}
  H(d) &=& 1 \quad \text{ for } \quad 1\lqs |d|<k, \label{ak1}\\
  H(d) &=& 0 \quad \text{ for } \quad  |d|> 2k^3. \label{ak0}
\end{eqnarray}
\end{lemma}
\begin{proof}
Since $(h,k)=1$ there exist $n_0$, $l_0$ such that $n_0 h+l_0 k =1$. Then for all $t\in \Z$
$$(n_0+t k)h+(l_0-th)k =1
$$
and we may choose $n_0$, $l_0$ satisfying $1\lqs n_0<k$ and $-h<l_0\lqs -1$. Similarly, for fixed $h$, $k$, $d$, all solutions $(n,l)$ of $n h+ l k=d$ are given by
\begin{equation}\label{nk0}
n=d n_0+t k, \quad l=d l_0 -t h \quad (t\in \Z).
\end{equation}
Hence, for $k \nmid d$, there is exactly one solution $(n,l)$ with $1\lqs n \lqs k-1$. Then $H(d)=1$ if the corresponding $l$ satisfies $-L \lqs l \lqs L$ and $H(d)=0$ otherwise.

In \eqref{nk0}, if $1\lqs n \lqs k-1$ then $t=-\lfloor d n_0/k\rfloor$. Therefore
$$
l=d l_0 -t h = d l_0 +h\lfloor d n_0/k\rfloor
$$
and $l$ satisfies $-k^2 < l < k^2$ for $|d|<k$. This proves \eqref{ak1}.
Finally, to show \eqref{ak0}, note that $|n|<k$, $|l|\lqs L$ implies $|nh+lk|<k(h+L)<2k^3$.
\end{proof}

The sum \eqref{sak} with indices $d$ restricted to $|d|<k$ is
\begin{equation} \label{sumvqx}
\sum_{-k < d < k,\ d \neq 0} \frac{\sin(2\pi m d/k)}{(d h^{-1} \bmod k)d}.
\end{equation}
Replacing $d$ by $dh \bmod k$ if $d>0$, and   $d$ by $-(dh \bmod k) \equiv (-dh) \bmod k$ if $d<0$, allows us to  write \eqref{sumvqx} as
\begin{equation*}
\sum_{-k < d < k,\ d \neq 0} \frac{\sin(2\pi m d h/k)}{(d h \bmod k)|d|} = S(m;h,k).
\end{equation*}

\begin{proof}[\bf Proof of Proposition \ref{sap}]
With Lemmas \ref{er} and \ref{aks} we have demonstrated that
\begin{equation} \label{dems}
    \sum_{j=1}^m \clp(2\pi j h/k)= \frac{k}{2\pi} S(m;h,k) +\frac{k}{2\pi} \sum_{d \in \Z \ : \ k<|d|<2k^3 } H(d) \frac{\sin(2\pi m d/k)}{(d h^{-1} \bmod k)d}  + O\left(\log k\right).
\end{equation}
To estimate the  sum on the right of \eqref{dems}, write $d=uk+r$ and use Lemma \ref{aks} to see  that it is bounded by
\begin{equation}\label{ssm}
 \sum_{\substack{-2k^2\lqs u \lqs 2k^2 \\ u\neq 0,-1}}
\sum_{r=1}^{k-1} \frac{1}{|uk+r|(r h^{-1} \bmod k)}.
\end{equation}
For $u\gqs 1$ the inner sum is less than
$$
\sum_{r=1}^{k-1} \frac{1}{uk(r h^{-1} \bmod k)}=\frac{1}{uk}\sum_{r=1}^{k-1} \frac{1}{r} <\frac{1+\log k}{uk}.
$$
Similarly for $u\lqs -2$ and therefore \eqref{ssm} is bounded by
$$
2\frac{1+\log k}{k}\sum_{u=1}^{2k^2} \frac 1u \ll \frac{\log^2 k}{k}. \qedhere
$$
\end{proof}

\subsection{Relating  $S(m;h,k)$ to Clausen's integral}

With \eqref{cata} and Proposition \ref{sap} we have proved that
\begin{equation} \label{cqor}
    \frac{1}{k}\log \left| \spr{h/k}{m} \right| = \frac{S(m;h,k)}{2\pi}  +O\left(\frac{\log^2 k}{k}\right).
\end{equation}

\begin{remark}
{\rm The implied constant in \eqref{cqor} is absolute and we may find it explicitly. In Corollary \ref{rome} the error is bounded by $2(1+\log(k/2))$. In \eqref{mul} the implied constant can be $1/2+1/\pi$ which follows (see \cite[Eq. (4.18)]{IwKo}) from
\begin{equation}\label{ikex}
    \left| x-\lfloor x \rfloor -1/2+\sum_{j=1}^L \frac{\sin(2\pi j x)}{\pi j}\right| \lqs \frac{T}{1+L \norm{x}} \qquad (T=1/2+1/\pi).
\end{equation}
 To prove \eqref{ikex}, show that the left is bounded by $1/2$ and, with a similar proof to Lemma \ref{lema}, also bounded by $1/(\pi L \norm{x})$. This yields \eqref{ikex}. (It seems that $T=1/2$ should be possible.) Hence the error in Lemma \ref{er} is bounded by
$$
2(1+\log(k/2))+ 1 + \log k +4\pi T(1+\log(k^2/2)).
$$
For Proposition \ref{sap} we add $(1+\log k)(1+\log(2k^2))/\pi$. Altogether this shows the error in \eqref{cqor} is bounded by
\begin{equation} \label{khbd}
    (5.31 +24.75\log k + 2/\pi \log^2 k)/k < 40.18(\log^2 k)/k  \qquad (k\gqs 2).
\end{equation}
}
\end{remark}

For the proof of Theorem \ref{mainest} we therefore need to estimate $S(m;h,k)$ in \eqref{cqor}.
To do this, note that the largest terms in the sum \eqref{smvx} should occur when $|\beta|$ and $\gamma$ are both small.
We introduce a parameter $R$ to the set $Z(h,k)$ to control the size of the elements:
\begin{equation}\label{brhk}
Z_{R}(h,k):=\Bigl\{ (\beta,\gamma)\in \Z\times \Z \ : \ 1\lqs |\beta| < R, \ 1\lqs \gamma < R,  \ \beta h \equiv \gamma \bmod k  \Bigr\}.
\end{equation}
Then $Z(h,k)$ is $Z_k(h,k)$ in this notation.

\begin{lemma} \label{l1}
For an absolute implied constant
\begin{equation} \label{lsk}
\sum_{(\beta,\gamma)\in Z_{k}(h,k)-Z_{R}(h,k)} \frac{\sin(2\pi m \gamma/k)}{|\beta \gamma|}  =O\left(\frac{\log R}{R}\right).
\end{equation}
\end{lemma}
\begin{proof}
We may partition the terms of the sum on the left of \eqref{lsk} into the three cases where $|\beta|\gqs R$ or $\gamma \gqs R$ or both. The first two corresponding sums are each bounded by $2(1+\log R)/R$. With the Cauchy-Schwarz inequality, the third is bounded by
\begin{equation*}
    2\left( \sum_{\beta=R}^{k-1} \frac 1{\beta^2}\right)^{1/2}\left( \sum_{\gamma=R}^{k-1} \frac 1{\gamma^2}\right)^{1/2} < 2\left( \sum_{d=R}^{\infty} \frac 1{d^2}\right) < \frac{2}{R}\left(1+ \frac{1}{R}\right). \qedhere
\end{equation*}
\end{proof}


\begin{lemma} \label{l2}
 Suppose $Z_{R}(h,k)$ is non-empty and $k \gqs 2R^2$. Let $(\beta_1,\gamma_1)$ be a pair in $Z_{R}(h,k)$ with $|\beta_1 \gamma_1|$ minimal. Then for each $(\beta,\gamma) \in Z_{R}(h,k)$ there exists a positive integer $\lambda$ such that $(\beta,\gamma) = (\lambda\beta_1,\lambda\gamma_1)$.
\end{lemma}
\begin{proof}
The number $\beta$  may not have an inverse mod $k$ so write $\beta=\beta'k'$ with $k' | k$ and $\gcd(\beta',k)=1$. Necessarily we also have  $\gamma=\gamma'k'$ with  $\gcd(\gamma',k)=1$. Similarly, there exists $k_1 | k$ so that
$$
\beta_1 =\beta'_1 k_1 , \quad \gamma_1  =\gamma'_1 k_1 , \quad \gcd(\beta'_1,k)=\gcd(\gamma'_1,k)=1.
$$
Then
$$
h \equiv (\beta')^{-1} \gamma' \bmod k/k', \quad h \equiv (\beta'_1)^{-1} \gamma'_1 \bmod k/k_0
$$
and letting $k^*=\gcd(k/k',k/k_1)$ we obtain
$$
(\beta')^{-1} \gamma' \equiv (\beta'_1)^{-1} \gamma'_1 \bmod k^*
$$
so that
\begin{equation}\label{bga}
\beta'_1 \gamma' - \beta'\gamma'_1 \equiv 0 \bmod k^*.
\end{equation}
Now
\begin{equation}\label{bga2}
|\beta'_1 \gamma' - \beta'\gamma'_1|<\frac{2R^2}{k_1 k'} \lqs \frac{k}{k_1 k'} \lqs k^*
\end{equation}
so that \eqref{bga} and \eqref{bga2} imply
$$
\beta'_1 \gamma' - \beta'\gamma'_1 = 0
$$
which, in turn, shows that $\beta/\beta_1=\gamma/\gamma_1$. Hence  $(\beta,\gamma)=(\mu \beta_1, \mu \gamma_1)$ for $\mu:=\gamma/\gamma_1 \in \Q_{>0}$.
Write $\mu=\lambda+\delta$ with $\lambda \in \Z$ and $0 \lqs \delta <1$. If $0< \delta <1$
then
$$
(\beta,\gamma)-\lambda(\beta_1,  \gamma_1) = (\beta-\lambda \beta_1, \gamma -\lambda \gamma_1)=(\delta \beta_1, \delta\gamma_1)\in Z_{k}(h,k),
$$
but $|\delta^2 \beta_1 \gamma_1|< |\beta_1 \gamma_1|$ and $|\beta_1 \gamma_1|$ was supposed to be minimal. We must have $\delta=0$, as required.
\end{proof}

\begin{prop} \label{l3}
Let $(\beta_0,\gamma_0)$ be a pair in $Z_k(h,k)$ with $|\beta_0 \gamma_0|$ minimal, and so equalling $D(h,k)$. Then for an absolute implied constant
\begin{equation} \label{mo3}
    S(m;h,k)=\frac{\cl(2\pi m \gamma_0/k)}{|\beta_0 \gamma_0|}+O\left(\frac{\log k}{\sqrt{k}}\right).
\end{equation}
\end{prop}
\begin{proof}
By Lemma \ref{l1} with $R=\sqrt{k/2}$
\begin{equation}\label{smfx}
    S(m;h,k)=\sum_{(\beta,\gamma)\in Z_{\sqrt{k/2}}(h,k)} \frac{\sin(2\pi m \gamma/k)}{|\beta \gamma|}+O\left(\frac{\log k}{\sqrt{k}}\right).
\end{equation}
\begin{enumerate}
\item[{\bf Case $(i)$}] Assume first that $Z_{\sqrt{k/2}}(h,k)$ is empty. If $(\beta_0,\gamma_0)\notin Z_{\sqrt{k/2}}(h,k)$  it follows that $|\beta_0 \gamma_0| \gqs \sqrt{k/2}$ and so
\begin{equation} \label{mo2}
    \frac{\cl(2\pi m \gamma_0/k)}{|\beta_0 \gamma_0|} = O\left( \frac{1}{\sqrt{k}}\right).
\end{equation}
Then \eqref{mo3} follows from \eqref{smfx} and \eqref{mo2}.
\item[{\bf Case $(ii)$}]Assume now that $Z_{\sqrt{k/2}}(h,k)$ is not empty.
Apply Lemma \ref{l2} with the same $R=\sqrt{k/2}$, and $(\beta_1,\gamma_1) \in  Z_{\sqrt{k/2}}(h,k)$ with $|\beta_1 \gamma_1|$ minimal, to get
\begin{align}
    \sum_{(\beta,\gamma)\in Z_{\sqrt{k/2}}(h,k)} \frac{\sin(2\pi m \gamma/k)}{|\beta \gamma|}
    & =\frac{1}{|\beta_1 \gamma_1|}\sum_{1\lqs \lambda < \sqrt{k/2}/\max\{|\beta_1|,\gamma_1\}} \frac{\sin(2\pi m \lambda \gamma_1/k)}{\lambda^2} \notag\\
    & = \frac{\cl(2\pi m \gamma_1/k)}{|\beta_1 \gamma_1|} + O\left( \frac{1}{|\beta_1 \gamma_1|}\sum_{ \lambda \gqs \sqrt{k/2}/\max\{|\beta_1|,\gamma_1\}} \frac{1}{\lambda^2}\right)\notag\\
    & = \frac{\cl(2\pi m \gamma_1/k)}{|\beta_1 \gamma_1|} + O\left( \frac{1}{\sqrt{k}}\right). \label{urk}
\end{align}
\item[{\bf Case $(iia)$}]If $(\beta_0,\gamma_0) \in  Z_{\sqrt{k/2}}(h,k)$ then necessarily $(\beta_0,\gamma_0)=(\beta_1,\gamma_1)$ and so \eqref{smfx} and \eqref{urk} prove the proposition in this case.
\item[{\bf Case $(iib)$}]In the final case, $Z_{\sqrt{k/2}}(h,k)$ is not empty and doesn't contain $(\beta_0,\gamma_0)$. Since $|\beta_1 \gamma_1| \gqs |\beta_0 \gamma_0| \gqs \sqrt{k/2}$ we find
\begin{equation} \label{mo4}
    \frac{\cl(2\pi m \gamma_1/k)}{|\beta_1 \gamma_1|} = O\left( \frac{1}{\sqrt{k}}\right)
\end{equation}
so that \eqref{mo3} follows from \eqref{smfx}, \eqref{mo2}, \eqref{urk} and \eqref{mo4}.
\end{enumerate}
We see that both sides of \eqref{mo3} are $O((\log k)/\sqrt{k})$ except in Case $(iia)$, and in this case the pair $(\beta_0,\gamma_0) \in  Z_{\sqrt{k/2}}(h,k)$ is unique.
\end{proof}

\begin{proof}[\bf Proof of Theorem \ref{mainest}]
The proof now follows directly from combining \eqref{cqor} and  Proposition \ref{l3}. Treating the error  in \eqref{mo3} of Proposition \ref{l3} more carefully, we find it is bounded by
\begin{equation*}
(2\sqrt{2}(5-\log 2 +\cl(\pi/3))+2\sqrt{2}\log k)/\sqrt{k} < (15.06+2\sqrt{2}\log k)/\sqrt{k}.
\end{equation*}
Combining this with the estimate \eqref{khbd} for the error in \eqref{cqor} shows that the error term in \eqref{logpx} of Theorem \ref{mainest} is bounded by $(16.05+\sqrt{2}/\pi \log k)/\sqrt{k}$.
\end{proof}

\begin{cor} \label{needc}
There exists an absolute constant $\tau$ such that for all integers $m$ with $0\lqs m \lqs k-1$
\begin{equation*}
\frac{1}{k} \Bigl| \log \bigl| \spn{h/k}{m} \bigr| \Bigr| \lqs  \frac{\cl(\pi/3)}{2\pi D(h,k)}  + \tau \frac{\log k}{\sqrt{k}}.
\end{equation*}
\end{cor}
\begin{proof}
We may take $\tau$ to be the absolute implied constant of Theorem \ref{mainest}   and note that $|\cl(\theta)|\lqs \cl(\pi/3)$ for all $\theta \in \R$. Hence we may take any $\tau > \sqrt{2}/\pi$ for $k$ large enough.
\end{proof}

\section{Bounds for most $Q_{hk\sigma}(N)$} \label{sec-bndq}
In this section we continue to assume that $h$ and $k$ are integers with $1\lqs h<k$ and $(h,k)=1$.

\subsection{Initial estimates}
The next result, mentioned in the introduction, is proved in this subsection.
\begin{mainp}
For $2 \lqs k \lqs N$, $\sigma \in \R$ and $s := \lfloor N/k \rfloor$
\begin{equation} \label{flgo}
    |Q_{hk\sigma}(N)|  \lqs \frac{3}{k^3} \exp\left( N \frac{2  +  \log \left(1+ 3k/4 \right)}{k} +\frac{|\sigma|}{N} \right)\left|\spr{h/k}{N-s k} \right|.
\end{equation}
\end{mainp}
\begin{proof}
From definition \eqref{qhkl},
\begin{equation}
Q_{hk\sigma}(N)   =      \int_{\mathcal L} \frac{e^{2\pi i \sigma z}}{(1-e^{2\pi i z})(1-e^{2\pi i2 z}) \cdots (1-e^{2\pi i N z})} \, dz \label{tg}
\end{equation}
where $z$ traces a loop  $\mathcal L$  of radius $1/(2\pi N k \lambda)$ around $h/k$, i.e.
 \begin{equation*}
    z=h/k+w, \qquad |w|=\frac 1{2\pi N k \lambda}
\end{equation*}
and $\lambda$ is large enough that only the pole of the integrand at $h/k$ is inside $\mathcal  L$. This is ensured when $\lambda> 1/2\pi$, since if $a/b$ is any other pole ($1 \lqs b \lqs N$) we have
$$
\left| \frac ab - \frac hk\right| = \left| \frac{ak-bh}{bk} \right| \gqs \frac 1{bk} \gqs \frac 1{Nk} > |w|.
$$
Therefore, letting $e^{2\pi i \sigma z} I_N(z)$ denote the integrand in \eqref{tg},
\begin{equation}\label{tg2}
|Q_{hk\sigma}(N)| \lqs \int_{\mathcal L} \left| e^{2\pi i \sigma z} I_N(z) \right| \, dz \lqs 2\pi \left( \frac 1{2\pi  N k \lambda}\right) \sup \bigl\{|e^{2\pi i \sigma z} I_N(z)| : z \in \mathcal L\bigr\}.
\end{equation}
It is easy to see that if $\lambda \gqs 1/k$ then
\begin{equation}\label{ets}
    |e^{2\pi i \sigma z} | \lqs e^{|\sigma|/N} \qquad (z \in \mathcal L, \ \sigma \in \R).
\end{equation}
Now write $I_N(z)=I_N^*(z) \cdot I_N^{**}(z)$ for
\begin{equation*}
    I_N^{*}(z):= \prod_{\substack{1 \lqs j \lqs N \\ k \mid j}}\frac{1}{(1-e^{2\pi i j z})}, \qquad
    I_N^{**}(z):= \prod_{\substack{1 \lqs j \lqs N \\ k \nmid j}}\frac{1}{(1-e^{2\pi i j z})}.
\end{equation*}
We use the following simple bounds, (better ones are proved in Lemma \ref{bbnds}). For all $z \in \C$ with $|z| \lqs 1$
\begin{align}
\left| 1-e^z  \right| & \lqs    2|z| \label{ine},\\
\left| 1-e^z  \right|^{-1} & \lqs    2/|z|  \label{iner},\\
\left|\log (1-z/2) \right| & \lqs    3|z|/4 \label{inel}.
\end{align}

\begin{lemma} \label{scroo1}
For $z \in \mathcal L$ and $\lambda \gqs 1/k$ we have
\begin{equation}\label{i1}
    \left| I_N^{*}(z) \right| \lqs \frac{e}{\sqrt{2\pi}}  \left(\frac{k}{ N}\right)^{1/2} (2 e k \lambda)^{s}.
\end{equation}
\end{lemma}
\begin{proof}
Clearly
$$
I_N^{*}(z) = \prod_{\substack{1 \lqs j \lqs N \\ k \mid j}}\frac{1}{(1-e^{2\pi i j (h/k+w)})} =
\prod_{1 \lqs m \lqs s}\frac{1}{(1-e^{2\pi i k m w})}.
$$
Also
\begin{equation} \label{kl1}
    |2\pi i k m w| = \frac{2\pi k m}{2\pi N k \lambda} \lqs \frac s{N\lambda} \lqs \frac{1}{k \lambda},
\end{equation}
so assuming $\lambda \gqs 1/k$,
we can apply \eqref{iner} to get
$$
\left| I_N^{*}(z) \right| \lqs \prod_{1 \lqs m \lqs s} \frac{2}{2\pi  k m |w|}= \prod_{1 \lqs m \lqs s} \frac{2 N \lambda}{m}
=\frac{(2 N \lambda)^s}{s !}.
$$
It follows from Stirling's formula that
$
1/a! < \frac{1}{\sqrt{2\pi a}}\left(\frac{e}{a}\right)^a$ for $a \in \Z_{\gqs 1}$.
Hence the lemma is obtained with
$$
\frac{1}{s!} = \frac{s+1}{(s+1)!}< \frac{s+1}{\sqrt{2\pi(s+1)}}\left(\frac{e}{s+1}\right)^{s+1} = \frac{e}{\sqrt{2\pi(s+1)}}\left(\frac{e}{s+1}\right)^{s} < \frac{e}{\sqrt{2\pi N/k}}\left(\frac{e k}{N}\right)^{s}. \qedhere
$$
\end{proof}

\begin{lemma} \label{scroo2}
For $z \in \mathcal L$ and $\lambda \gqs 1$  we have
\begin{equation*}
    \left|  I_N^{**}(z) \right| \lqs \exp \left(\frac{N}{2 k \lambda} + \frac {3N}{8 \lambda}\right) \frac{1}{k^{s}} \left|\spr{h/k}{N-s k} \right|.
\end{equation*}
\end{lemma}
\begin{proof}
Write
\begin{align}
    I_N^{**}(z) & = \prod_{\substack{1 \lqs j \lqs N \\ k \nmid j}}\frac{1}{(1-e^{2\pi i j (h/k+w)})} \notag\\
    & =  \prod_{\substack{1 \lqs j \lqs N \\ k \nmid j}} e^{-2\pi i j w}
    \prod_{\substack{1 \lqs j \lqs N \\ k \nmid j}}\frac{1}{(1-e^{2\pi i j h/k} -1+ e^{-2\pi i j w})} \notag\\
    & =   e^{-\pi i  w(N(N+1)-k s(s+1))}
    \prod_{\substack{1 \lqs j \lqs N \\ k \nmid j}} \frac{1}{(1-e^{2\pi i j h/k})}
    \prod_{\substack{1 \lqs j \lqs N \\ k \nmid j}}\frac{1}{(1 - \eta_{h/k}(j,w))} \label{3p}
\end{align}
for
$$
\eta_{h/k}(j,w):=\frac{1 - e^{-2\pi i j w}}{1-e^{2\pi i j h/k}}.
$$
To estimate the parts of \eqref{3p}, we start with
\begin{equation}\label{nvl}
N(N+1)-k s(s+1) \lqs N^2, \qquad (k \lqs N)
\end{equation}
to see that
\begin{equation}\label{i2}
\left|   e^{-\pi i  w(N(N+1)-k s(s+1))} \right|  \lqs \exp \left(\frac{N}{2 k \lambda}\right).
\end{equation}
With $(1-\zeta)(1-\zeta^2) \cdots (1-\zeta^{k-1})= k$
for $\zeta$  a primitive $k$th root of unity, (by \cite[Lemma 4.4]{OS} for example), the middle product satisfies
\begin{equation}\label{i3}
  \prod_{\substack{1 \lqs j \lqs N \\ k \nmid j}} \frac{1}{\left| (1-e^{2\pi i j h/k})\right|}   = \frac{1}{k^{s}} \left|\spr{h/k}{N-s k} \right| .
\end{equation}

Next we estimate the right-hand product  of \eqref{3p}. By \eqref{ine}
\begin{equation}\label{ala1}
    \left|1 - e^{-2\pi i j w} \right|  \lqs 2 \cdot 2\pi  j |w| = \frac{2j}{ N k \lambda}
\end{equation}
provided $\lambda \gqs 1/k$.
We have
\begin{equation*}
    \frac 1{|1-e^{-2\pi i \theta}|} = \frac 1{2|\sin (\pi \theta)|} \lqs \frac 1{4|\theta|} \qquad (-1/2 \lqs \theta \lqs 1/2)
\end{equation*}
and it follows that
\begin{equation} \label{ala2}
    \frac 1{\left|1 - e^{-2\pi i j h/k} \right|}  \lqs \frac 1{\left|1 - e^{-2\pi i /k} \right|} \lqs \frac{k}{4} \qquad (k\gqs 2).
\end{equation}
Consequently, \eqref{ala1}, \eqref{ala2} show
\begin{equation}\label{allll}
|\eta_{h/k}(j,w)| \lqs \frac{j}{2  N \lambda}.
\end{equation}
If
$
\lambda \gqs 1
$
 then $|\eta_{h/k}(j,w)| \lqs 1/2$ for all $j \lqs N$ and we may apply \eqref{inel}:
\begin{align}
    \prod_{\substack{1 \lqs j \lqs N \\ k \nmid j}}\frac{1}{\left|1 - \eta_{h/k}(j,w) \right|} & = \exp\left(-\sum_{1 \lqs j \lqs N, \  k \nmid j} \log \left|1 - \eta_{h/k}(j,w) \right| \right) \notag\\
    & \lqs \exp\left(\sum_{1 \lqs j \lqs N, \  k \nmid j} \left|\log (1 - \eta_{h/k}(j,w)) \right| \right) \notag\\
    & \lqs \exp\left(\frac 32\sum_{1 \lqs j \lqs N, \  k \nmid j}  \left| \eta_{h/k}(j,w) \right| \right) \lqs
    \exp\left(\frac {3N}{8 \lambda} \right) \label{i4}
\end{align}
where we used \eqref{nvl} in the last inequality.
Combining the estimates \eqref{i2}, \eqref{i3} and \eqref{i4} for \eqref{3p} finishes the proof.
\end{proof}

Inserting the bounds from \eqref{ets} and Lemmas \ref{scroo1}, \ref{scroo2}  into \eqref{tg2},  we obtain
\begin{equation} \label{i5}
    |Q_{hk\sigma}(N)|  \lqs  \frac{e}{\sqrt{2\pi} N^{3/2} k^{1/2} \lambda} \exp\left( N \left[ \frac {1}{2  k \lambda} +\frac {3}{8 \lambda} + \frac{1+\log 2\lambda}{k}\right] +\frac{|\sigma|}{N}\right)\left|\spr{h/k}{N-s k} \right|.
\end{equation}
For fixed $k$, the expression
$$
\frac {1}{2  k \lambda} +\frac {3}{8 \lambda} +  \frac{1+\log 2\lambda}{k}
$$
has its minimum at $\lambda = 1/2+3k/8$. We may set  $\lambda$ to this value in \eqref{i5} since all the conditions $\lambda \gqs 1/(2\pi)$, $1/k$, $1$ are satisfied when $k\gqs 2$. This completes the proof of Proposition \ref{prc}.
\end{proof}

An example of Proposition \ref{prc} is given in Figure \ref{prcfig} for $h=\sigma=1$ and $N=50$ where we denote the right side of \eqref{flgo} as $Q^*_{hk\sigma}(N)$. The numbers $Q_{hk\sigma}(N)$ are calculated using the methods of \cite[Sect. 5]{OS} as follows. For $N$, $k\gqs 1$, $m\gqs 0$ and $0\lqs r \lqs k-1$ define the rational numbers $E_k(N,m;r)$ recursively with $E_k(0,m;r)$ set as $1$  if $m=r=0$ and $0$ otherwise. Also
\begin{equation*}
E_k(N,m;r):=\sum_{a=0}^m \frac{N^a k^{a-1}}{a!} \sum_{j=0}^{k-1} E_k \bigl(N-1,m-a;(r-Nj)\bmod k \bigr) \cdot B_a(j/k) \qquad (N\gqs 1)
\end{equation*}
for $B_a(x)$ the Bernoulli polynomial. Then
\begin{equation} \label{niceformula}
    Q_{hk\sigma}(N) = \frac{(-1)^N}{N!} \sum_{r=0}^{k-1} e^{2\pi i (r+\sigma) h/k} \sum_{j=0}^{N-1} \frac{\sigma^j}{j!} E_k(N,N-1-j;r).
\end{equation}
In particular, we see from \eqref{niceformula} that $e^{-2\pi i \sigma h/k} Q_{hk\sigma}(N)$ is a polynomial in $\sigma$ of degree $N-1$.



\SpecialCoor
\psset{griddots=5,subgriddiv=0,gridlabels=0pt}
\psset{xunit=0.2cm, yunit=0.04cm}
\psset{linewidth=1pt}
\psset{dotsize=2pt 0,dotstyle=*}

\begin{figure}[h]
\begin{center}
\begin{pspicture}(-5,-20)(60,50) 

\psaxes[linecolor=gray,Ox=0,Oy=-10,Dx=10,dx=10,Dy=10,dy=10]{->}(0,-10)(-3,-15)(52,50)

\newrgbcolor{darkgrn}{0 0.5 0}
\newrgbcolor{lightgrn}{0 1 0}

\savedata{\mydata}[
{ {3, 49.7017}, {4, 38.2487}, {5, 31.8717}, {6,
  26.0669}, {7, 22.7983}, {8, 19.4631}, {9, 15.2377}, {10,
  14.9112}, {11, 11.2459}, {12, 12.2498}, {13, 6.94855}, {14,
  6.73311}, {15, 7.98948}, {16, 8.27486}, {17, 3.37705}, {18,
  1.10356}, {19, 1.28922}, {20, 2.21845}, {21, 3.29326}, {22,
  4.15042}, {23, 4.47975}, {24, 3.88688}, {25,
  1.42829}, {26, -3.68176}, {27, -5.31492}, {28, -5.91948}, {29,
-5.9734}, {30, -5.69085}, {31, -5.19492}, {32, -4.56561}, {33,
-3.859}, {34, -3.11657}, {35, -2.37026}, {36, -1.64564}, {37,
-0.963863}, {38, -0.343114}, {39, 0.200362}, {40, 0.651407}, {41,
  0.995184}, {42, 1.21641}, {43, 1.29847}, {44, 1.22227}, {45,
  0.964511}, {46,
  0.494896}, {47, -0.229066}, {48, -1.27321}, {49, -2.75632}, {50,
-4.9668}}
  ]
\dataplot[linecolor=orange,linewidth=0.8pt,plotstyle=line]{\mydata}
\dataplot[linecolor=orange,linewidth=0.8pt,plotstyle=dots]{\mydata}

\savedata{\mydatab}[
{{2, -2.33607}, {3, -3.21032}, {4, -3.99423}, {5, -4.47876}, {6,
-4.8469}, {7, -4.93827}, {8, -5.15406}, {9, -6.44913}, {10,
-5.74208}, {11, -6.97482}, {12, -5.26902}, {13, -8.36114}, {14,
-7.84872}, {15, -5.98189}, {16, -5.19056}, {17, -8.34424}, {18,
-10.1263}, {19, -9.43756}, {20, -8.05185}, {21, -6.56582}, {22,
-5.33812}, {23, -4.67156}, {24, -4.94005}, {25, -6.94793}, {26,
-11.1969}, {27, -12.5013}, {28, -12.797}, {29, -12.5602}, {30,
-12.0033}, {31, -11.248}, {32, -10.373}, {33, -9.43328}, {34,
-8.46933}, {35, -7.51217}, {36, -6.58657}, {37, -5.71296}, {38,
-4.90886}, {39, -4.18994}, {40, -3.57081}, {41, -3.06582}, {42,
-2.68982}, {43, -2.459}, {44, -2.3921}, {45, -2.51208}, {46,
-2.84891}, {47, -3.44481}, {48, -4.36533}, {49, -5.72901}, {50,
-7.82405}}
  ]
\dataplot[linecolor=darkgrn,linewidth=0.8pt,plotstyle=line]{\mydatab}
\dataplot[linecolor=darkgrn,linewidth=0.8pt,plotstyle=dots]{\mydatab}

\rput(53,-20){$k$}
\rput(60,35){$\log Q^*_{1k1}(50)$}
\rput(60,15){$\log |Q_{1k1}(50)|$}
\psline[linecolor=orange](50,35)(54,35)
\psline[linecolor=darkgrn](50,15)(54,15)


\end{pspicture}
\caption{Bounding $Q_{hk\sigma}(50)$ for $h=\sigma=1$ and $2 \leqslant k \lqs 50$}\label{prcfig}
\end{center}
\end{figure}
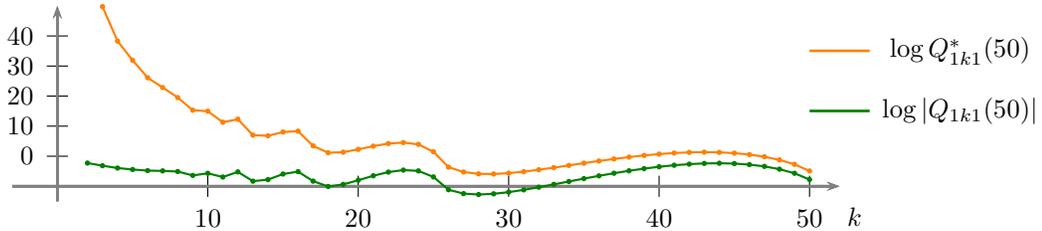


\subsection{Improved estimates}

By tightening up the bounds \eqref{ine},  \eqref{iner}, \eqref{inel} and restricting the range of $k$ we can improve Proposition \ref{prc} a little as follows.

\begin{lemma} \label{bbnds}
For $z \in \C$ and $|z| \lqs Y$ we have
\begin{alignat}{2}
\left|\frac{ 1-e^z }{z} \right| & \lqs   \alpha(Y):=\frac{e^Y-1}{Y}   & &  \label{ine2}\\
\left|\frac{z}{ 1-e^z } \right| & \lqs   \beta(Y):=2+\frac{Y}{2} \left(1-\cot \left(\frac{Y}{2}\right)\right) \qquad \qquad  & & (Y < 2\pi)  \label{iner2}\\
\left|\frac{\log(1-z)}{z} \right| & \lqs   \gamma(Y):=\frac{1}{Y} \log \left(\frac{1}{1-Y} \right)   & & (Y  <1) \label{inel2}.
\end{alignat}
\end{lemma}
\begin{proof}
For $|z| \lqs Y < 2\pi$ we have
\begin{equation*}
    \left|\frac{z}{ 1-e^z } \right|  = \left|\sum_{n=0}^\infty B_n \frac{z^n}{n!} \right|
     \lqs \sum_{n=0}^\infty |B_n| \frac{Y^n}{n!}
     = 1+\frac{Y}{2} +\left( 1- \frac{Y}{2}\cot \left(\frac{Y}{2}\right) \right),
\end{equation*}
using \cite[Eq. (11.1)]{Ra}. The other two inequalities have similar proofs. Note that for $Y=0$ we have $\alpha(0)=\beta(0)=\gamma(0)=1$ in the limit, with $\alpha(Y)$, $\beta(Y)$ and $\gamma(Y)$ increasing for $Y\gqs 0$.
\end{proof}

Start with a parameter $K\gqs 2$. We assume
\begin{equation}\label{snow}
    k \gqs K, \qquad \lambda \gqs 1/2+  K/8.
\end{equation}
The quantity $1/(k\lambda)$ in \eqref{kl1} then satisfies
\begin{equation*}
    \frac{1}{k\lambda} \lqs \frac{1}{K(1/2+ K/8)} < 2\pi.
\end{equation*}
With \eqref{iner2} we may therefore replace the factor $2$ in \eqref{i1} by
\begin{equation} \label{xi1}
    \xi_1=\xi_1(K):=\beta\left( \frac{1}{K(1/2+ K/8)}\right).
\end{equation}
Similarly, the factor $2$ in \eqref{ala1} may be replaced by
\begin{equation} \label{xi2}
    \xi_2=\xi_2(K):=\alpha\left( \frac{1}{K(1/2+ K/8)}\right).
\end{equation}
This improves the bound \eqref{allll} to
\begin{equation*}
|\alpha_{h/k}(j,w)| \lqs \frac{\xi_2 j}{4  N \lambda}
\end{equation*}
so that for all $j \lqs N$ we have $|\alpha_{h/k}(j,w)| \lqs \xi_2 /(4   \lambda)<1$. The factor $3/2$ in \eqref{i4} can now be
replaced by
\begin{equation} \label{xi3}
    \xi_3=\xi_3(K):=\gamma\left( \frac{\xi_2}{4(1/2+ K/8)}\right)
\end{equation}
and we obtain
\begin{equation*}
    \prod_{\substack{1 \lqs j \lqs N \\ k \nmid j}}\frac{1}{\left|1 - \alpha_{h/k}(j,w) \right|} \lqs \exp\left(\frac{\xi_2 \xi_3 N}{8 \lambda} \right).
\end{equation*}
Hence
\begin{equation} \label{ii6}
    |Q_{hk\sigma}(N)|  \lqs  \frac{e}{\sqrt{2\pi} N^{3/2} k^{1/2} \lambda} \exp\left( N \left[ \frac {1}{2  k \lambda} +\frac {\xi_2 \xi_3}{8 \lambda} + \frac{1+\log \xi_1 \lambda}{k}\right] +\frac{|\sigma|}{N}\right)\left|\spr{h/k}{N-s k} \right|
\end{equation}
and setting $\lambda = 1/2+\xi_2 \xi_3 k/8$ minimizes \eqref{ii6}. Note that $\xi_2 \xi_3 \gqs 1$ so that our initial inequality \eqref{snow} for $\lambda$ is true. We have proved

\begin{prop} \label{prc22}
For $2 \lqs K \lqs k \lqs N$ and $s := \lfloor N/k \rfloor$ we have
\begin{equation} \label{i7}
    |Q_{hk\sigma}(N)|  \lqs \frac{9}{k^3} \exp\left( N \frac{2  +  \log \left(\xi_1/2+ \xi_1 \xi_2 \xi_3 k/8 \right)}{k} +\frac{|\sigma|}{N} \right)\left|\spr{h/k}{N-s k} \right|
\end{equation}
for $\xi_1$, $\xi_2$, $\xi_3$ defined in \eqref{xi1}, \eqref{xi2}, \eqref{xi3} and depending on $K$.
\end{prop}

Some examples of triples $(K, \ \xi_1, \ \xi_1 \xi_2 \xi_3)$ are
\begin{alignat}{3}
    K & =2: \quad & \xi_1 &\approx 1.37065, & \quad \xi_1 \xi_2 \xi_3 &  \approx 2.64070 \label{kx1}\\
    K & =61: \quad & \xi_1 &\approx 1.00101, & \quad \xi_1 \xi_2 \xi_3 &  \approx 1.01778 \label{kx2}\\
    K & =82: \quad & \xi_1 &\approx 1.00057, & \quad \xi_1 \xi_2 \xi_3 &  \approx 1.01297 \label{kx3}\\
    K & =101: \quad & \xi_1 &\approx 1.00038, & \quad \xi_1 \xi_2 \xi_3 &  \approx 1.01041. \label{kx4}
\end{alignat}

\subsection{Final bounds}

Define $\mathcal B(K,N)$ to be the set
\begin{subequations} \label{subbe}
\begin{align}
 \Bigl\{ h/k \ :  \ K  \lqs k \lqs N,  & \ 0\lqs h < k, \ (h,k)=1 \Bigr\} \label{subbe1}\\
 \intertext{but with the restrictions}
    h \not\equiv \pm 1 \bmod k & \qquad \text{if} \qquad N/3  <k \lqs N/2, \label{subbe2}\\
    h \not\equiv \pm 1, \pm 2, (k\pm 1)/2 \bmod k & \qquad \text{if} \qquad N/2  <k \lqs N. \label{subbe3}
\end{align}
\end{subequations}

\begin{theorem} \label{ewu}
There exists $W<U:=-\log |w_0| \approx 0.068076$ so that
\begin{equation*}
\sum_{h/k \in \mathcal B(101,N)}  Q_{hk\sigma}(N) = O(e^{WN}).
\end{equation*}
We may take any $W > \cl(\pi/3)/(6\pi) \approx 0.0538$ and the implied constant depends only on $\sigma$ and $W$.
\end{theorem}
\begin{proof}
Recall from Corollary \ref{needc}
that
there exists an absolute constant $\tau$ such that for all $m$, $h$, $k\in \Z$ with $1 \lqs  h < k$, $(h,k)=1$ and $0 \lqs m < k$ we have
\begin{equation}\label{logp}
\log \left| \spr{h/k}{m} \right| \ \ \lqs \ \ \frac{\cl(\pi/3)}{2\pi D(h,k)}  \cdot  k + \tau \sqrt{k} \log k.
\end{equation}
 It follows from Proposition \ref{prc22} and  \eqref{logp} that
  \begin{equation*}
    Q_{hk\sigma}(N)  \ll \frac{1}{k^3} \exp\left( N \frac{2  +  \log \left(\xi_1/2+ \xi_1 \xi_2 \xi_3 k/8 \right)}{k}  + \frac{\cl(\pi/3)}{2\pi D(h,k)}   \cdot k + \tau \sqrt{N} \log N \right)
  \end{equation*}
 where $k \gqs K=101$ and $\xi_1$, $\xi_1 \xi_2 \xi_3$ are as in \eqref{kx4}.
 Given any $\epsilon>0$ we have $ \tau \sqrt{N} \log N  \lqs \epsilon N$ for $N$ large enough.
 For $k$ in a range $0<a \lqs k \lqs b$ where we know $D(h,k) \gqs D^*$, the expression
 \begin{equation} \label{mabel}
    N \frac{2  +  \log \left(\xi_1/2+ \xi_1 \xi_2 \xi_3 k/8 \right)}{k}  + \frac{\cl(\pi/3)}{2\pi D^*}   \cdot k
 \end{equation}
  has possible maxima only at the end points $k=a$ or $k=b$. For $h/k \in \mathcal B(101,N)$ with $101 \lqs k\lqs N/3$ we know $D(h,k) \gqs 1$ and see the end points are bounded by
 \begin{align}
     N \frac{2  +  \log \left(\xi_1/2+ \xi_1 \xi_2 \xi_3 101/8 \right)}{101}  + \frac{\cl(\pi/3)}{2\pi \cdot 1}   \cdot 101 & < 0.0454 N + 16.315, \label{numb}\\
      N \frac{2  +  \log \left(\xi_1/2+ \xi_1 \xi_2 \xi_3 (N/3)/8 \right)}{N/3}  + \frac{\cl(\pi/3)}{2\pi \cdot 1}   \cdot \frac{N}{3} & < 6+  \epsilon N + \frac{\cl(\pi/3)}{6\pi} N. \notag
 \end{align}
 Therefore
 \begin{equation*}
   Q_{hk\sigma}(N)  \ll \frac{1}{k^3} \exp\left( N \left[   \frac{\cl(\pi/3)}{6\pi} + 2\epsilon \right] \right) \qquad (h/k \in \mathcal B(101,N), \ k\lqs N/3).
 \end{equation*}
Similarly, for $h/k \in \mathcal B(101,N)$ with $N/3<k\lqs N/2$ we have $D(h,k) \gqs 2$ by \eqref{egd}. Hence \eqref{mabel} is bounded by the maximum of
\begin{equation*}
   6+\epsilon N + \frac{\cl(\pi/3)}{2\pi \cdot 2}     \cdot \frac{N}{3}, \qquad 4+\epsilon N + \frac{\cl(\pi/3)}{2\pi \cdot 2}   \cdot \frac{N}{2}.
\end{equation*}
For $h/k \in \mathcal B(101,N)$ with $N/2<k\lqs N$  we have $D(h,k) \gqs 3$ by \eqref{egd2}. Hence \eqref{mabel} is  bounded by the maximum of
\begin{equation*}
  4+ \epsilon N + \frac{\cl(\pi/3)}{2\pi \cdot 3}    \cdot \frac{N}{2}, \qquad 2+ \epsilon N + \frac{\cl(\pi/3)}{2\pi \cdot 3}    \cdot N
\end{equation*}

It follows that for any $W > \cl(\pi/3)/(6\pi)$
\begin{equation*}
    Q_{hk\sigma}(N)  \ll   e^{WN}/k^3 \qquad (h/k \in \mathcal B(101,N)).
\end{equation*}
 Finally,
\begin{align*}
    \sum_{h/k \in \mathcal B(101,N)}  Q_{hk\sigma}(N) & \ll  \sum_{h/k \in \mathcal B(101,N)}  e^{WN}/k^3 \\
     &  \ll e^{WN} \sum_{k=1}^N \sum_{h=1}^k 1/k^3 = e^{WN} \sum_{k=1}^N  1/k^2 \ll e^{WN}. \qedhere
\end{align*}
\end{proof}

\begin{remark} \label{K=2} \rm Theorem \ref{ewu} is still true if we enlarge $\mathcal B(101,N)$ to  $\mathcal B(82,N)$, i.e. allowing all $k \gqs 82$. This is because we obtain $0.0535 N + \dots$ on the right side of \eqref{numb} when we replace $101$ by $K=82$ on the left (and use the corresponding $\xi_i$s as in \eqref{kx3}).  Furthermore, with $K=61$ we find
 \begin{equation*}
\sum_{h/k \in \mathcal B(61,N)}  Q_{hk\sigma}(N)
 = O(e^{WN}),
\end{equation*}
needing $W \approx 0.067403$, very close to $U$ (see \eqref{kx2}). We expect that $K$ can be pushed all the way back to 2 and that with improved techniques it should be possible to prove that for some $W<U$
\begin{equation*}
\sum_{h/k \in \mathcal B(2,N)}  Q_{hk\sigma}(N)
 = O(e^{WN}).
\end{equation*}
This would eliminate the $\sum_{0<h/k \in \farey_{100}}$ term in \eqref{maineq2} of Theorem \ref{mainthmb}.
\end{remark}

What remains from $\farey_N- \bigl(\farey_{100} \cup \mathcal A(N) \cup \mathcal B(101,N)\bigr)$
are the subsets $\mathcal C(N)$, $\mathcal D(N)$ and $\mathcal E(N)$ as defined in \eqref{cnsub}, \eqref{dnsub} and \eqref{ensub}.
In the following sections we find the asymptotics for each of the corresponding $Q_{hk\sigma}(N)$ sums.

\section{Further required results}
We gather here some more results from \cite{OS1} we will require for developing the asymptotic expansions in the next sections.
Throughout we write $z=x+iy \in \C$.

\subsection{Some dilogarithm results} \label{dilogg}

In  \cite[Sect. 2.3]{OS1} we saw the identity
\begin{equation}\label{ss}
    \li \left(e^{-2\pi i z} \right) = -\li \left(e^{2\pi i z} \right) +2\pi^2\left(z^2-(2m+1)z +m^2+m+1/6 \right)
\end{equation}
for $m< \Re(z) < m+1$ where $m\in\Z$.
Also
\begin{equation}\label{dli}
\cl\left(2\pi  z\right) = -i \li\left( e^{2\pi i z}\right)+i\pi^2\left(z^2-(2m+1)z+m^2+m+1/6 \right)
\end{equation}
for $m\lqs z \lqs m+1$.

\begin{lemma} \label{dil1}
Consider $\Im(\li(e^{2\pi i z}))$  as a function of $y \in \R$. It is positive and decreasing for fixed $x\in (0,1/2)$ and negative and increasing for fixed $x\in (1/2,1)$.
\end{lemma}

\begin{lemma} \label{dil1a}
Consider $\Re(\li(e^{2\pi i z}))$  as a function of $y \gqs 0$. It is positive and decreasing for fixed $x$ with $|x| \lqs 1/6$. It is negative and increasing for fixed $x$ with $1/4 \lqs |x| \lqs 3/4$.
\end{lemma}

\begin{lemma} \label{dil2}
For $y\gqs 0$ we have $|\li(e^{2\pi i z})| \lqs \li(1)$.
\end{lemma}

\subsection{Approximating products of sines}
In the following, let $h$ and $k$ be relatively prime integers with $1\lqs h<k$.
From \cite[Sect. 2.1]{OS1} we have
\begin{prop} \label{simple}
For  $N/2<k \lqs N$
\begin{equation*}
    Q_{hk\sigma}(N)  = \frac{(-1)^{k+1}}{k^2}  \exp\left(\frac{-\pi i h(N^2+N-4\sigma)}{2k}\right)
     \exp\left(\frac{\pi i}{2}(2Nh+N+h+k-hk)\right) \spr{h/k}{N-k}.
\end{equation*}
\end{prop}

So estimating $Q_{hk\sigma}(N)$ requires these
further results on sine products  from \cite[Sect. 3]{OS1}:

\begin{prop} \label{prrh}
For $m$, $L \in \Z_{\gqs 1}$ and $-1/m<\theta<1/m$ with $\theta \neq 0$ we have
\begin{multline} \label{pimv2}
    \spn{\theta}{m} = \left(\frac{ \theta}{|\theta|}\right)^{m}\left(\frac{2 \sin(\pi m \theta)}{\theta}\right)^{1/2}  \exp\left(-\frac{\cl(2\pi m \theta)}{2\pi \theta} \right) \\
      \quad \times \exp\left(\sum_{\ell=1}^{L-1} \frac{B_{2\ell}}{(2\ell)!} ( \pi \theta)^{2\ell-1} \cot^{(2\ell-2)} (\pi m \theta)\right) \exp \bigl(T_{L}(m,\theta)\bigr)
  \end{multline}
for
\begin{equation*}
    T_{L}(m,\theta) := \left( \pi \theta\right)^{2L} \int_0^m \frac{B_{2L}-B_{2L}(x-\lfloor x \rfloor)}{(2L)!} \rho^{(2L)}(\pi x \theta) \, dx + \int_0^\infty \frac{B_{2L}-B_{2L}(x - \lfloor x \rfloor)}{2L(x+m)^{2L}} \, dx.
\end{equation*}
  \end{prop}

\begin{lemma} \label{melrh2}
For  $1 \lqs m <k/h$ we have
\begin{equation} \label{1/12}
    |T_1(m,h/k)| \lqs \pi^2 h/18 + 1/12.
\end{equation}
\end{lemma}

\begin{prop} \label{gprop}
Let $W>0$. For  $\delta$ satisfying $0<\delta \lqs 1/e$ and $\delta \log(1/\delta) \lqs W$ we have
\begin{equation*}
    \spr{h/k}{m}   \lqs c(h) \exp\left( \frac{kW}h \right) \quad \text{for} \quad 0 \lqs \frac {mh}k \lqs \delta,
    \quad  \frac 12 - \delta \lqs \frac {mh}k <1
\end{equation*}
and
\begin{equation*}
    c(h):=h^{1/2} \exp( \pi^2 h/18 +1/6)/2.
\end{equation*}
\end{prop}

\begin{prop} \label{dela}
Suppose   $\Delta$ and $W$ satisfy $0.0048 \lqs \Delta \lqs 0.0079$ and $\Delta \log 1/\Delta \lqs W$.
For the integers $h$, $k$, $s$ and $m$ we require
\begin{equation*}
    0<h<k \lqs s, \quad  R_\Delta \lqs s/h, \quad \Delta s/h \lqs m \lqs k/(2h).
\end{equation*}
 Then for $L:=\lfloor \pi e \Delta \cdot s/h \rfloor$ we have
\begin{align}
\left| \spr{h/k}{m}  T_L(m,h/k) \right| & \lqs (\pi^3/2) c(h) \cdot e^{sW/h}, \label{slma}\\
    \left| T_L(m,h/k) \right| & \lqs \pi^3/2. \label{slmb}
\end{align}
\end{prop}

See \cite[Sect. 3.4]{OS1} for the definition of $R_\Delta$. We will only use it in the case when $\Delta = 0.006$ and then $R_\Delta \approx 130.7$.

\begin{cor}\label{delacor}
Let $W, \Delta, s, h,k,m$ and $L$ be as in Proposition \ref{dela}. Suppose also that $0< u/v\lqs h/k$. Then
\begin{align}
\left| \spr{h/k}{m}  T_L(m,u/v) \right| & \lqs (\pi^3/2) c(h) \cdot e^{sW/h}, \label{slmax}\\
    \left| T_L(m,u/v) \right| & \lqs \pi^3/2. \label{slmbx}
\end{align}
\end{cor}

The main consequence of Propositions \ref{prrh} and \ref{dela} is:
\begin{prop}\label{hard}
For $W, \Delta, s, h,k,m$ and $L$  as in Proposition \ref{dela} we have
\begin{multline} \label{pimv9}
    \spr{h/k}{m}  = \left(\frac{h}{2k \sin(\pi m h/k)}\right)^{1/2}  \exp\left(\frac k{2\pi h} \cl\bigl(2\pi  m h/k \bigr) \right) \\
      \quad \times \exp\left(-\sum_{\ell=1}^{L-1} \frac{B_{2\ell}}{(2\ell)!} \left( \frac{\pi h}{k}\right)^{2\ell-1} \cot^{(2\ell-2)}\left(\frac{\pi m h}{k}\right)\right) +O\left(e^{sW/h}\right)
  \end{multline}
for an  implied constant depending only on $h$.
\end{prop}

\section{The sum $\mathcal C_1(N,\sigma)$} \label{c1ns}

Let $\sigma \in \Z$. In this section and the next we prove Theorem \ref{thmc}, giving the asymptotic expansion as $N \to \infty$ of
\begin{equation}\label{grig}
\mathcal C_1(N,\sigma)  := \sum_{h/k \in \mathcal C(N)} Q_{hk\sigma}(N) = 2 \Re \sum_{\frac{N}{2}  <k \lqs N, \ k \text{ odd}}  Q_{2k\sigma}(N).
\end{equation}
 Setting $h=2$ in Proposition \ref{simple}  yields
\begin{equation} \label{azi}
    Q_{2k\sigma}(N)=\frac 1{k^2} \exp\left(-\pi i \frac{N^2+N-4\sigma}k \right)\exp\left( \frac{\pi i}2 (5N+2-k)\right) \spr{2/k}{N-k}.
\end{equation}
The sum \eqref{grig} corresponds to $2N/k \in [2,4)$ and we break it into two parts:
 $\mathcal C_2(N,\sigma)$ for $2N/k \in [2,3)$ and $\mathcal C^*_2(N,\sigma)$  with $2N/k \in [3,4)$.


\SpecialCoor
\psset{griddots=5,subgriddiv=0,gridlabels=0pt}
\psset{xunit=0.8cm, yunit=0.8cm}
\psset{linewidth=1pt}
\psset{dotsize=2pt 0,dotstyle=*}

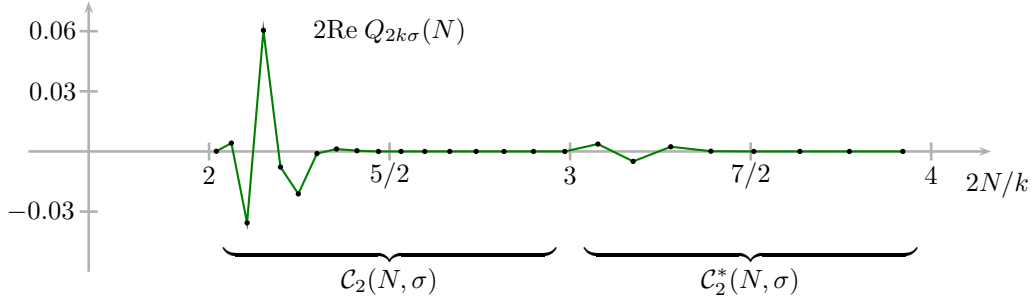
\begin{figure}[h]
\begin{center}
\begin{pspicture}(-2,-2.3)(16,2.5) 

\newgray{xgray}{.7}
\psline[linecolor=xgray]{->}(0,-2)(0,2.5)
\psline[linecolor=xgray]{->}(-1,0)(15,0)
\psline[linecolor=xgray](-0.15,1)(0.15,1)
\psline[linecolor=xgray](-0.15,2)(0.15,2)
\psline[linecolor=xgray](-0.15,-1)(0.15,-1)
\psline[linecolor=xgray](2,-0.15)(2,0.15)
\psline[linecolor=xgray](5,-0.15)(5,0.15)
\psline[linecolor=xgray](8,-0.15)(8,0.15)
\psline[linecolor=xgray](11,-0.15)(11,0.15)
\psline[linecolor=xgray](14,-0.15)(14,0.15)

\rput(2,-0.4){$2$}
  \rput(5,-0.4){$5/2$}
  \rput(8,-0.4){$3$}
  \rput(11,-0.4){$7/2$}
  \rput(14,-0.4){$4$}

  \rput(-0.6,2){$0.06$}
  \rput(-0.6,1){$0.03$}
  \rput(-0.8,-1){$-0.03$}
  \rput(15.1,-0.55){$2N/k$}

\savedata{\mydata}[
{{13.5294, 0}, {12.6415, 0}, {11.8182, 0}, {11.0526, 0}, {10.339,
  0.00411513}, {9.67213, 0.078065}, {9.04762, -0.1652}, {8.46154,
  0.123153}, {7.91045, 0}, {7.3913, 0}, {6.90141, 0}, {6.43836,
  0}, {6., 0}, {5.58442, 0}, {5.18987, 0}, {4.81481, 0}, {4.45783,
  0.0124227}, {4.11765,
  0.0408909}, {3.7931, -0.0349745}, {3.48315, -0.705286}, {3.18681,
-0.2596}, {2.90323, 2.01378}, {2.63158, -1.187}, {2.37113,
  0.13912}, {2.12121, 0.003367}}
]
\newrgbcolor{darkred}{0 0.5 0}
\dataplot[linecolor=darkred,linewidth=0.8pt,plotstyle=line]{\mydata}
\dataplot[linecolor=black,linewidth=0.8pt,plotstyle=dots]{\mydata}

\rput(5,-1.9){$\underbrace{\phantom{xxxxxxxxxxxxxxxxxxxxxx}}_{\displaystyle \mathcal C_2(N,\sigma)}$}
\rput(11,-1.9){$\underbrace{\phantom{xxxxxxxxxxxxxxxxxxxxxx}}_{\displaystyle \mathcal C_2^*(N,\sigma)}$}

\rput(5,2){$2 \Re \ Q_{2k\sigma}(N)$}

\end{pspicture}
\caption{$2 \Re \ Q_{2k\sigma}(N)$ for $\sigma=1$ and $N=100$}\label{b2}
\end{center}
\end{figure}


\subsection{First results for $\mathcal C_2(N,\sigma)$}

With \eqref{azi} we have
\begin{multline}\label{b1n}
\mathcal C_2(N,\sigma) = \Re \sum_{k \text{ odd}, \ 2N/k \in [2,3)} \frac{-2}{k^2}
\exp\left( N\left[ \frac{\pi i}{2}\left(-\frac{2N}{k}+5 -2\frac k{2N} \right)\right]\right) \\
\times \exp\left(\frac{-\pi i}{2}\frac{2N}{k}\right) \exp\left( \frac 1N\left[ 2 \pi i \sigma\frac{2N}{k}\right]\right)
 \spr{2/k}{N-k}.
\end{multline}
Define
\begin{equation}\label{grz}
g_\ell(z):=-\frac{B_{2\ell}}{(2\ell)!} \left( \pi z\right)^{2\ell-1} \cot^{(2\ell-2)}\left(\pi z\right)
\end{equation}
and set $z=z(N,k):=2N/k$. The analog of the sine product approximation,  \cite[Thm. 4.1]{OS1}, we need here  is:
\begin{theorem} \label{cplxb2}
Fix $W>0$. Let $\Delta$ be in the range $0.0048 \lqs \Delta \lqs 0.0079$ and set $\alpha = \Delta \pi e$.  Suppose $\delta$ and $\delta'$ satisfy
\begin{equation*}
    \frac{\Delta}{1-\Delta} < \delta \lqs \frac 1e, \ 0<\delta' \lqs \frac 1e \quad \text{ and  } \quad \delta \log 1/\delta, \ \ \delta' \log 1/\delta' \lqs W.
\end{equation*}
Then for all $N \gqs 2 \cdot R_\Delta$  we have
\begin{equation}\label{ott1b2}
    \spr{2/k}{N-k} = O\left(e^{WN/2}\right) \quad \text{ for } \quad  z \in \left[2, \ 2+\delta \right] \cup  \left[5/2 -\delta' , \ 3 \right)
\end{equation}
and
\begin{multline} \label{ott2b2}
    \spr{2/k}{N-k} =  \frac{1}{N^{1/2}}  \exp\left(N\frac {\cl(2\pi  z)}{2\pi z}  \right) \left( \frac{z}{2\sin(\pi z)}\right)^{1/2}
      \\
      \quad \times \exp\left(\sum_{\ell=1}^{L-1} \frac{g_{\ell}(z)}{N^{2\ell-1}} \right) +O\left(e^{WN/2}\right)
        \quad \text{ for } \quad z \in  \left(2+\delta, \ 5/2 - \delta' \right)
  \end{multline}
with $L=\lfloor \alpha \cdot N/2\rfloor$. The implied constants in \eqref{ott1b2}, \eqref{ott2b2} are absolute.
\end{theorem}
\begin{proof}
The bound \eqref{ott1b2} follows directly from Proposition \ref{gprop} with $m=N-k$ and $h=2$. Next, in Proposition \ref{hard}, we set  $s=N$ and again $m=N-k$ and $h=2$.  The condition on $m$ in Proposition \ref{hard} is equivalent to
$$
2+\frac{\Delta}{1-\Delta/2} \lqs \frac{2 N}k \lqs \frac 52.
$$
So \eqref{ott2b2} follows from Proposition \ref{hard} if
\begin{equation}\label{dcv}
    \frac{\Delta}{1-\Delta/2} \lqs \delta.
\end{equation}
The inequality \eqref{dcv} is equivalent to $1/\Delta - 1/\delta \gqs 1/2$. Since our assumption
$\Delta/(1-\Delta)<\delta$  is equivalent to $1/\Delta - 1/\delta >1$, we have that \eqref{dcv} is true.
\end{proof}

With \eqref{dli} for $m=2$ we obtain
\begin{equation*}
    \cl(2\pi  z) = -i \li(e^{2\pi i z})+i \pi^2(z^2-5z+37/6) \qquad (2<z<3).
\end{equation*}
Therefore
\begin{equation} \label{pqzo}
    \frac {\cl(2\pi  z)}{2\pi z} + \frac{\pi i}{2}\left(-z+5 -\frac 2z \right) = \frac{1}{2\pi i z} \Bigl[  \li(e^{2\pi i z}) -\li(1) -4\pi^2 \Bigr],
\end{equation}
with the right side of \eqref{pqzo} now holomorphic in the strip $2<\Re(z)<3$.

To combine \eqref{b1n} and \eqref{ott2b2} we set, initially with $z\in (2,3)$,
\begin{align}
r_\mathcal C(z)
&  := \frac{1}{2\pi i z} \Bigl[  \li(e^{2\pi i z}) -\li(1) -4\pi^2 \Bigr], \label{rbw}\\
  q_\mathcal C(z) & := \left( \frac{z }{2\sin(\pi z)}\right)^{1/2} \exp(-\pi i z/2 ),  \label{qbw}\\
  v_\mathcal C(z;N,\sigma)  & := \frac{2\pi i \sigma z}N +\sum_{\ell=1}^{L-1} \frac{g_{\ell}(z)}{N^{2\ell-1}} \qquad (L= \lfloor \alpha \cdot N/2\rfloor). \label{vbw}
\end{align}
Then define
\begin{equation}\label{c3sum}
    \mathcal C_3(N,\sigma)  := \frac{-2}{N^{1/2}} \Re \sum_{k  \text{ odd} : \  z \in  (2+\delta ,  5/2 -\delta')}
    \frac{1}{k^2} \exp \bigl(N \cdot r_\mathcal C \left(z\right)\bigr) q_\mathcal C \left(z\right) \exp \bigl(v_\mathcal C (z;N,\sigma)\bigr),
\end{equation}
and it follows from \eqref{b1n} and Theorem \ref{cplxb2} that for  an absolute implied constant
\begin{equation}\label{cosmos}
\mathcal C_2(N,\sigma) = \mathcal C_3(N,\sigma) +O(e^{WN/2}).
\end{equation}

\subsection{Expressing $\mathcal C_3(N,\sigma)$ as an integral}

\begin{prop} \label{addx}
Suppose $3/2 \lqs \Re(z) \lqs 5/2$ and $|z-2| \gqs \varepsilon >0$. Also assume that
\begin{equation} \label{asmx}
    \max \Bigl\{1+\frac{1}{\varepsilon}, \ 16\Bigr\} \ < \frac{\pi e}{\alpha}.
\end{equation}
Then, for an implied constant depending only on $\varepsilon$, $\alpha$ and $d$,
\begin{equation} \label{we}
    \sum_{\ell=d}^{L-1} \frac{g_{\ell}(z)}{N^{2\ell-1}} \ll \frac{1}{N^{2d-1}}e^{-\pi|y|} \qquad (d \gqs 2, \ L=\lfloor \alpha \cdot N/2 \rfloor).
\end{equation}
\end{prop}
\begin{proof}
For $z$ in this range, Theorem 3.3 of \cite{OS1} bounding derivatives of the cotangent allows us to show
 \begin{equation}\label{hoot}
    \frac{g_{\ell}(z)}{N^{2\ell-1}} \ll F_{N,\varepsilon}(2\ell-1)\cdot e^{-\pi|y|}, \qquad \text{for} \qquad F_{N,\varepsilon}(\ell):=
    \left(\frac{\ell}{2\pi e N} \right)^\ell\left(\left(1+\frac 2\varepsilon\right)^\ell +32^\ell \right).
 \end{equation}
 This bound gets very large for $\ell$  large. The condition \eqref{asmx} ensures  $L$ is small enough that $g_{\ell}(z)/N^{2\ell-1}$ remains small. See \cite[Sect. 4.2]{OS1} for the details.
\end{proof}

We now fix some of the parameters in Theorem \ref{cplxb2}  and take
\begin{equation}\label{fix}
    W=0.05, \quad \alpha = 0.006 \pi e  \approx 0.0512, \quad 0.0061 \lqs \delta, \ \delta' \lqs 0.01, \quad N\gqs 400.
\end{equation}
Also, with $\varepsilon = 0.0061$, condition \eqref{asmx} is satisfied and Proposition \ref{addx} implies:

\begin{cor} \label{acdcx}
With $\delta, \delta' \in [0.0061, 0.01]$ and  $z \in \C$ such that $2+\delta \lqs \Re(z) \lqs 5/2-\delta'$  we have
\begin{equation*}
    v_\mathcal C(z;N,\sigma) =  \frac{2\pi i \sigma z}N +\sum_{\ell=1}^{d-1} \frac{g_{\ell}(z)}{N^{2\ell-1}} + O\left( \frac 1{N^{2d-1}}\right)
\end{equation*}
 for $2 \lqs d\lqs L=\lfloor 0.006 \pi e \cdot N/2 \rfloor$ and an implied constant depending only on  $d$.
\end{cor}

In the next theorem we assemble the results we need to convert the sum $\mathcal C_3(N,\sigma)$ in \eqref{c3sum} into an integral.

\begin{theorem} \label{rzjb}
The functions $r_\mathcal C(z)$, $q_\mathcal C(z)$ and $v_\mathcal C(z;N,\sigma)$ are holomorphic for $2< \Re(z)<5/2$. In this strip,
\begin{alignat}{2}
    \Re\left(r_\mathcal C \left(z \right) +\frac{2\pi i j}{z}\right) & \lqs \frac{1}{2\pi |z|^2} \left(x \cl(2\pi x) +\pi^2 |y|\left[\frac 13 +4(j+1) \right]  \right) \qquad & & (y\gqs 0) \label{efb1}\\
   \Re\left(r_\mathcal C \left(z \right) +\frac{2\pi i j}{z}\right) & \lqs \frac{1}{2\pi |z|^2} \left(x \cl(2\pi x) +\pi^2 |y|\left[\frac 13 -4j \right] \right) \qquad & & (y\lqs 0) \label{efb2}
\end{alignat}
for $j\in \R$.
Also, in the box with $2+\delta \lqs \Re(z) \lqs 5/2-\delta'$ and $-1\lqs \Im(z) \lqs 1$,
\begin{equation} \label{efb3}
    q_\mathcal C(z), \quad \exp \bigl(v_\mathcal C (z;N,\sigma)\bigr) \ll 1
\end{equation}
for an  implied constant depending only on $\sigma \in \R$.
\end{theorem}
\begin{proof}
Since $\li(e^{2\pi i z})$ is holomorphic away from the vertical branch cuts $(-i\infty,n]$ for $n \in \Z$, we see that $r_\mathcal C(z)$ is holomorphic for $2< \Re(z)<5/2$. Then in this strip, using \eqref{ss},
\begin{align}
     r_\mathcal C \left(z \right) +\frac{2\pi i j}{z} & = \frac 1{2\pi i z} \Bigl[ \li(e^{2\pi i z}) -\li(1) -4\pi^2(j+1) \Bigr] \notag \\
     & = \frac 1{2\pi i z} \Bigl[- \li(e^{-2\pi i z}) +\li(1)-4\pi^2 (j-2) \Bigr]-  \pi i(z-5). \label{fdrb}
\end{align}
The inequalities \eqref{efb1} and \eqref{efb2} follow, as in \cite[Sect. 4.3]{OS1}.

Check that for $w\in \C$,
\begin{equation*}
    -\pi/2< \arg \bigl(\sin(\pi w) \bigr) <\pi/2 \quad \text{for} \quad 0<\Re(w)<1.
\end{equation*}
Consequently, $-\pi<\arg\bigl( z/\sin(\pi z)\bigr) <\pi$ for $2< \Re(z)<5/2$ and so $q_\mathcal C(z)$ is holomorphic in this strip. Also $v_\mathcal C(z;N,\sigma)$ is holomorphic here since the only poles of $g_\ell(z)$ are at $z\in \Z$.

Finally, $q_\mathcal C(z)$ is bounded on the compact box, as is $\exp \bigl(v_\mathcal C (z;N,\sigma)\bigr)$ by Corollary \ref{acdcx}.
\end{proof}

By the calculus of residues, see for example \cite[p. 300]{Ol},
\begin{equation} \label{tantan}
\sum_{a \lqs k \lqs b, \ k \text{ odd}}  \varphi(k) =  \frac 12\int_C \frac{\varphi(z)}{2i\tan(\pi (z-1)/2)} \, dz
\end{equation}
for $\varphi(z)$  a holomorphic function and $C$ a positively oriented closed contour surrounding the interval $[a,b]$ and not surrounding any integers outside this interval.
Hence
\begin{equation*}
    \sum_{a \lqs k \lqs b, \ k \text{ odd}} \frac{1}{k^2} \varphi(2N/k) =  \frac{-1}{4N}\int_C \frac{\varphi(z)}{2i\tan(\pi (2N/z-1)/2)} \, dz
\end{equation*}
for $C$ now surrounding $\{2N/k \ | \ a\lqs k \lqs b\}$ with $a>0$.
Therefore
\begin{equation}\label{shir}
    \mathcal C_3(N,\sigma)  = \frac{1}{2N^{3/2}} \Re \int_{C_1}
     \exp \bigl(N \cdot r_\mathcal C(z)\bigr) \frac{q_\mathcal C(z)}{2i\tan \bigl(\pi (2N/z-1)/2 \bigr)} \exp \bigl(v_\mathcal C(z;N,\sigma)\bigr)\, dz
\end{equation}
where $C_1$ is the positively oriented rectangle with horizontal sides $C_1^+$, $C_1^-$ having imaginary parts $1/N^2$, $-1/N^2$ and vertical sides $C_{1,L}$, $C_{1,R}$ having real parts $2+\delta$ and $5/2-\delta'$ respectively, as shown in Figure \ref{c1c2}.
\SpecialCoor
\psset{griddots=5,subgriddiv=0,gridlabels=0pt}
\psset{xunit=0.8cm, yunit=0.8cm}
\psset{linewidth=1pt}
\psset{dotsize=2pt 0,dotstyle=*}
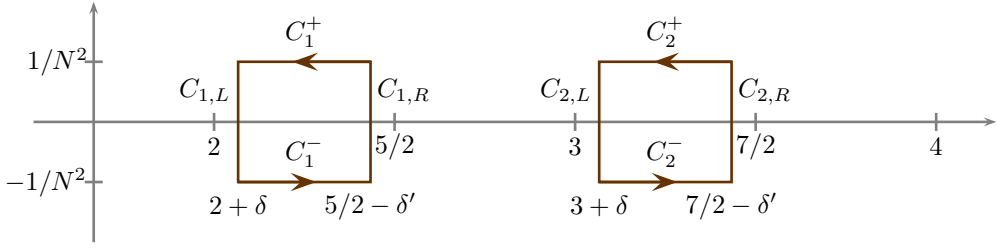
\begin{figure}[h]
\begin{center}
\begin{pspicture}(-2,-2)(16,2) 

\newrgbcolor{darkbrn}{0.4 0.2 0}

\psline[linecolor=gray]{->}(0,-2)(0,2)
\psline[linecolor=gray]{->}(-1,0)(15,0)
\psline[linecolor=gray](-0.15,1)(0.15,1)
\psline[linecolor=gray](-0.15,-1)(0.15,-1)
\psline[linecolor=gray](2,-0.15)(2,0.15)
\psline[linecolor=gray](5,-0.15)(5,0.15)
\psline[linecolor=gray](8,-0.15)(8,0.15)
\psline[linecolor=gray](11,-0.15)(11,0.15)
\psline[linecolor=gray](14,-0.15)(14,0.15)

\psset{arrowscale=2}
\psline[linecolor=darkbrn]{->}(4.6,1)(3.3,1)
\psline[linecolor=darkbrn]{->}(2.4,-1)(3.7,-1)
\pspolygon[linecolor=darkbrn](2.4,-1)(4.6,-1)(4.6,1)(2.4,1)

\rput(2,-0.4){$2$}
  \rput(5,-0.4){$5/2$}
  \rput(8,-0.4){$3$}
  \rput(11,-0.4){$7/2$}
  \rput(14,-0.4){$4$}

  \rput(-0.6,1){$1/N^2$}
  \rput(-0.8,-1){$-1/N^2$}
  \rput(2.4,-1.4){$2+\delta$}
  \rput(4.6,-1.4){$5/2-\delta'$}
  \rput(3.5,1.5){$C^+_1$}
  \rput(3.5,-0.5){$C^-_1$}
  \rput(1.85,0.5){$C_{1,L}$}
  \rput(5.15,0.5){$C_{1,R}$}

\psset{arrowscale=2}
\psline[linecolor=darkbrn]{->}(10.6,1)(9.3,1)
\psline[linecolor=darkbrn]{->}(8.4,-1)(9.7,-1)
\pspolygon[linecolor=darkbrn](8.4,-1)(10.6,-1)(10.6,1)(8.4,1)

  \rput(8.4,-1.4){$3+\delta$}
  \rput(10.6,-1.4){$7/2-\delta'$}
  \rput(9.5,1.5){$C^+_2$}
  \rput(9.5,-0.5){$C^-_2$}
  \rput(7.85,0.5){$C_{2,L}$}
  \rput(11.15,0.5){$C_{2,R}$}

\end{pspicture}
\caption{The rectangles $C_1$ and $C_2$}\label{c1c2}
\end{center}
\end{figure}
The next result shows that the integrals over $C_{1,L}$, $C_{1,R}$ are small.

\begin{prop} \label{cpcmb}
For $N$ greater than an absolute constant, we may choose  $\delta$, $\delta' \in [0.0061,0.01]$ so that
\begin{equation*}
\mathcal C_3(N,\sigma)  = \frac{1}{2N^{3/2}} \Re \int_{C_1^+ \cup C_1^-}
     \exp \bigl(N \cdot r_\mathcal C(z)\bigr) \frac{q_\mathcal C(z)}{2i\tan \bigl(\pi (2N/z-1)/2 \bigr)} \exp \bigl(v_\mathcal C(z;N,\sigma)\bigr)\, dz + O(e^{WN/2})
\end{equation*}
for $W = 0.05$ and an  implied constant depending only on $\sigma$.
\end{prop}
\begin{proof}
The proposition follows from \eqref{shir} if we can show $\int_{C_{1,L} \cup C_{1,R}} = O(e^{WN/2})$.
For $N$ large enough, we may choose $\delta$ and $\delta'$ so that $C_{1,L}$ and $C_{1,R}$ pass midway between the poles of $1/\tan \bigl(\pi (2N/z-1)/2 \bigr)$. Hence
\begin{equation} \label{sha}
    \frac{1}{\tan \bigl(\pi (2N/z-1)/2 \bigr)} \ll 1 \qquad (z \in C_{1,L} \cup C_{1,R}).
\end{equation}
The bound \eqref{efb3} from Theorem \ref{rzjb} implies
\begin{equation} \label{shb}
    q_\mathcal C(z) \exp \bigl(v_\mathcal C(z;N,\sigma)\bigr) \ll 1 \qquad (z \in C_{1,L} \cup C_{1,R}).
\end{equation}
Theorem \ref{rzjb} with $j=0$ also implies
\begin{equation*}
    \Re\bigl(r_\mathcal C(z) \bigr)  < \frac{1}{8\pi} \left(x \cl(2\pi x) + \frac{5\pi^2}{N^2} \right) \qquad (z \in C_{1,L} \cup C_{1,R}).
\end{equation*}
Note that
\begin{equation}\label{futcl}
    \cl(2\pi x) < 0.24 \quad \text{if} \quad 2 \lqs x \lqs 2.01, \qquad \cl(2\pi x) < 0.05 \quad \text{if} \quad 2.49 \lqs x \lqs 2.5.
\end{equation}
Therefore
\begin{equation} \label{klr}
    \Re\bigl(r_\mathcal C(z) \bigr)  < \frac{1}{8\pi} \left(2.01 \times 0.24 + \frac{5\pi^2}{N^2} \right) < 0.025 \qquad (z \in C_{1,L}, \ N \gqs 25)
\end{equation}
and we obtain \eqref{klr} for $z \in C_{1,R}$ in the same way. Consequently
\begin{equation} \label{shc}
    \exp\bigl(N \cdot r_\mathcal C(z) \bigr) \ll  \exp(0.025 N) \qquad (z \in C_{1,L} \cup C_{1,R}).
\end{equation}
The proposition now follows from the bounds \eqref{sha}, \eqref{shb} and \eqref{shc}.
\end{proof}

We have
\begin{equation}\label{tantan2}
\frac{1}{2i\tan(\pi (2N/z-1)/2)} = \begin{cases}
1/2+\sum_{j \lqs -1} (-1)^j e^{2\pi i j N/z} & \text{ \ if \ } \Im z >0 \\
-1/2-\sum_{j \gqs 1} (-1)^j e^{2\pi i j N/z} & \text{ \ if \ } \Im z <0
\end{cases}
\end{equation}
and therefore
\begin{align}
    \int_{C_1^+} & =  \sideset{}{'}{\sum}_{j \lqs 0} (-1)^j \int_{C_1^+}
    \exp \bigl(N [ r_\mathcal C(z)+2\pi i j/z] \bigr) q_\mathcal C(z)\exp \bigl(v_\mathcal C(z;N,\sigma)\bigr)\, dz, \label{gas+}\\
    \int_{C_1^-} & =  -\sideset{}{'}{\sum}_{j \gqs 0} (-1)^j \int_{C_1^-}
    \exp \bigl(N [ r_\mathcal C(z)+2\pi i j/z] \bigr) q_\mathcal C(z)\exp \bigl(v_\mathcal C(z;N,\sigma)\bigr)\, dz \label{gas-}
\end{align}
where $\sum'$ indicates the $j=0$ term is taken with a $1/2$  factor.
The  terms with $j=0$, $-1$ are the largest:

\begin{prop} \label{lgst}
For $W = 0.05$ and an  implied constant depending only on $\sigma$
\begin{equation} \label{gas}
\mathcal C_3(N,\sigma)  = \frac{-1}{2N^{3/2}} \sum_{j=0,-1} (-1)^j \Re \int_{2.01}^{2.49}
     \exp \bigl(N [ r_\mathcal C(z)+2\pi i j/z] \bigr) q_\mathcal C(z)\exp \bigl(v_\mathcal C(z;N,\sigma)\bigr) \, dz + O(e^{WN/2}).
\end{equation}
\end{prop}
\begin{proof}
As in \cite[Sect. 4.5]{OS1}, the total contribution to \eqref{gas+}, \eqref{gas-} for all $j$ with $|j|>N^2$ can be shown to be $O(N)$.
Let $D_1^{+}$ be the three lines which, when added to $C_1^+$, make a rectangle with top side having imaginary part $1$. Orient the path $D_1^{+}$ so that it has the same starting and ending points as $C_1^+$. Since the integrand is holomorphic we see that $\int_{C_1^+ } = \int_{D_1^+ }$. For integers $j$ with $-N^2 \lqs j<0$ we consider
\begin{equation} \label{nnj1b}
      \int_{D_1^+} \exp \bigl(N [ r_\mathcal C(z)+2\pi i j/z] \bigr) q_\mathcal C(z)\exp \bigl(v_\mathcal C(z;N,\sigma)\bigr) \, dz .
\end{equation}
We have $q_\mathcal C(z)\exp \left(v_\mathcal C(N,z)\right) \ll 1$ for $z \in D_1^+$ by Theorem \ref{rzjb}. On the vertical sides of $D_1^+$ we have
\begin{equation*}
    \Re\left(r_\mathcal C(z) +\frac{2\pi i j}{z}\right) < \frac{x \cl(2\pi x)}{8\pi} <0.02
\end{equation*}
by Theorem \ref{rzjb} and \eqref{futcl} if $j\lqs -2$. On the horizontal side of $D_1^+$, with $y=1$, Theorem \ref{rzjb} implies
\begin{equation*}
    \Re\left(r_\mathcal C\left(z \right) +\frac{2\pi i j}{z}\right) \lqs \frac{1}{2\pi |z|^2} \left(2.5 \cl(\pi/3) +\pi^2 \left[\frac 13 +4(j+1) \right]  \right) <0
\end{equation*}
if $j\lqs -2$. Hence, for each integer $j$ with $-N^2 \lqs j \lqs -2$, \eqref{nnj1b} is $O\bigl(\exp(0.02 N)\bigr)$. In a similar way, the terms in \eqref{gas-} for $1 \lqs j \lqs N^2$ are $O\bigl(\exp(0.02 N)\bigr)$. Moving the lines of integration from $C_1^-$ and $C_1^+$ to $[2.01,2.49]$ is valid with \eqref{shb}, \eqref{shc} and this completes the proof.
\end{proof}

A slightly more detailed argument shows that the $j=0$ term in \eqref{gas} is also $O(e^{WN/2})$:

\begin{prop} \label{lgst2}
For $W = 0.05$ and an  implied constant depending only on $\sigma$
\begin{equation} \label{gast}
 \int_{2.01}^{2.49}
     \exp \bigl(N \cdot r_\mathcal C(z) \bigr) q_\mathcal C(z)\exp \bigl(v_\mathcal C(z;N,\sigma)\bigr) \, dz = O(e^{WN/2}).
\end{equation}
\end{prop}
\begin{proof}
Change the path of integration to the lines joining $2.01$, $2.01-i$, $2.49-i$ and $2.49$.
The result follows if we can show $\Re (r_\mathcal C(z) ) \lqs W/2$ on these lines. For $y\lqs 0$, by \eqref{fdrb},
\begin{align*}
\Re \bigl(r_\mathcal C(z) \bigr) & = \pi y-\frac{y}{2\pi|z|^2}\left(\li(1)-\Re (\li (e^{-2\pi i z}))+ 8\pi^2 \right)  -\frac{x  \Im(\li (e^{-2\pi i z}))}{2\pi|z|^2}\\
& \lqs
\pi y-\frac{y}{2\pi|z|^2}\left(\li(1)-\Re (\li (e^{-2\pi i z}))+ 8\pi^2 \right)  +\frac{x \cl(2\pi x)}{2\pi|z|^2}
\end{align*}
using Lemma \ref{dil1}. Recalling \eqref{futcl} we obtain the following bounds on each segment:
\begin{itemize}
\item $x=2.01$, $-1\lqs y \lqs 0$. By Lemma \ref{dil1a} we have $-\Re (\li (e^{-2\pi i z})) \lqs 0$ so that
\begin{equation*}
    \Re \bigl(r_\mathcal C(z) \bigr) \lqs \pi y+\frac{1}{2\pi(x^2+y^2)}\left(-y(\li(1)+ 8\pi^2) +0.24 x \right) < 0.025.
\end{equation*}
\item $x=2.49$, $-1\lqs y \lqs 0$. By Lemma \ref{dil2} we have $-\Re (\li (e^{-2\pi i z})) \lqs \li(1)$ so that
\begin{equation*}
    \Re \bigl(r_\mathcal C(z) \bigr) \lqs \pi y+\frac{1}{2\pi(x^2+y^2)}\left(-y(2\li(1)+ 8\pi^2) + 0.05 x \right) < 0.01.
\end{equation*}
\item $2\lqs x \lqs 2.5$, $y= -1$. With Lemma \ref{dil2} again
\begin{equation*}
    \Re \bigl(r_\mathcal C(z) \bigr) \lqs \pi y+\frac{1}{2\pi(2^2+y^2)}\left(-y(2\li(1)+ 8\pi^2) + 2.5 \cl(\pi/3) \right) < 0. \qedhere
\end{equation*}
\end{itemize}
\end{proof}

Since
$
    p(z) = -(r_\mathcal C(z)-2\pi i /z)
$,
and recalling \eqref{cosmos}, we have therefore shown
\begin{equation}\label{cos2}
    \mathcal C_2(N,\sigma)  = \mathcal C_4(N,\sigma) + O(e^{WN/2})
\end{equation}
for $W=0.05$, an implied constant depending only on $\sigma$, and
\begin{equation}\label{newe}
    \mathcal C_4(N,\sigma)  := \frac{1}{2N^{3/2}}  \Re \int_{2.01}^{2.49}
     \exp \bigl(-N \cdot p(z) \bigr) q_\mathcal C(z)\exp \bigl(v_\mathcal C(z;N,\sigma)\bigr) \, dz.
\end{equation}

\subsection{A path through the saddle-point}
To apply the saddle-point method, Theorem \ref{sdle}, to $\mathcal C_4(N,\sigma)$ we first locate the unique solution to $p'(z)=0$ for $3/2<\Re(z)<5/2$ as
\begin{equation*}
    z_1:=2+\frac{\log \bigl(1-w(0,-2)\bigr)}{2\pi i} \approx 2.20541 + 0.345648 i
\end{equation*}
by Theorem \ref{disol}. Then we replace the path of integration $[2.01,2.49]$ in \eqref{newe} with one passing through $z_1$.

Let $v=\Im(z_1)/\Re(z_1) \approx 0.156728$ and $c=1+i v $. The path we take through the saddle-point $z_1$ is $\mathcal Q := \mathcal Q_1 \cup \mathcal Q_2 \cup \mathcal Q_3$, the polygonal path between the points $2.01$, $2.01 c$, $2.49 c$ and $2.49$ as shown in Figure \ref{pthb}.

\SpecialCoor
\psset{griddots=5,subgriddiv=0,gridlabels=0pt}
\psset{xunit=0.8cm, yunit=0.6cm}
\psset{linewidth=1pt}
\psset{dotsize=4pt 0,dotstyle=*}

\newrgbcolor{darkbrn}{0.4 0.2 0}

\begin{figure}[h]
\begin{center}
\begin{pspicture}(-1,-1)(8,4.5) 

\psline[linecolor=gray]{->}(0,-1)(0,4.3)
\psline[linecolor=gray]{->}(-1,0)(8,0)
\psline[linecolor=gray](-0.15,3.4884)(0.15,3.4884)
\psline[linecolor=gray](1.8,-0.15)(1.8,0.15)
\psline[linecolor=gray](7.2,-0.15)(7.2,0.15)
\psline[linecolor=gray](4.06,-0.15)(4.06,0.15)

\psset{arrowscale=2,arrowinset=0.5}
\psline[linecolor=darkbrn]{->}(2,0)(2,1.9)
\psline[linecolor=darkbrn]{->}(7,3.9)(7,1.6)
\psline[linecolor=darkbrn]{->}(5,3.62)(5.1,3.634)
\psline[linecolor=darkbrn](2,0)(2,3.2)(7,3.9)(7,0)

\rput(1.8,-0.6){$2$}
  \rput(7.2,-0.6){$5/2$}
  \rput(-1,3.4884){$0.346$}
  \rput(4.06,-0.6){$2.205$}
  \rput(1.4,1.55){$\mathcal Q_1$}
  \rput(7.6,1.9){$\mathcal Q_3$}
  \rput(5,2.9){$\mathcal Q_2$}
  \rput(4.06,4){$z_1$}

\psdots(4.06,3.4884)

\end{pspicture}
\caption{The path $\mathcal Q = \mathcal Q_1 \cup \mathcal Q_2 \cup \mathcal Q_3$ through $z_1$}\label{pthb}
\end{center}
\end{figure}
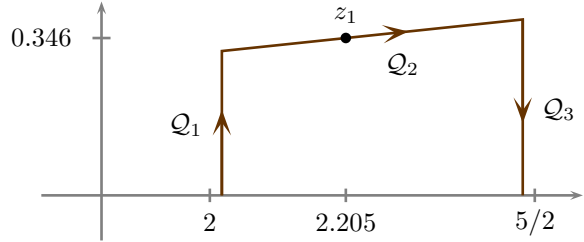
%
%
For Theorem \ref{sdle} we require the next result.
\begin{theorem} \label{sdleverc}
 For the path $\mathcal Q$ above, passing through the saddle point $z_1$, we have $\Re(p(z)-p(z_1))>0$ for $z\in \mathcal Q$ except at $z=z_1$.
 \end{theorem}

 Theorem \ref{sdleverc} seems apparent from Figure \ref{rpzb}. We  prove it by approximating $\Re(p(z))$ and its derivatives by the first terms in their series expansions and  reducing the issue to a finite computation. This method was used in \cite[Sect. 5.2]{OS1} and we repeat the results from there. To take into account that we are using an approximation to $z_1$, we give proofs valid in a range $0.15\lqs v \lqs 0.16$.


\SpecialCoor
\psset{griddots=5,subgriddiv=0,gridlabels=0pt}
\psset{xunit=0.7cm, yunit=0.8cm}
\psset{linewidth=1pt}
\psset{dotsize=4pt 0,dotstyle=*}

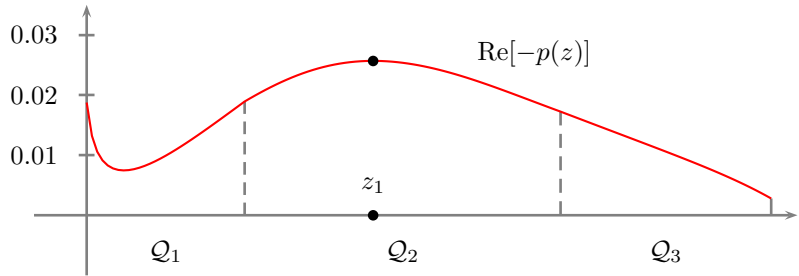
\begin{figure}[h]
\begin{center}
\begin{pspicture}(-1,-1)(14,3.5) 

\psline[linecolor=gray]{->}(0,-1)(0,3.5)
\psline[linecolor=gray]{->}(-1,0)(13.5,0)
\psline[linecolor=gray,linestyle=dashed](3,0)(3,1.89142)
\psline[linecolor=gray,linestyle=dashed](9,0)(9,1.723)
\psline[linecolor=gray,linestyle=dashed](13,0)(13,0.278306)

  \rput(-1,1){$0.01$}
  \rput(-1,3){$0.03$}
   \rput(-1,2){$0.02$}
  \rput(1.5,-0.6){$\mathcal Q_1$}
  \rput(11,-0.6){$\mathcal Q_3$}
  \rput(6,-0.6){$\mathcal Q_2$}

\savedata{\mydata}[
{{0., 1.8743}, {0.1, 1.31941}, {0.2, 1.05684}, {0.3, 0.912648}, {0.4,
  0.827755}, {0.5, 0.778529}, {0.6, 0.753293}, {0.7, 0.74538}, {0.8,
  0.750589}, {0.9, 0.766077}, {1., 0.789811}, {1.1, 0.820283}, {1.2,
  0.85633}, {1.3, 0.89704}, {1.4, 0.941676}, {1.5, 0.989637}, {1.6,
  1.04042}, {1.7, 1.09361}, {1.8, 1.14885}, {1.9, 1.20582}, {2.,
  1.26428}, {2.1, 1.32399}, {2.2, 1.38475}, {2.3, 1.44638}, {2.4,
  1.50874}, {2.5, 1.57169}, {2.6, 1.6351}, {2.7, 1.69887}, {2.8,
  1.7629}, {2.9, 1.82711}, {3., 1.89142}, {3.1, 1.94233}, {3.2,
  1.99193}, {3.3, 2.04008}, {3.4, 2.08665}, {3.5, 2.1315}, {3.6,
  2.17452}, {3.7, 2.2156}, {3.8, 2.25465}, {3.9, 2.29159}, {4.,
  2.32634}, {4.1, 2.35885}, {4.2, 2.38906}, {4.3, 2.41694}, {4.4,
  2.44245}, {4.5, 2.46557}, {4.6, 2.4863}, {4.7, 2.50462}, {4.8,
  2.52054}, {4.9, 2.53407}, {5., 2.54524}, {5.1, 2.55407}, {5.2,
  2.56059}, {5.3, 2.56484}, {5.4, 2.56686}, {5.5, 2.5667}, {5.6,
  2.56442}, {5.7, 2.56006}, {5.8, 2.5537}, {5.9, 2.54538}, {6.,
  2.53518}, {6.1, 2.52317}, {6.2, 2.50941}, {6.3, 2.49399}, {6.4,
  2.47696}, {6.5, 2.45841}, {6.6, 2.43841}, {6.7, 2.41705}, {6.8,
  2.39439}, {6.9, 2.37052}, {7., 2.34551}, {7.1, 2.31944}, {7.2,
  2.29239}, {7.3, 2.26444}, {7.4, 2.23567}, {7.5, 2.20615}, {7.6,
  2.17596}, {7.7, 2.14517}, {7.8, 2.11385}, {7.9, 2.08209}, {8.,
  2.04996}, {8.1, 2.01752}, {8.2, 1.98485}, {8.3, 1.95201}, {8.4,
  1.91907}, {8.5, 1.8861}, {8.6, 1.85317}, {8.7, 1.82033}, {8.8,
  1.78765}, {8.9, 1.75519}, {9., 1.723}, {9.1, 1.68965}, {9.2,
  1.65633}, {9.3, 1.62304}, {9.4, 1.58977}, {9.5, 1.55653}, {9.6,
  1.52331}, {9.7, 1.49012}, {9.8, 1.45696}, {9.9, 1.42381}, {10.,
  1.39068}, {10.1, 1.35756}, {10.2, 1.32444}, {10.3, 1.29131}, {10.4,
  1.25816}, {10.5, 1.22498}, {10.6, 1.19175}, {10.7, 1.15845}, {10.8,
  1.12506}, {10.9, 1.09155}, {11., 1.05791}, {11.1, 1.02409}, {11.2,
  0.990075}, {11.3, 0.955813}, {11.4, 0.921269}, {11.5,
  0.886396}, {11.6, 0.851145}, {11.7, 0.815463}, {11.8,
  0.779289}, {11.9, 0.742561}, {12., 0.70521}, {12.1,
  0.667162}, {12.2, 0.628336}, {12.3, 0.588648}, {12.4,
  0.548008}, {12.5, 0.506317}, {12.6, 0.463474}, {12.7,
  0.419369}, {12.8, 0.373888}, {12.9, 0.32691}, {13., 0.278306}}
]

\multirput(-0.15,1)(0,1){3}{\psline[linecolor=gray](0,0)(0.3,0)}

\dataplot[linecolor=red,linewidth=0.8pt,plotstyle=line]{\mydata}

  \rput(5.44259,0.5){$z_1$}

\psdots(5.44259,2.56706)(5.44259,0)

\rput(8.5,2.7){$\Re[-p(z)]$}

\end{pspicture}
\caption{Graph of $\Re[-p(z)]$ for $z \in \mathcal Q$}\label{rpzb}
\end{center}
\end{figure}


Generalizing to $p_d(z)$, we examine $\Re(p_d(z))$ for $z$ on the ray $z=ct$ for $c=1+ i v$ with $v>0$. We also
write
\begin{equation*}
    c = \rho e^{i \theta}  \qquad (0<\rho, \ 0<\theta <\pi/2).
\end{equation*}
For the second derivative we have
\begin{equation*}
    \frac{d^2}{dt^2}\Re [p_d(ct) ] = R_2(L;t)+R^*_2(L;t)
\end{equation*}
for
\begin{align*}
    R_2(L;t) & :=  -  \frac{\pi (24d+1)\sin \theta}{6\rho t^3} + \sum_{m=1}^{L-1}\Bigl( A_m(t) \cos(2\pi m t) + B_m(t) \sin(2\pi m t)\Bigr),\\
    A_m(t) & := e^{-2\pi m v t} \left(\frac{2}{m t^2}+\sin \theta\left(\frac{2\pi \rho}{t}+\frac{1}{m^2 \pi \rho t^3} \right) \right),\\
    B_m(t) & := e^{-2\pi m v t} \cos \theta \left(\frac{2\pi \rho}{t}-\frac{1}{m^2 \pi \rho t^3}  \right)
\end{align*}
and
\begin{equation*}
    |R^*_2(L;t)| \lqs E_2(L;t) := \frac{e^{-2\pi L v t}}{1-e^{-2\pi v t}}\left( \frac{  1}{\pi \rho L^2 t^3}
    + \frac{2}{L t^2} + \frac{2\pi \rho }{t}\right).
\end{equation*}
We see that $E_2(L;t)$ is a decreasing function of $L$ and $t$. We have $A_m(t)$ a positive and decreasing function of $t$. Also  $B_m(t)$ is a positive and decreasing function of $t$ when $t>\frac{\sqrt{3}}{\sqrt{2} \pi \rho m}$. The above formulas for $R_2(L;t)$ use \eqref{def0}, which is valid since $|e^{2\pi i z}| \lqs 1$ when $\Im(z) \gqs 0$. For a ray $z=ct$ with $\Im(c)<0$, the functional equation \eqref{ss} must be applied first and then similar formulas are found.

Let $v_1=0.15$ and $v_2=0.16$. Writing $\rho_1 e^{i \theta_1}=1+i v_1$ and $\rho_2 e^{i \theta_2}=1+i v_2$ we have
\begin{equation*}
    1<\rho_1\lqs \rho \lqs \rho_2, \quad 0< \theta_1 \lqs \theta \lqs \theta_2<\pi/2.
\end{equation*}
For $v$ in the interval $[v_1,v_2]$, we may bound $A_m(t)$, $B_m(t)$ and $E_2(L;t)$   from above and below by replacing $v$, $\rho$ and $\theta$ appropriately by $v_j$, $\rho_j$ and $\theta_j$, $j=1,2$. For example
\begin{equation*}
    0 < A^-_m(t)  \lqs A_m(t) \lqs A^+_m(t) \qquad (v \in [v_1,v_2])
\end{equation*}
with
\begin{align*}
    A^-_m(t) & := e^{-2\pi m v_2 t} \left(\frac{2}{m t^2}+\sin \theta_1\left(\frac{2\pi \rho_1}{t}+\frac{1}{m^2 \pi \rho_2 t^3} \right) \right),\\
    A^+_m(t) & := e^{-2\pi m v_1 t} \left(\frac{2}{m t^2}+\sin \theta_2\left(\frac{2\pi \rho_2}{t}+\frac{1}{m^2 \pi \rho_1 t^3} \right) \right)
\end{align*}
and similarly write $0< B^-_m(t)  \lqs B_m(t) \lqs B^+_m(t)$ and  $0< E^-_2(L;t) \lqs E_2(L;t) \lqs E^+_2(L;t)$.

\begin{lemma} \label{msy1c}
Let $c=1+i v$ with $0.15 \lqs v \lqs 0.16$. Then  $\frac{d^2}{dt^2}\Re [p(ct) ] >0$ for $t \in [2,2.35]$.
\end{lemma}
\begin{proof}
Break up $[2,2.35]$ into $n$ equal segments $[x_{j-1},x_j]$. Then
\begin{equation} \label{dabel}
    \frac{d^2}{dt^2}\Re [p(ct) ] \gqs \min_{1\lqs j \lqs n} \left(\left(\min_{t \in [x_{j-1},x_j]} R_2(L;t) \right) -E^+_2(L;x_{j-1})\right).
\end{equation}
Let $t=x^*_{j,m}$ correspond to the minimum value of $\cos(2\pi m t)$ for $t\in [x_{j-1},x_j]$ (so that $x^*_{j,m}$ equals $x_{j-1}$, $x_j$ or a local minimum $k/2m$ for $k$ odd). Similarly, let $t=x^{**}_{j,m}$ correspond to the minimum value of $\sin(2\pi m t)$ for $t\in [x_{j-1},x_j]$. Then
\begin{equation} \label{dabell}
    \min_{t \in [x_{j-1},x_j]} R_2(L;t) \gqs  -  \frac{\pi \sin \theta_2}{6\rho_1 x_{j-1}^3} + \sum_{m=1}^{L-1}\Bigl( A^-_m(x_j) \cos(2\pi m x^*_{j,m}) + B^-_m(x_j) \sin(2\pi m x^{**}_{j,m})\Bigr)
\end{equation}
where we must replace $A^-_m(x_j)$  in \eqref{dabell} by $A^+_m(x_{j-1})$  if $\cos(2\pi m x^*_{j,m})<0$ and replace $B^-_m(x_j)$ in \eqref{dabell} by $B^+_m(x_{j-1})$ if $\sin(2\pi m x^{**}_{j,m})<0$.

A computation using \eqref{dabel} and \eqref{dabell} with $n=10$ and $L=3$ for example shows $\frac{d^2}{dt^2}\Re [p(ct) ] > 0.09$.
\end{proof}

We may analyze the first derivative in a similar way. We have
\begin{equation*}
    \frac{d}{dt}\Re [p_d(ct) ] = R_1(L;t)+R^*_1(L;t)
\end{equation*}
for
\begin{align*}
    R_1(L;t) & :=   \frac{\pi (24d+1)\sin \theta}{12\rho t^2} + \sum_{m=1}^{L-1}\Bigl( -C_m(t) \cos(2\pi m t) + D_m(t) \sin(2\pi m t)\Bigr),\\
    C_m(t) & := e^{-2\pi m v t} \left(\frac{1}{m t}+\frac{\sin \theta}{m^2 2\pi \rho t^2} \right), \qquad D_m(t)  := e^{-2\pi m v t} \frac{\cos \theta}{m^2 2\pi \rho t^2}\\
\end{align*}
and
\begin{equation*}
    |R^*_1(L;t)| \lqs E_1(L;t) := \frac{e^{-2\pi L v t}}{1-e^{-2\pi v t}}\left( \frac{  1}{2\pi \rho L^2 t^2}
    + \frac{1}{L t} \right).
\end{equation*}
We see that $E_1(L;t)$ is a decreasing function of $L$ and $t$. Also $C_m(t)$ and $D_m(t)$ are positive and decreasing functions of $t$.

\begin{lemma} \label{msy2c}
Let $c=1+i v$ with $0.15 \lqs v \lqs 0.16$. Then  $\frac{d}{dt}\Re [p(ct) ] >0$ for $t \in [2.35, 2.5]$.
\end{lemma}
\begin{proof}
  Break $[2.35, 2.5]$ into $n$ equal segments and, as in the proof of Lemma \ref{msy1c}, bound $\frac{d}{dt}\Re [p(ct) ]$ from below on each piece. Taking $n=10$ and $L=3$ shows $\frac{d}{dt}\Re [p(ct) ] > 0.03$ for example.
\end{proof}

\begin{cor} \label{bflc}
Let $c=1+i v$ with $0.15 \lqs v \lqs 0.16$. There is a unique solution to $\frac{d}{dt}\Re [p(ct)]=0 $    for $t\in [2,2.5]$ that we label as $t_0$. We then have $\Re [p(ct)-p(c t_0) ] >0$ for $t\in [2,2.5]$ except at $t=t_0$.
\end{cor}
\begin{proof}
 Check that $\frac{d}{dt}\Re [p(ct) ]<0$ when $t=2$ and $\frac{d}{dt}\Re [p(ct) ]>0$ when $t=2.35$. By Lemma \ref{msy1c} we see that $\frac{d}{dt}\Re [p(ct) ]$ is strictly increasing for $t\in [2,2.35]$. It necessarily has a unique zero that we label $t_0$. By Lemma \ref{msy2c}, $\frac{d}{dt}\Re [p(ct) ]$ remains $>0$ for $t\in [2.35, 2.5]$ . Hence $\Re [p(ct)-p(c t_0)]$ is strictly decreasing on $[2,t_0)$ and strictly increasing on $(t_0,2.5]$ as required.
\end{proof}

\begin{prop} \label{bxlc}
For $0.15 \lqs v \lqs 0.16$  we have $\Re[-p(z)]<0.024$ for $z \in \mathcal Q_1 \cup \mathcal Q_3$.
\end{prop}
\begin{proof} We have $x$ fixed as $2.01$ on $\mathcal Q_1$ and $2.49$ on $\mathcal Q_3$.
Write
\begin{equation*}
    \Re[-p(z)] = \frac{f(y)+g(y)}{2\pi |z|^2}
\end{equation*}
for
\begin{equation*}
    f(y):=y\left(\li(1)-\Re(\li(e^{2\pi i z}))\right), \quad g(y)= x \Im (\li(e^{2\pi i z})).
\end{equation*}
 If $x=2.01$ or $2.49$ it follows from Lemma \ref{dil1} that $g(y)$ is positive and decreasing. Similarly, it follows from Lemma \ref{dil1a} that $f(y)$ is always positive and increasing for $y>0$.

For $z\in \mathcal Q_1$, so that $x=2.01$ and $0\lqs y \lqs Y:=2.01\times 0.16=0.3216$,
\begin{equation*}
    \Re[-p(z)] \lqs \begin{cases} (f(Y/3)+g(0))/(2\pi 2.01^2) \approx 0.0232 & \quad y \in [0,Y/3]\\
    (f(Y)+g(Y/3))/(2\pi( 2.01^2+(Y/3)^2) \approx 0.0226 & \quad y \in [Y/3,Y].
    \end{cases}
\end{equation*}
For $z\in \mathcal P_3$, so that $x=2.49$ and $0\lqs y \lqs Y:=2.49 \times 0.16 = 0.3984$,
\begin{equation*}
    \Re[-p(z)] \lqs (f(Y)+g(0))/(2\pi 2.49^2) \approx 0.021,  \qquad y \in [0,Y]. \qedhere
\end{equation*}
\end{proof}

\begin{proof}[\bf Proof of Theorem \ref{sdleverc}]
Let $v$ be given by $\Im(z_1)/\Re(z_1)$.   Then
\begin{equation*}
    \left.\frac{d}{dt}\Re [p(ct) ]\right|_{t=\Re(z_1)} = \Re [c p'(c\Re(z_1))] = \Re [c p'(z_1)]=0.
\end{equation*}
It follows from Corollary \ref{bflc} that $\Re [p(z)-p(z_1) ] >0$ for $z \in \mathcal Q_2$ and $z \neq z_1$. We also note that $\Re [-p(z_1) ] \approx 0.0256706$.

For $z \in \mathcal Q_1 \cup \mathcal Q_3$,  Proposition \ref{bxlc} implies $\Re [p(z)-p(z_1) ] > -0.024 +0.0256 >0$.
\end{proof}

\subsection{Applying the saddle-point method}

For $j \in \Z_{\gqs 0}$ put
\begin{equation} \label{uiz}
    u_{\sigma,j}(z):=\sum_{m_1+3m_2+5m_3+ \dots =j}\frac{(2\pi i \sigma z +g_1(z))^{m_1}}{m_1!}\frac{g_2(z)^{m_2}}{m_2!} \cdots \frac{g_j(z)^{m_j}}{m_j!},
\end{equation}
with $u_{\sigma,0}=1$. Recalling the definition of $g_\ell(z)$ in \eqref{grz}, we see that $u_{\sigma,j}(z)$ is holomorphic for $z \not\in \Z$. The proof of the next proposition uses Corollary \ref{acdcx}, see \cite[Sect. 5.3]{OS1}.

\begin{prop} \label{gasc}
For $2.01 \lqs \Re(z) \lqs 2.49$ and $|\Im(z)|\lqs 1$, say, there is a holomorphic function $\zeta_d(z;N,\sigma)$ of $z$ so that
\begin{equation*}
    \exp\bigl(v_{\mathcal C}(z;N,\sigma)\bigr) = \sum_{j=0}^{d-1} \frac{u_{\sigma,j}(z)}{N^j} + \zeta_d(z;N,\sigma) \quad \text{for} \quad \zeta_d(z;N,\sigma) = O\left(\frac{1}{N^d} \right)
\end{equation*}
with an implied constant depending only on $\sigma$ and  $d$ where $1 \lqs d \lqs 2L-1$ and $L=\lfloor 0.006 \pi e \cdot N/2 \rfloor$.
\end{prop}

We now have everything in place to get the asymptotic expansion of $\mathcal C_2(N,\sigma)$.

\begin{theorem} \label{c2se} With $c_{0}=-z_1 e^{-\pi i z_1}/2$ and  explicit  $c_{1}(\sigma)$, $c_{2}(\sigma), \dots $ depending on $\sigma \in \R$ we have
\begin{equation} \label{pres}
    \mathcal C_2(N,\sigma) = \Re\left[\frac{w(0,-2)^{-N}}{N^{2}} \left( c_{0}+\frac{c_{1}(\sigma)}{N}+ \dots +\frac{c_{m-1}(\sigma)}{N^{m-1}}\right)\right] + O\left(\frac{|w(0,-2)|^{-N}}{N^{m+2}}\right)
\end{equation}
for an implied constant depending only on  $\sigma$ and $m$.
\end{theorem}
\begin{proof}
Recall from \eqref{pzlogw} that
$e^{p(z_1)} = w(0,-2)$.
Proposition \ref{gasc} implies
\begin{equation}\label{umand}
    \mathcal C_4(N,\sigma)  = \Re\Biggl[ \sum_{j=0}^{d-1} \frac{1}{2N^{3/2+j}}  \int_\mathcal Q e^{-N \cdot p(z)} \cdot q_{\mathcal C}(z) \cdot u_{\sigma,j}(z) \, dz + \frac{1}{2N^{3/2}}  \int_\mathcal Q e^{-N \cdot p(z)} \cdot q_{\mathcal C}(z) \cdot \zeta_d(z;N,\sigma) \, dz \Biggr]
\end{equation}
where the last term in \eqref{umand} is
\begin{equation*}
    \ll \frac{1}{N^{3/2}}  \int_\mathcal Q \left|e^{-N \cdot p(z)}\right| \cdot 1 \cdot \frac{1}{N^d} \, dz
    \ll \frac{1}{N^{d+3/2}} e^{-N \Re(p(z_1))} = \frac{|w(0,-2)|^{-N}}{N^{d+3/2}}
\end{equation*}
by Theorem \ref{sdleverc}, \eqref{efb3} and Proposition \ref{gasc}.
Applying Theorem \ref{sdle} to each integral  in the first part of \eqref{umand}
 we obtain
\begin{equation} \label{wmand}
    \int_\mathcal Q e^{-N \cdot p(z)} \cdot q_{\mathcal C}(z) \cdot u_{\sigma,j}(z) \, dz = 2 e^{-N p(z_1)}\left(\sum_{s=0}^{S-1} \G(s+1/2)\frac{a_{2s}(q_{\mathcal C} \cdot u_{\sigma,j})}{N^{s+1/2}}+O\left( \frac{1}{N^{S+1/2}}\right) \right).
\end{equation}
The error term in \eqref{wmand} corresponds to an error for $\mathcal C_4(N,\sigma)$ of size $O(|w(0,-2)|^{-N}/N^{s+j+2})$.
We choose $S=d$ so that this error  is less than $O(|w(0,-2)|^{-N}/N^{d+3/2})$ for all $j \gqs 0$.
Therefore
\begin{align*}
    \mathcal C_4(N,\sigma) & = \Re \left[
    \sum_{j=0}^{d-1} \frac{1}{N^{j+3/2}}   e^{-N \cdot p(z_1)} \sum_{s=0}^{d-1} \frac{\G(s+1/2) a_{2s}(q_{\mathcal C} \cdot u_{\sigma,j})}{N^{s+1/2}}
    \right]+ O\left( \frac{|w(0,-2)|^{-N}}{N^{d+3/2}}\right) \\
    & = \Re \left[  w(0,-2)^{-N}
    \sum_{t=0}^{2d-2} \frac{1}{N^{t+2}}    \sum_{s=0}^{\min(t,d-1)} \G(s+1/2) a_{2s}(q_{\mathcal C} \cdot u_{\sigma, t-s})
    \right]+ O\left( \frac{|w(0,-2)|^{-N}}{N^{d+3/2}}\right) \\
    & = \Re \left[  w(0,-2)^{-N}
    \sum_{t=0}^{d-2} \frac{1}{N^{t+2}}    \sum_{s=0}^{t} \G(s+1/2) a_{2s}(q_{\mathcal C} \cdot u_{\sigma, t-s})
    \right]+ O\left( \frac{|w(0,-2)|^{-N}}{N^{d+1}}\right).
\end{align*}
Hence, recalling \eqref{cos2} and with
\begin{equation} \label{ctys}
    c_t(\sigma):=   \sum_{s=0}^t \G(s+1/2) a_{2s}(q_{\mathcal C} \cdot u_{\sigma,t-s}),
\end{equation}
we obtain \eqref{pres} in the statement of the theorem.

Use the formula \eqref{a2sb} for $a_0$   to get
\begin{equation*}
    c_0(\sigma)=\G(1/2)a_0(q_{\mathcal C} \cdot u_{\sigma,0})= \sqrt{\pi} \frac{\omega}{2(\omega^2 p_0)^{1/2}} q_0
\end{equation*}
which is independent of $\sigma$. The terms $p_0$ and $q_0$ are defined in \eqref{psp}, \eqref{psq}. Using the identity
\begin{equation} \label{pzzz22}
    p''(z)=-\frac 1z \left(2 p'(z) + \frac{2\pi i \cdot e^{2\pi i z}}{1-e^{2\pi i z}} \right)
\end{equation}
we obtain
\begin{equation} \label{oad2}
    p_0  =p''(z_1)/2 = \frac{-\pi i e^{2\pi i z_1}}{z_1 w(0,-2)},\qquad
    q^2_0  =q_{\mathcal C}(z_1)^2 = \frac{-i z_1}{w(0,-2)}.
\end{equation}
Therefore
\begin{equation*}
    c_0^2 = \frac{\pi q_0^2}{4p_0} = \frac{z_1^2}{4 e^{2\pi i z_1}}.
\end{equation*}
We may take $\omega=z_1$ since the path $\mathcal Q_2$ is a segment of the ray from the origin through $z_1$. A numerical check then gives us the correct square root:
\begin{equation*} \label{a0q}
    c_0 = \sqrt{\pi}\frac{\omega}{2(\omega^2 p_0)^{1/2}} q_0 = -\frac{z_1}{2 e^{\pi i z_1}}. \qedhere
\end{equation*}
\end{proof}

For example, Table \ref{c2n1} compares both sides of \eqref{pres} in Theorem \ref{c2se} with $\sigma=1$ and some different values of $m$ and $N$. For other values of $\sigma$ we get similar agreement.
\begin{table}[h]
\begin{center}
\begin{tabular}{c|cccc|c}
$N$ & $m=1$ & $m=2$ & $m=3$ & $m=4$ & $\mathcal C_2(N,1)$  \\ \hline
$800$ & $293.204$ &  $301.757$ & $303.016$ &  $303.119$ &  $303.112$ \\
$1000$ & $-263123.$ &  $-261461.$ & $-261486.$ &  $-261493.$ &  $-261493.$
\end{tabular}
\caption{Theorem \ref{c2se}'s approximations to $\mathcal C_2(N,1)$.} \label{c2n1}
\end{center}
\end{table}

\section{The asymptotic behavior of $\mathcal C_2^*(N,\sigma)$}

We find the asymptotic expansion of
\begin{equation*}
\mathcal C_2^*(N,\sigma)   = 2 \Re \sum_{k \text{ odd} \ : \ 2N/k \in [3,4)}  Q_{2k\sigma}(N),
\end{equation*}
the second component of $\mathcal C_1(N,\sigma)$, in this section.

\subsection{Approximating the sine product}

From \eqref{azi}, $Q_{2k\sigma}(N)$ contains the sine product $\spr{2/k}{m}$ for $m=N-k$ and $k/2<m<k$. The next result expresses this product in terms of a new variable $a$.

\begin{prop} \label{apb} Let $k$ be an odd positive integer.
Write $m=a +(k-1)/2$ for $1 \lqs a \lqs (k-1)/2$. Then
\begin{equation*}
    \spr{2/k}{m} = \frac{(-1)^a}{\sqrt{k}}\frac{\spr{1/k}{2a}}{\spr{2/k}{a}} \qquad \text{for} \quad k/2<m<k.
\end{equation*}
\end{prop}
\begin{proof}
The formula
\begin{equation*}
    \spr{h/k}{k-1} = (-1)^{(h-1)(k-1)/2}\frac 1k
\end{equation*}
from \cite[Sect. 2.2]{OS1} implies $\spr{2/k}{k-1} =(-1)^{(k-1)/2}/k$ and therefore, by symmetry,
\begin{equation*}
    \spr{2/k}{(k-1)/2} = \frac{1}{\sqrt{k}}.
\end{equation*}
Hence
\begin{align*}
  \spr{2/k}{m} & =  \spr{2/k}{(k-1)/2} \prod_{j=1}^a \frac{1}{2 \sin(\pi(j+(k-1)/2)2/k)} \\
  & =  \frac 1{\sqrt{k}} \prod_{j=1}^a \frac{1}{2 \sin(\pi(2j-1)/k + \pi)} \\
  & =  \frac {(-1)^a}{\sqrt{k}} \prod_{j=1}^a \frac{1}{2 \sin(\pi(2j-1)/k)}
\end{align*}
and the result follows.
\end{proof}

In this subsection we  define $z=z(N,k):=2(N+1/2)/k$ because of \eqref{cvn} below. The next result is the sine product approximation  we need here.
\begin{theorem} \label{cplxb2s}
Fix $W>0$. Let $\Delta$ be in the range $0.0048 \lqs \Delta \lqs 0.0079$ and set $\alpha = \Delta \pi e$.  Suppose $\delta$ and $\delta'$ satisfy
\begin{equation*}
    \frac{\Delta}{1-\Delta} < \delta \lqs \frac{1}{e}, \ 0<\delta' \lqs \frac{1}{e} \quad \text{ and  } \quad \delta \log 1/\delta, \ \ \delta' \log 1/\delta' \lqs W.
\end{equation*}
Then for all $N \gqs 3 \cdot R_\Delta$  we have
\begin{equation}\label{ott1b2s}
    \spr{2/k}{N-k} = O\left(e^{WN/3}\right) \quad \text{ for } \quad  z \in \left[3, \ 3+\delta \right] \cup  \left[7/2 -\delta' , \ 4 \right)
\end{equation}
and also for $ z \in  \left(3+\delta, \ 7/2 - \delta' \right)$
\begin{multline} \label{ott2b2s}
    \spr{2/k}{N-k} =  \frac{(-1)^{N+(k+1)/2}}{\sqrt{2k}}  \exp\left(\frac {N+1/2}{2\pi z} \cl(2\pi  z) \right)  \\
      \times \exp\left(\sum_{\ell=1}^{L-1} \frac{g_{\ell}(z)}{(2(N+1/2))^{2\ell-1}}-\sum_{\ell=1}^{L^*-1} \frac{g_{\ell}(z)}{(N+1/2)^{2\ell-1}} \right) +O\left(e^{WN/3}\right)
  \end{multline}
with $L=\lfloor \alpha \cdot 2N/3\rfloor$ and $L^*=\lfloor \alpha \cdot N/3\rfloor$. The implied constants in \eqref{ott1b2s}, \eqref{ott2b2s} are absolute.
\end{theorem}
\begin{proof}
We have $2N/k \in [3,4)$ so that $N/2 < k \lqs 2N/3$. For $m=N-k$ this corresponds to $k/2<m<k$. In terms of $a$ this means
\begin{equation}\label{cvn}
1 \lqs a<k/2, \qquad 2a/k=z-3.
\end{equation}

For $L=1$, Proposition \ref{prrh} implies
\begin{align*}
    \spr{1/k}{2a} & = \left( \frac{1}{2k \sin(2\pi a/k)}\right)^{1/2} \exp\left(\frac{k}{2\pi}\cl(4\pi a/k) \right)\exp\left( -T_1(2a,1/k) \right),\\
    \spr{2/k}{a} & = \left( \frac{1}{k \sin(2\pi a/k)}\right)^{1/2} \exp\left(\frac{k}{4\pi}\cl(4\pi a/k) \right)\exp\left( -T_1(a,2/k) \right)
\end{align*}
so that
\begin{align*}
    \frac{\spr{1/k}{2a}}{\spr{2/k}{a}} & = \frac{1}{2^{1/2}} \exp\left(\frac{k}{4\pi}\cl(4\pi a/k) \right)\exp\left(-T_1(2a,1/k) + T_1(a,2/k) \right)\\
     & = \left( \frac{k}{2} \sin(2\pi a/k)\right)^{1/2} \exp\left(-T_1(2a,1/k) + 2 T_1(a,2/k) \right) \cdot \spr{2/k}{a}.
\end{align*}
Therefore, employing Lemma \ref{melrh2},
\begin{equation} \label{oo1}
    \frac{\spr{1/k}{2a}}{\spr{2/k}{a}} \ll k^{1/2} \spr{2/k}{a} \qquad (1 \lqs a<k/2)
\end{equation}
with an absolute implied constant. A similar argument proves
\begin{equation} \label{oo2}
    \frac{\spr{1/k}{2a}}{\spr{2/k}{a}} \ll k^{1/4} \left(\spr{1/k}{2a}\right)^{1/2} \qquad (1 \lqs a<k/2).
\end{equation}
Then using Proposition \ref{gprop} to bound $\spr{2/k}{a}$ on the right of \eqref{oo1} and noting that $k \lqs 2N/3$ proves \eqref{ott1b2s}.

For positive integers $L_1$, $L_2$, Proposition \ref{prrh} implies
\begin{multline} \label{eva2}
\frac{\spr{1/k}{2a}}{\spr{2/k}{a}} = \frac{1}{\sqrt 2}\exp\left(\frac{k}{4\pi}\cl(4\pi a/k) \right)
 \exp\left(-\sum_{\ell=1}^{L_1-1} \frac{B_{2\ell}}{(2\ell)!} \left( \frac{\pi }{k}\right)^{2\ell-1} \cot^{(2\ell-2)}\left(\frac{2a\pi}{k}\right)\right) \\
\times\exp\left(\sum_{\ell=1}^{L_2-1} \frac{B_{2\ell}}{(2\ell)!} \left( \frac{2\pi }{k}\right)^{2\ell-1} \cot^{(2\ell-2)}\left(\frac{2a\pi}{k}\right)\right) \exp\left(-T_{L_1}(2a,1/k)+T_{L_2}(a,2/k) \right).
\end{multline}
Use Proposition \ref{dela} with $h=2$, $m=a$ and $s=2N/3$ to show that, for $\Delta N/3 \lqs a \lqs k/4$,
\begin{align}
     \spr{2/k}{a} T_{L_2}(a,2/k)  & \ll e^{WN/3} \label{tlmab2s}\\
     T_{L_2}(a,2/k)  & \ll 1  \label{tlmbb2s}
\end{align}
with absolute implied constants, $L_2:=\lfloor \pi e \Delta \cdot N/3 \rfloor$ and $N \gqs 3 \cdot R_\Delta$. The above inequality \eqref{tlmab2s} is valid with $\spr{2/k}{a} $ replaced by $\spr{1/k}{2a}  \big/ \spr{2/k}{a}$ using \eqref{oo1}:
\begin{equation}\label{tlmab2sy}
    \left(\spr{1/k}{2a} \big/ \spr{2/k}{a}\right) T_{L_2}(a,2/k)  \ll N^{1/2} e^{WN/3}.
\end{equation}

Use Proposition \ref{dela} with $h=1$, $m=2a$ and $s=2N/3$ to show that, also for $\Delta N/3 \lqs a \lqs k/4$,
\begin{align}
     \spr{1/k}{2a} T_{L_1}(2a,1/k)  & \ll e^{2WN/3} \label{tlmab2s1}\\
     T_{L_1}(2a,1/k)  & \ll 1  \label{tlmbb2s1}
\end{align}
with absolute implied constants, $L_1:=\lfloor \pi e \Delta \cdot 2N/3 \rfloor$  and $N \gqs 3 \cdot R_\Delta/2$. Taking square roots of both sides of inequality \eqref{tlmab2s1} and using \eqref{oo2} and that $\left| T_{L_1}(2a,1/k) \right| \ll \left| T_{L_1}(2a,1/k) \right|^{1/2}$ shows
\begin{equation}\label{tlmab2s1x}
    \left(\spr{1/k}{2a} \big/ \spr{2/k}{a}\right) T_{L_1}(2a,1/k)  \ll N^{1/4} e^{WN/3}.
\end{equation}

 With the inequalities \eqref{tlmab2s} - \eqref{tlmab2s1x} established, the arguments of Proposition \ref{hard} now go through, applied to \eqref{eva2}. This allows us to remove the factor $\exp\left(-T_{L_1}(2a,1/k)+T_{L_2}(a,2/k)\right)$ in \eqref{eva2} at the expense of adding an $O(e^{WN/3})$ error. The interval $\Delta N/3 \lqs a \lqs k/4$ corresponds to
\begin{equation*}
    3 + \frac{\Delta}{1-\Delta/3} \lqs z \lqs \frac 72
\end{equation*}
so we require
\begin{equation}\label{dcv2}
    \frac{\Delta}{1-\Delta/3} < \delta.
\end{equation}
The inequality \eqref{dcv2} is equivalent to $1/\Delta - 1/\delta \gqs 1/3$. Since our assumption
$\Delta/(1-\Delta)<\delta$  is equivalent to $1/\Delta - 1/\delta >1$, we have that \eqref{dcv2} is true. This completes the proof of \eqref{ott2b2s}.
\end{proof}

We rewrite \eqref{azi} as
\begin{equation*}
    Q_{2k\sigma}(N)=\frac {e^{-\pi i/4}}{k^2} \exp\left( (N+1/2)\frac{-\pi i}2 (z-5+2/z)\right) \exp\left(\frac{\pi i(16\sigma+1) z}{8(N+1/2)} \right)\spr{2/k}{N-k}
\end{equation*}
and combine with \eqref{ott2b2s} from Theorem \ref{cplxb2s} as follows.
With \eqref{dli} for $m=3$ we obtain
\begin{equation*}
    \cl(2\pi  z) = -i \li(e^{2\pi i z})+i \pi^2(z^2-7z+73/6) \qquad (3<z<4).
\end{equation*}
Hence
\begin{equation*}
    \frac {\cl(2\pi  z)}{2\pi z} -\frac{\pi i}2 (z-5+2/z) = -\pi i  +\frac 1{2\pi i z} \Bigl[ \li(e^{2\pi i z}) -\li(1) -10\pi^2 \Bigr].
\end{equation*}
Define the following  functions
\begin{align*}
  r_\mathcal C^*(z) &:= \frac 1{2\pi i z} \Bigl[\li(e^{2\pi i z}) -\li(1) -10\pi^2 \Bigr],\\
  q_\mathcal C^*(z) &:= e^{-3\pi i/4} \sqrt{z},\\
  v_\mathcal C^*(z;N,\sigma) &:= \frac{\pi i(16\sigma+1) z}{8(N+1/2)} + \sum_{\ell=1}^{L-1} \frac{g_{\ell}(z)}{(2(N+1/2))^{2\ell-1}}-\sum_{\ell=1}^{L^*-1} \frac{g_{\ell}(z)}{(N+1/2)^{2\ell-1}}
\end{align*}
for $L=\lfloor \alpha \cdot 2N/3\rfloor$ and $L^*=\lfloor \alpha \cdot N/3\rfloor$. Set
\begin{equation*}
    \mathcal C_3^*(N,\sigma)  := \frac{1}{(N+1/2)^{1/2}} \Re \sum_{k \text{ odd} \ : \ z \in  (3+\delta, 7/2 -\delta')}
    \frac{(-1)^{(k+1)/2}}{k^2} \exp \bigl((N+1/2) r^*_\mathcal C(z) \bigr) q^*_\mathcal C(z) \exp \bigl(v_\mathcal C^*(z;N,\sigma)\bigr) .
\end{equation*}
It follows from Theorem \ref{cplxb2s} that
\begin{equation}\label{coss}
\mathcal C_2^*(N,\sigma) = \mathcal C_3^*(N,\sigma) +O(e^{WN/3}).
\end{equation}

\subsection{Expressing $\mathcal C_3^*(N,\sigma)$ as an integral}
Similarly to Proposition \ref{addx} we have
\begin{prop} \label{addy}
Suppose $5/2 \lqs \Re(z) \lqs 7/2$ and $|z-3| \gqs \varepsilon >0$ and assume
$
    \max\{1+\frac{1}{\varepsilon}, \ 16 \} \ < \frac{\pi e}{\alpha}
$.
Then for $d \gqs 2$,
\begin{equation} \label{wey}
    \sum_{\ell=d}^{L-1} \frac{g_{\ell}(z)}{(2(N+1/2))^{2\ell-1}}-\sum_{\ell=d}^{L^*-1} \frac{g_{\ell}(z)}{(N+1/2)^{2\ell-1}} \ll \frac{1}{N^{2d-1}}e^{-\pi|y|}
\end{equation}
where $L=\lfloor \alpha \cdot 2N/3\rfloor$, $L^*=\lfloor \alpha \cdot N/3\rfloor$ and the implied constant depends only on $\varepsilon$, $\alpha$ and $d$.
\end{prop}

Fixing the choice of constants in \eqref{fix} and with $\varepsilon = 0.0061$
and
\begin{equation} \label{gazi}
    g_{\mathcal C, \ell}(z):=g_\ell(z)(2^{-(2\ell-1)}-1)
\end{equation}
 we obtain:
\begin{cor} \label{acdcy}
With $\delta, \delta' \in [0.0061, 0.01]$ and  $z \in \C$ such that $3+\delta \lqs \Re(z) \lqs 7/2-\delta'$  we have
\begin{equation*}
    v_\mathcal C^*(z;N,\sigma) = \frac{\pi i(16\sigma+1) z}{8(N+1/2)} +\sum_{\ell=1}^{d-1} \frac{g_{\mathcal C, \ell}(z)}{(N+1/2)^{2\ell-1}} + O\left( \frac 1{N^{2d-1}}\right)
\end{equation*}
 for $2 \lqs d\lqs L^*=\lfloor 0.006 \pi e \cdot N/3 \rfloor$ and an implied constant depending only on  $d$.
\end{cor}

Next,
\begin{align}
     r_\mathcal C^* \left(z \right) +\frac{2\pi i}{z}(j-1/2) & = \frac 1{2\pi i z} \Bigl[-\li(1)+ \li(e^{2\pi i z}) -4\pi^2(j+2) \Bigr], \notag \\
     r_\mathcal C^* \left(z \right) +\frac{2\pi i}{z}(j+1/2) & = \frac 1{2\pi i z} \Bigl[\li(1)- \li(e^{-2\pi i z}) -4\pi^2 (j-3) \Bigr]-  \pi i(2z-7) \label{tip}
\end{align}
where \eqref{tip} follows from \eqref{ss} when $3<\Re(z)<4$.
Then with a similar proof to  Theorem \ref{rzjb} we have
\begin{theorem}  \label{rzjbs}
The functions $r_\mathcal C^* (z)$, $q_\mathcal C^*(z)$ and $v_\mathcal C^*(z;N,\sigma)$ are holomorphic for $3<\Re(z)<7/2$. In this strip
\begin{alignat}{2}
    \Re\left(r_\mathcal C^*\left(z \right) +\frac{2\pi i }{z}(j-1/2)\right) & \lqs \frac{1}{2\pi |z|^2} \left(x \cl(2\pi x) +\pi^2 |y|\left[\frac 13 +4(j+2) \right]  \right) \qquad & & (y\gqs 0) \label{efbs1}\\
   \Re\left(r_\mathcal C^* \left(z \right) +\frac{2\pi i }{z}(j+1/2)\right) & \lqs \frac{1}{2\pi |z|^2} \left(x \cl(2\pi x) +\pi^2 |y|\left[\frac 13 -4(j+3/2) \right] \right) \qquad & & (y\lqs 0) \label{efbs2}
\end{alignat}
for $j\in \R$.
Also, in the box with $3+\delta \lqs \Re(z) \lqs 7/2-\delta'$ and $-1\lqs \Im(z) \lqs 1$,
\begin{equation} \label{efbs3}
    q_\mathcal C^*(z), \quad \exp \bigl(v_\mathcal C^* (z;N,\sigma)\bigr) \ll 1
\end{equation}
for an  implied constant depending only on $\sigma \in \R$.
\end{theorem}

By the calculus of residues,
\begin{equation*}
\sum_{a \lqs k \lqs b, \ k \text{ odd}} (-1)^{(k+1)/2} \varphi(k) =  \frac 12\int_C \frac{\varphi(z)}{2i\cos(\pi z/2)} \, dz
\end{equation*}
for $\varphi(z)$  a holomorphic function and $C$ a positively oriented closed contour surrounding the interval $[a,b]$ and not surrounding any integers outside this interval.
Hence
\begin{equation} \label{coscos}
    \sum_{a \lqs k \lqs b, \ k \text{ odd}} \frac{(-1)^{(k+1)/2}}{k^2} \varphi(2(N+1/2)/k) =  \frac {-1}{4(N+1/2)}\int \frac{\varphi(z)}{2i\cos(\pi (N+1/2)/z)} \, dz,
\end{equation}
for $C$ now surrounding $\{2(N+1/2)/k \ | \ a\lqs k \lqs b\}$ with $a>0$. Therefore
\begin{equation} \label{5the}
    \mathcal C_3^*(N,\sigma)  = \frac{-1}{4(N+1/2)^{3/2}} \Re \int_{C_2}
    \exp \bigl((N+1/2) r^*_\mathcal C(z) \bigr) \frac{ q^*_\mathcal C(z)}{2i\cos(\pi (N+1/2)/z)} \exp \bigl(v_\mathcal C^*(z;N,\sigma)\bigr) \, dz
\end{equation}
where $C_2$ is the positively oriented rectangle with horizontal sides $C_2^+$, $C_2^-$ having imaginary parts $1/N^2$, $-1/N^2$ and vertical sides $C_{2,L}$, $C_{2,R}$ having real parts $3+\delta$ and $7/2-\delta'$ respectively, as shown in Figure \ref{c1c2}.

Arguing as in Proposition \ref{cpcmb} proves the contribution to \eqref{5the} from  integrating over the vertical sides $C_{2,L}$, $C_{2,R}$ is $O(e^{0.016N})$.
We have
\begin{equation}\label{coscos2}
\frac{1}{2i\cos(\pi(N+1/2)/z)} = -i \times \begin{cases}
\sum_{j \lqs 0} (-1)^j \exp\left(\frac{2\pi i}{z}(N+1/2)(j-1/2)\right) & \text{ \ if } \Im z >0 \\
\sum_{j \gqs 0} (-1)^j \exp\left(\frac{2\pi i}{z}(N+1/2)(j+1/2)\right) & \text{ \ if } \Im z <0.
\end{cases}
\end{equation}
Therefore
\begin{multline*}
    -4(N+1/2)^{3/2}\mathcal C_3^*(N,\sigma) \\
     = \sum_{j \lqs 0} (-1)^j \Im \int_{C_2^+}  \exp  \left((N+1/2) \left[r^*_\mathcal C(z) +\frac{2\pi i}{z}(j-1/2)\right] \right) q^*_\mathcal C(z) \exp \bigl(v_\mathcal C^*(z;N,\sigma)\bigr) \, dz\\
    +  \sum_{j \gqs 0} (-1)^j \Im \int_{C_2^-}  \exp \left((N+1/2) \left[r^*_\mathcal C(z) +\frac{2\pi i}{z}(j+1/2)\right] \right) q^*_\mathcal C(z) \exp \bigl(v_\mathcal C^*(z;N,\sigma)\bigr) \, dz +O(e^{0.016N}).
\end{multline*}
A similar proof to Proposition \ref{lgst}'s, employing Theorem \ref{rzjbs}, shows that the total size of all but the $j=-1$, $-2$ terms  above
is $O(e^{0.013N})$. Let $d=j+2$ and we see
$
    p_d(z) = -(r^*_\mathcal C(z)+2\pi i (j-1/2)/z)
$
so that
\begin{multline} \label{murrs}
    4(N+1/2)^{3/2}\mathcal C_3^*(N,\sigma) \\
     =  \sum_{d=0,1} (-1)^d  \Im \int_{3.01}^{3.49}  \exp\bigl(-(N+1/2) p_d(z) \bigr) q^*_\mathcal C(z) \exp \bigl(v_\mathcal C^*(z;N,\sigma)\bigr) \, dz
    +O(e^{0.016N}).
\end{multline}

\subsection{Paths through the saddle-points}
We treat the $d=0$ case of \eqref{murrs} first. The unique solution to $p'(z)=0$ for $5/2<\Re(z)<7/2$ is
\begin{equation*}
    z_2:=3+\frac{\log \bigl(1-w(0,-3)\bigr)}{2\pi i} \approx 3.21625 + 0.402898 i
\end{equation*}
by Theorem \ref{disol}.
Let $v=\Im(z_2)/\Re(z_2) \approx 0.125269$ and $c=1+i v $. The path we take through the saddle point $z_2$ is $\mathcal R := \mathcal R_1 \cup \mathcal R_2 \cup \mathcal R_3$, the polygonal path between the points $3.01$, $3.01 c$, $3.49 c$ and $3.49$.

A similar proof to that of Theorem \ref{sdleverc} shows that $\Re[p(z)-p(z_2)]>0$ for $z \in \mathcal R$ except at $z=z_2$, as seen in Figure \ref{rpzbs0}. Hence
\begin{equation*}
\Re[-p(z)]\lqs \Re[-p(z_2)] \approx 0.013764 \qquad (z \in \mathcal R)
\end{equation*}
and it follows that the term corresponding to $d=0$ in \eqref{murrs} is $O(e^{0.014 N})$.


\SpecialCoor
\psset{griddots=5,subgriddiv=0,gridlabels=0pt}
\psset{xunit=0.7cm, yunit=0.8cm}
\psset{linewidth=1pt}
\psset{dotsize=4pt 0,dotstyle=*}

\begin{figure}[h]
\begin{center}
\begin{pspicture}(-1,-0.8)(14,2.5) 

\psline[linecolor=gray]{->}(0,-1)(0,2.5)
\psline[linecolor=gray]{->}(-1,0)(13.5,0)
\psline[linecolor=gray,linestyle=dashed](4,0)(4,1.04245)
\psline[linecolor=gray,linestyle=dashed](8,0)(8,0.978138)
\psline[linecolor=gray,linestyle=dashed](13,0)(13,0.305678)

  \rput(-1,1){$0.01$}

   \rput(-1,2){$0.02$}
  \rput(2,-0.6){$\mathcal R_1$}
  \rput(10.5,-0.6){$\mathcal R_3$}
  \rput(6,-0.6){$\mathcal R_2$}

\savedata{\mydata}[
{{0., 1.25161}, {0.1, 0.904793}, {0.2, 0.726667}, {0.3,
  0.621323}, {0.4, 0.552809}, {0.5, 0.506318}, {0.6, 0.474508}, {0.7,
  0.453217}, {0.8, 0.439846}, {0.9, 0.432657}, {1., 0.430424}, {1.1,
  0.432248}, {1.2, 0.437444}, {1.3, 0.445481}, {1.4, 0.455935}, {1.5,
  0.468461}, {1.6, 0.482776}, {1.7, 0.498645}, {1.8, 0.515868}, {1.9,
  0.534275}, {2., 0.553722}, {2.1, 0.574082}, {2.2, 0.595246}, {2.3,
  0.617118}, {2.4, 0.639614}, {2.5, 0.66266}, {2.6, 0.68619}, {2.7,
  0.710145}, {2.8, 0.734474}, {2.9, 0.759129}, {3., 0.784069}, {3.1,
  0.809257}, {3.2, 0.834658}, {3.3, 0.860242}, {3.4, 0.885982}, {3.5,
  0.911852}, {3.6, 0.937831}, {3.7, 0.963898}, {3.8, 0.990035}, {3.9,
  1.01622}, {4., 1.04245}, {4.1, 1.07689}, {4.2, 1.11022}, {4.3,
  1.14222}, {4.4, 1.1727}, {4.5, 1.20149}, {4.6, 1.22842}, {4.7,
  1.25335}, {4.8, 1.27617}, {4.9, 1.29677}, {5., 1.31509}, {5.1,
  1.33105}, {5.2, 1.34462}, {5.3, 1.35577}, {5.4, 1.36451}, {5.5,
  1.37083}, {5.6, 1.37477}, {5.7, 1.37635}, {5.8, 1.37565}, {5.9,
  1.37271}, {6., 1.3676}, {6.1, 1.36043}, {6.2, 1.35126}, {6.3,
  1.34021}, {6.4, 1.32738}, {6.5, 1.31288}, {6.6, 1.29683}, {6.7,
  1.27934}, {6.8, 1.26053}, {6.9, 1.24054}, {7., 1.21948}, {7.1,
  1.19749}, {7.2, 1.17469}, {7.3, 1.1512}, {7.4, 1.12716}, {7.5,
  1.10269}, {7.6, 1.07792}, {7.7, 1.05296}, {7.8, 1.02793}, {7.9,
  1.00295}, {8., 0.978138}, {8.1, 0.964204}, {8.2, 0.950304}, {8.3,
  0.936439}, {8.4, 0.922612}, {8.5, 0.908823}, {8.6, 0.895075}, {8.7,
  0.881366}, {8.8, 0.8677}, {8.9, 0.854077}, {9., 0.840497}, {9.1,
  0.826961}, {9.2, 0.81347}, {9.3, 0.800025}, {9.4, 0.786624}, {9.5,
  0.773269}, {9.6, 0.759959}, {9.7, 0.746693}, {9.8, 0.733471}, {9.9,
  0.720292}, {10., 0.707155}, {10.1, 0.694057}, {10.2,
  0.680997}, {10.3, 0.667973}, {10.4, 0.654982}, {10.5,
  0.64202}, {10.6, 0.629084}, {10.7, 0.61617}, {10.8,
  0.603274}, {10.9, 0.590389}, {11., 0.577511}, {11.1,
  0.564633}, {11.2, 0.551747}, {11.3, 0.538847}, {11.4,
  0.525923}, {11.5, 0.512967}, {11.6, 0.499968}, {11.7,
  0.486916}, {11.8, 0.473798}, {11.9, 0.460602}, {12.,
  0.447315}, {12.1, 0.433922}, {12.2, 0.420408}, {12.3,
  0.406755}, {12.4, 0.392946}, {12.5, 0.378962}, {12.6,
  0.364784}, {12.7, 0.350389}, {12.8, 0.335756}, {12.9,
  0.32086}, {13., 0.305678}}
]

\multirput(-0.15,1)(0,1){2}{\psline[linecolor=gray](0,0)(0.3,0)}

\dataplot[linecolor=red,linewidth=0.8pt,plotstyle=line]{\mydata}

  \rput(5.71878,0.5){$z_2$}

\psdots(5.71878,1.3764)(5.71878,0)

\rput(8.5,1.7){$\Re[-p(z)]$}

\end{pspicture}
\caption{Graph of $\Re[-p(z)]$ for $z \in \mathcal R$}\label{rpzbs0}
\end{center}
\end{figure}
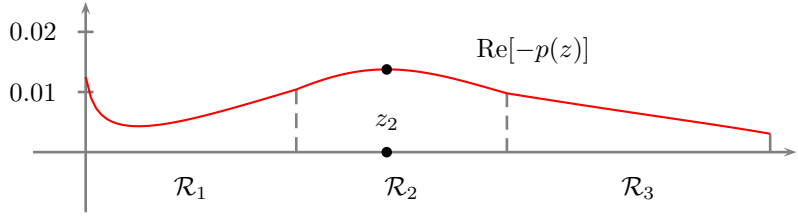


Define
\begin{equation} \label{c4s}
    \mathcal C_4^*(N,\sigma)  := \frac{-1}{4(N+1/2)^{3/2}}   \Im \int_{3.01}^{3.49}  \exp\bigl(-(N+1/2) p_1(z) \bigr) q^*_\mathcal C(z) \exp \bigl(v_\mathcal C^*(z;N,\sigma)\bigr) \, dz
\end{equation}
and we now know from \eqref{coss}, \eqref{murrs} and the above that
\begin{equation}\label{cosm4}
    \mathcal C_2^*(N,\sigma) = \mathcal C_4^*(N,\sigma) +O(e^{WN/3}).
\end{equation}
The unique solution to $p_1'(z)=0$ for $5/2<\Re(z)<7/2$ is
\begin{equation*}
    z_3:=3+\frac{\log \bigl(1-w(1,-3)\bigr)}{2\pi i} \approx 3.08382 - 0.0833451  i
\end{equation*}
by Theorem \ref{disol}.
Let $v=\Im(z_3)/\Re(z_3) \approx -0.027027 $ and $c=1+i v $. The path we take through the saddle point $z_3$ is $\mathcal S := \mathcal S_1 \cup \mathcal S_2 \cup \mathcal S_3$, the polygonal path between the points $3.01$, $3.01 c$, $3.49 c$ and $3.49$.
A similar proof to that of Theorem \ref{sdleverc} shows that $\Re[p_1(z)-p_1(z_3)]>0$ for $z \in \mathcal S$ except at $z=z_3$. This is seen in Figure \ref{rpzbs1}.


\SpecialCoor
\psset{griddots=5,subgriddiv=0,gridlabels=0pt}
\psset{xunit=0.7cm, yunit=0.4cm}
\psset{linewidth=1pt}
\psset{dotsize=4pt 0,dotstyle=*}

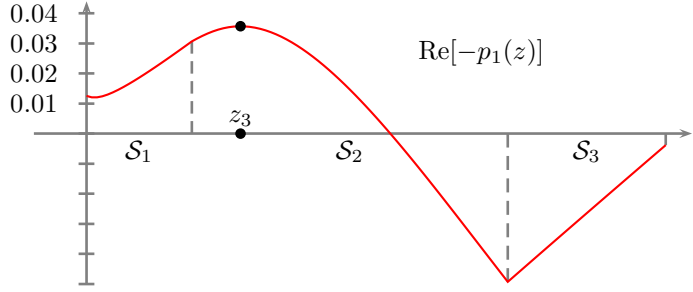
\begin{figure}[h]
\begin{center}
\begin{pspicture}(-1,-4.8)(12,4) 

\psline[linecolor=gray]{->}(0,-5)(0,4.4)
\psline[linecolor=gray]{->}(-1,0)(11.5,0)
\psline[linecolor=gray,linestyle=dashed](2,0)(2,3.06538)
\psline[linecolor=gray,linestyle=dashed](8,0)(8,-4.93269)
\psline[linecolor=gray,linestyle=dashed](11,0)(11,-0.384733)

  \rput(-1,1){$0.01$}
\rput(-1,2){$0.02$}
   \rput(-1,3){$0.03$}
   \rput(-1,4){$0.04$}
  \rput(1,-0.6){$\mathcal S_1$}
  \rput(9.5,-0.6){$\mathcal S_3$}
  \rput(5,-0.6){$\mathcal S_2$}

\savedata{\mydata}[
{{0., 1.25161}, {0.1, 1.20411}, {0.2, 1.20419}, {0.3, 1.23846}, {0.4,
  1.29584}, {0.5, 1.3691}, {0.6, 1.45369}, {0.7, 1.54669}, {0.8,
  1.64614}, {0.9, 1.75067}, {1., 1.8593}, {1.1, 1.97128}, {1.2,
  2.08605}, {1.3, 2.20317}, {1.4, 2.32229}, {1.5, 2.4431}, {1.6,
  2.56537}, {1.7, 2.6889}, {1.8, 2.8135}, {1.9, 2.93904}, {2.,
  3.06538}, {2.1, 3.16021}, {2.2, 3.24837}, {2.3, 3.32799}, {2.4,
  3.3976}, {2.5, 3.45608}, {2.6, 3.50266}, {2.7, 3.53689}, {2.8,
  3.55855}, {2.9, 3.56763}, {3., 3.56426}, {3.1, 3.54869}, {3.2,
  3.52123}, {3.3, 3.48227}, {3.4, 3.43222}, {3.5, 3.3715}, {3.6,
  3.30056}, {3.7, 3.21983}, {3.8, 3.12975}, {3.9, 3.03076}, {4.,
  2.92327}, {4.1, 2.80769}, {4.2, 2.68441}, {4.3, 2.55381}, {4.4,
  2.41627}, {4.5, 2.27212}, {4.6, 2.12172}, {4.7, 1.96539}, {4.8,
  1.80344}, {4.9, 1.63618}, {5., 1.46389}, {5.1, 1.28686}, {5.2,
  1.10536}, {5.3, 0.919649}, {5.4, 0.729973}, {5.5, 0.536576}, {5.6,
  0.339693}, {5.7,
  0.139548}, {5.8, -0.0636407}, {5.9, -0.269661}, {6., -0.478309},
{6.1, -0.689385}, {6.2, -0.902698}, {6.3, -1.11806}, {6.4, -1.33529},
{6.5, -1.55421}, {6.6, -1.77465}, {6.7, -1.99644}, {6.8, -2.21941},
{6.9, -2.44342}, {7., -2.66829}, {7.1, -2.89388}, {7.2, -3.12004},
{7.3, -3.34662}, {7.4, -3.57348}, {7.5, -3.80047}, {7.6, -4.02746},
{7.7, -4.25431}, {7.8, -4.48089}, {7.9, -4.70706}, {8., -4.93269},
{8.1, -4.77966}, {8.2, -4.62674}, {8.3, -4.47394}, {8.4, -4.32124},
{8.5, -4.16866}, {8.6, -4.01618}, {8.7, -3.86381}, {8.8, -3.71154},
{8.9, -3.55937}, {9., -3.4073}, {9.1, -3.25533}, {9.2, -3.10346},
{9.3, -2.95169}, {9.4, -2.80001}, {9.5, -2.64842}, {9.6, -2.49692},
{9.7, -2.34552}, {9.8, -2.1942}, {9.9, -2.04297}, {10., -1.89182},
{10.1, -1.74075}, {10.2, -1.58977}, {10.3, -1.43887}, {10.4,
-1.28805}, {10.5, -1.13731}, {10.6, -0.986647}, {10.7, -0.836058},
{10.8, -0.685543}, {10.9, -0.535102}, {11., -0.384733}}
]

\multirput(-0.15,-5)(0,1){10}{\psline[linecolor=gray](0,0)(0.3,0)}

\dataplot[linecolor=red,linewidth=0.8pt,plotstyle=line]{\mydata}

  \rput(2.9227,0.5){$z_3$}

\psdots(2.9227,3.56795)(2.9227,0)

\rput(7.5,2.7){$\Re[-p_1(z)]$}

\end{pspicture}
\caption{Graph of $\Re[-p_1(z)]$ for $z \in \mathcal S$}\label{rpzbs1}
\end{center}
\end{figure}


\subsection{Applying the saddle-point method}

Recall \eqref{gazi}
and for $j \in \Z_{\gqs 0}$ put
\begin{equation*}
    u^*_{\sigma, j}(z):=\sum_{m_1+3m_2+5m_3+ \dots =j}\frac{(\pi i (16\sigma+1)z/8 +g_{\mathcal C, 1}(z))^{m_1}}{m_1!}\frac{g_{\mathcal C,  2}(z)^{m_2}}{m_2!} \cdots \frac{g_{\mathcal C, j}(z)^{m_j}}{m_j!},
\end{equation*}
with $u^*_{\sigma,  0}=1$.
Similarly to Proposition  \ref{gasc} we have

\begin{prop} \label{gascs}
For $3.01 \lqs \Re(z) \lqs 3.49$ and $|\Im(z)|\lqs 1$, say, there is a holomorphic function $\zeta_d^*(z;N,\sigma)$ of $z$ so that
\begin{equation*}
    \exp\bigl(v_{\mathcal C}^*(z;N,\sigma)\bigr) = \sum_{j=0}^{d-1} \frac{u^*_{\sigma,j}(z)}{(N+1/2)^j} + \zeta_d^*(z;N,\sigma) \quad \text{for} \quad \zeta_d^*(z;N,\sigma) = O\left(\frac{1}{N^d} \right)
\end{equation*}
with an implied constant depending only on $\sigma$ and  $d$ where $1 \lqs d \lqs 2L-1$ and $L=\lfloor 0.006 \pi e \cdot N/2 \rfloor$.
\end{prop}

\begin{theorem} \label{c2ses} With $c_{0}^*=-z_3 e^{-\pi i z_3}/4$ and  explicit  $c_{1}^*(\sigma)$, $c_{2}^*(\sigma), \dots $ depending on $\sigma \in \R$ we have
\begin{equation} \label{press}
   \mathcal C_2^*(N,\sigma) = \Re\left[\frac{w(1,-3)^{-N}}{N^{2}} \left( c_{0}^*+\frac{c_{1}^*(\sigma)}{N}+ \dots +\frac{c_{m-1}^*(\sigma)}{N^{m-1}}\right)\right] + O\left(\frac{|w(1,-3)|^{-N}}{N^{m+2}}\right)
\end{equation}
for an implied constant depending only on  $\sigma$ and $m$.
\end{theorem}
\begin{proof}
As in Theorem \ref{c2se}, applying the saddle-point method to \eqref{c4s}, with the path of integration moved to $\mathcal S$, yields
\begin{equation*}
    \mathcal C_4^*(N,\sigma)  = \Re \left[  e^{-(N+1/2) \cdot p_1(z_3)}
    \sum_{t=0}^{d-2} \frac{1}{(N+1/2)^{t+2}}    \sum_{s=0}^{t} \frac{\G(s+1/2)}2 a_{2s}(i q_{\mathcal C}^* \cdot u_{\sigma, t-s}^*)
    \right]+ O\left( \frac{|w(1,-3)|^{-N}}{N^{d+1}}\right).
\end{equation*}
From \eqref{pzlogw} we know that
$e^{p_1(z_3)} = w(1,-3)$.
Hence, set
\begin{equation} \label{ctys}
    c_t^{**}(\sigma):=   e^{-p_1(z_3)/2}\sum_{s=0}^t \G(s+1/2) a_{2s}(i q_{\mathcal C}^* \cdot u_{\sigma, t-s}^*)/2.
\end{equation}
We want to convert the above series in $1/(N+1/2)$ to one in $1/N$. With the Binomial Theorem we have $(1+z)^{-j}=\sum_{r=0}^\infty \binom{-j}{r} z^r$ for $|z|<1$. Also, by Taylor's Theorem,
\begin{equation}\label{tayl}
    \frac{1}{(1+z)^{j}}=\sum_{r=0}^{m-1}\binom{-j}{r} z^r +O(z^m) \qquad(|z|<1).
\end{equation}
With $z=1/(2N)$  above we find
\begin{equation*}
    \frac{\alpha_j}{(N+1/2)^{j+2}}=\sum_{r=0}^\infty \binom{-j-2}{r} \frac{2^{-r} \cdot \alpha_j}{N^{j+2+r}}
\end{equation*}
for any $\alpha_j$s, and can write
\begin{equation*}
    \frac{\alpha_0}{(N+1/2)^2}+\frac{\alpha_1}{(N+1/2)^3}+\cdots = \frac{\beta_0}{N^2}+\frac{\beta_1}{N^3}+\cdots
\end{equation*}
with
\begin{equation*}
    \beta_t=\sum_{j+r=t} \binom{-j-2}{r} 2^{-r} \cdot \alpha_j = \sum_{j=0}^t \binom{-j-2}{t-j} 2^{j-t} \cdot \alpha_j = \sum_{j=0}^t (-2)^{j-t}\binom{t+1}{j+1} \alpha_j.
\end{equation*}
So we set
\begin{equation*} 
     c_t^{*}(\sigma):= \sum_{j=0}^t (-2)^{j-t}\binom{t+1}{j+1}  c_j^{**}(\sigma)
\end{equation*}
and with \eqref{cosm4} we obtain \eqref{press} in the statement of the theorem.
Note that the omitted terms satisfy
\begin{equation*}
    \sum_{t=m}^\infty \frac{c_{t}^*(\sigma)}{N^t} = O\left(\frac 1{N^{m}}\right)
\end{equation*}
by \eqref{tayl} and can be incorporated into the error term of \eqref{press}.

A similar computation to that of $c_0$ in the proof of Theorem \ref{c2se} shows that
\begin{equation*}
    \bigl(c_0^*(\sigma)\bigr)^2 =  z_3^2 e^{-2\pi i z_3}/16
\end{equation*}
and a numerical check then indicates that the correct square root has a minus sign.
\end{proof}

For example, Table \ref{c2sn1} compares both sides of \eqref{press} in Theorem \ref{c2ses} for some different values of $m$ and $N$. This is for $\sigma=1$ and the results for other values of $\sigma$ are similar.
\begin{table}[h]
\begin{center}
\begin{tabular}{c|cccc|c}
$N$ & $m=1$ & $m=2$ & $m=3$ & $m=4$ & $\mathcal C^*_2(N,1)$  \\ \hline
$800$ & $1.43938 \times 10^{6}$ &  $1.39381 \times 10^{6}$ & $1.3934 \times 10^{6}$ &  $1.39341 \times 10^{6}$ &  $1.39341 \times 10^{6}$ \\
$1000$ & $1.7278 \times 10^{9}$ &  $1.74062 \times 10^{9}$ & $1.74028 \times 10^{9}$ &  $1.74028 \times 10^{9}$ &  $1.74028 \times 10^{9}$
\end{tabular}
\caption{Theorem \ref{c2ses}'s approximations to $\mathcal C^*_2(N,1)$.} \label{c2sn1}
\end{center}
\end{table}

A consequence of Theorem \ref{c2se} is that
\begin{equation}\label{c2bnd}
    \mathcal C_2(N,\sigma) = O(e^{U_\mathcal C N}/N^2) \qquad \text{for} \qquad U_\mathcal C:=-\log|w(0,-2)| \approx 0.0256706.
\end{equation}
Since $-\log|w(1,-3)| \approx 0.0356795$ we see  that $\mathcal C_2(N,\sigma)$ is much smaller than $\mathcal C_2^*(N,\sigma)$ (despite appearances in Figure \ref{b2}) and is bounded by the error term in \eqref{press}. Therefore,  Theorem \ref{thmc} on the asymptotic expansion of $\mathcal C_1(N,\sigma) = \mathcal C_2(N,\sigma)+\mathcal C_2^*(N,\sigma)$  follows from Theorems \ref{c2se} and \ref{c2ses}.

\section{The sum $\mathcal D_1(N,\sigma)$}
Let $\sigma \in \Z$. In this section  we prove Theorem \ref{c2sed}, giving the asymptotic expansion as $N \to \infty$ of
$$
\mathcal D_1(N,\sigma)  := \sum_{h/k \in \mathcal D(N)} Q_{hk\sigma}(N) = 2 \Re \sum_{\frac{N}{2}  <k \lqs N, \ k \text{ odd}}  Q_{(\frac{k-1}{2})k\sigma}(N).
$$
With $k$  odd, setting $h=(k-1)/2$ in Proposition \ref{simple}  yields
\begin{multline}\label{rgg}
Q_{\left(\frac{k-1}{2}\right) k\sigma}(N)=\frac{1}{k^2}
\exp\left( \frac{\pi i}{4}\left[ \frac{N^2+N-4\sigma}{k} \right]\right)\\
\times
\exp\left( \frac{-\pi i}{4}\Bigl[(N-k)(N-k+1)-3k-3 \Bigr]\right) \spr{(k-1)/2k}{N-k}.
\end{multline}

\subsection{$\mathcal D_1(N,\sigma)$ for $N$ odd}
If $N$ is odd then $N-k$ is even and $(N-k)(N-k+1) \equiv k-N \bmod 8$. Hence \eqref{rgg} becomes
\begin{multline}\label{nodd}
    Q_{\left(\frac{k-1}{2}\right) k\sigma}(N)= \frac{1}{k^2}
\exp\left( N\left[ \frac{\pi i}{4}\left(\frac{N}{k}+1+\frac{2k}N \right)\right]\right)
\\
\times
\exp\left( \frac{\pi i}{4}\left(\frac{N}{k}+3\right)\right)
\exp\left( \frac 1N\left[ -\pi i \sigma \frac{N}{k}\right]\right)
\spr{(k-1)/2k}{N-k}.
\end{multline}

We next get $\spr{(k-1)/2k}{N-k}$ into the right form to apply Proposition \ref{prrh}.
\begin{prop} \label{bigprod}
For $k$  odd and $m$  even with $0\lqs m <k$ we have
\begin{equation}\label{bpsp}
    \spr{(k-1)/2k}{m} = \frac{\spr{1/k}{m} }{\spr{1/2k}{m }} \times \frac{ \spp{1/k}{m/2}}{ \spr{2/k}{m/2} }.
\end{equation}
\end{prop}
\begin{proof}
Since
$$
\sin (\pi j (k-1)/2k) = \begin{cases}
(-1)^{j/2+1} \sin(\pi j/2k) & \quad j \text{ even} \\
(-1)^{(j-1)/2} \cos(\pi j/2k) & \quad j \text{ odd}
\end{cases}
$$
we have
\begin{equation}\label{sub}
\spr{(k-1)/2k}{m} = \prod_{\substack{1 \lqs j \lqs m \\ j \text{ even}}} \frac{(-1)^{j/2+1}}{2 \sin(\pi j/2k)}
\prod_{\substack{1 \lqs j \lqs m \\ j \text{ odd}}} \frac{(-1)^{(j-1)/2}}{2 \cos(\pi j/2k)}.
\end{equation}
Hence
\begin{align}
    \spr{(k-1)/2k}{m} & = \spr{1/k}{m/2}
\prod_{\substack{1 \lqs j \lqs m \\ j \text{ odd}}} \frac{1}{2 \cos(\pi j/2k)} \notag\\
& = \left. \spr{1/k}{m/2}
\prod_{1 \lqs j \lqs m } \frac{1}{2 \cos(\pi j/2k)}\right/ \prod_{1 \lqs j \lqs m/2 } \frac{1}{2 \cos(\pi j/k)}. \label{cos}
\end{align}
Use the identity $2\sin 2\theta = 2\sin \theta \cdot 2\cos \theta$ to convert the cosines in \eqref{cos} back to sines and complete the proof.
\end{proof}

Recall the definition of $g_\ell(z)$ in \eqref{grz}, define
\begin{equation}\label{grzs}
g_\ell^*(z):=-\frac{B_{2\ell}}{(2\ell)!} \left( \pi z/2\right)^{2\ell-1} \cot^{(2\ell-2)}\left(\pi (z-1)/2\right)
\end{equation}
and set $z=z(N,k):=N/k$. The sine product approximation we need is as follows.

\begin{theorem} \label{cplxc2}
Fix $W>0$. Let $\Delta$ be in the range $0.0048 \lqs \Delta \lqs 0.0079$ and set $\alpha = \Delta \pi e$.  Suppose $\delta$ and $\delta'$ satisfy
\begin{equation*}
    \frac{\Delta}{1-\Delta} < \delta \lqs \frac{1}{e}, \ 0<\delta'\lqs \frac{1}{e} \quad \text{ and  } \quad \delta \log 1/\delta, \ \ \delta' \log 1/\delta' \lqs W.
\end{equation*}
Then for all $N$ odd $\gqs 2 \cdot R_\Delta$  we have
\begin{equation}\label{ott1c2}
    \spr{(k-1)/2k}{N-k} = O\left(e^{WN/2}\right) \quad \text{ for } \quad  z \in \left[1, \ 1+ \delta \right] \cup  \left[3/2 - \delta' , \ 2 \right)
\end{equation}
and
\begin{multline} \label{ott2c2}
    \spr{(k-1)/2k}{N-k} =  \exp\left(N\frac {\cl(2\pi  z)}{4\pi z}  \right) \left( \frac{z}{2N\sin(\pi (z-1)/2)}\right)^{1/2}
    \exp\left(\sum_{\ell=1}^{L-1} \frac{g_{\ell}(z)-g^*_{\ell}(z)}{N^{2\ell-1}}\right)
      \\
      \quad \times \exp\left(\sum_{\ell=1}^{L^*-1} \frac{2 g^*_{\ell}(z) - g_{\ell}(z)}{(N/2)^{2\ell-1}}\right) +O\left(e^{WN/2}\right)
        \quad \text{ for } \quad z \in  \left(1+\delta, \ 3/2 - \delta'\right)
  \end{multline}
with $L=\lfloor \alpha \cdot N\rfloor$ and $L^*=\lfloor \alpha \cdot N/2\rfloor$. The implied constants in \eqref{ott1c2}, \eqref{ott2c2} are absolute.
\end{theorem}
\begin{proof}
Applying Proposition \ref{prrh} to each of the factors on the right of  \eqref{bpsp} shows
\begin{multline} \label{eva3}
\spr{(k-1)/2k}{m}  =  \left( \frac{1}{2k \sin(\pi m /(2k))}\right)^{1/2} \exp\left(\frac{k}{4\pi}\cl(2\pi m/k) \right)  \\
  \times \exp\left(-\sum_{\ell=1}^{L_1-1} \frac{B_{2\ell}}{(2\ell)!} \left( \frac{\pi }{k}\right)^{2\ell-1}
\cot^{(2\ell-2)}\left(\frac{\pi m}{k}\right)\right) \exp\left(-T_{L_1}(m,1/k)\right)
\\
 \times \exp\left(-2\sum_{\ell=1}^{L_2-1} \frac{B_{2\ell}}{(2\ell)!} \left( \frac{\pi }{k}\right)^{2\ell-1}
\cot^{(2\ell-2)}\left(\frac{\pi m}{2k}\right)\right) \exp\left(-2 T_{L_2}(m/2,1/k)\right)
\\
 \times \exp\left(\sum_{\ell=1}^{L_3-1} \frac{B_{2\ell}}{(2\ell)!} \left( \frac{2\pi }{k}\right)^{2\ell-1}
\cot^{(2\ell-2)}\left(\frac{\pi m}{k}\right)\right) \exp\left(T_{L_3}(m/2,2/k)\right)
\\
  \times \exp\left(\sum_{\ell=1}^{L_4-1} \frac{B_{2\ell}}{(2\ell)!} \left( \frac{\pi }{2k}\right)^{2\ell-1}
\cot^{(2\ell-2)}\left(\frac{\pi m}{2k}\right)\right) \exp\left(T_{L_4}(m,1/(2k))\right)
\end{multline}
for  $1\lqs m <k$ and positive integers $L_1$, $L_2$, $L_3$, $L_4$.

First we set each $L_i$ to $1$ in \eqref{eva3}  to see
\begin{multline} \label{eva3x}
\spr{(k-1)/2k}{m} =  \left( \frac{1}{2k \sin(\pi m /(2k))}\right)^{1/2} \exp\left(\frac{k}{4\pi}\cl(2\pi m/k) \right)
\\
\times \exp\left(-T_{1}(m,1/k)-2 T_{1}(m/2,1/k)+T_{1}(m/2,2/k)+T_{1}(m,1/(2k))\right).
\end{multline}
Comparing \eqref{eva3x} with the expansion of $\spr{2/k}{m/2}$ from Proposition \ref{prrh} then shows
\begin{multline} \label{eva3xx}
\spr{(k-1)/2k}{m} =  \left( \cos(\pi m /(2k))\right)^{1/2}
\\
\times \exp\left(-T_{1}(m,1/k)-2 T_{1}(m/2,1/k)+2T_{1}(m/2,2/k)+T_{1}(m,1/(2k))\right) \spr{2/k}{m/2}.
\end{multline}
It follows from \eqref{eva3xx} and Lemma \ref{melrh2} that for $0\lqs m <k$,
\begin{equation}\label{uu1}
    \spr{(k-1)/2k}{m}  \ll  \spr{2/k}{m/2}
\end{equation}
with an absolute implied constant. Similarly, by comparing \eqref{eva3x} with the expansion of $\spr{1/k}{m}$ from Proposition \ref{prrh},
\begin{equation}
    \spr{(k-1)/2k}{m}  \ll \left( \spr{1/k}{m}\right)^{1/2}. \label{uu2}
\end{equation}
Using Proposition \ref{gprop} to bound $\spr{2/k}{m/2}$ on the right of \eqref{uu1} and noting that $k \lqs N$ proves \eqref{ott1c2}.

To prove \eqref{ott2c2} we wish to apply the argument of Proposition \ref{hard} to \eqref{eva3}. This requires finding $L_1$, $L_2$, $L_3$ and $L_4$ so that, for $m=N-k$,
\begin{align}
    \spr{(k-1)/2k}{m} \Bigl(\left|  T_{L_1}(m,1/k) \right| +\left|  T_{L_2}(m/2,1/k) \right| +\left|  T_{L_3}(m/2,2/k) \right| +\left|  T_{L_4}(m,1/(2k)) \right|   \Bigr) & \ll e^{WN/2} \label{eon1}\\
    \left|  T_{L_1}(m,1/k) \right| +\left|  T_{L_2}(m/2,1/k) \right| +\left|  T_{L_3}(m/2,2/k) \right| +\left|  T_{L_4}(m,1/(2k)) \right|  & \ll 1  . \label{eon2}
\end{align}
We examine the four terms $T_{L_i}$ in \eqref{eon1} and \eqref{eon2} separately:
\begin{itemize}
\item{\bf The term $T_{L_3}(m/2,2/k)$.} Use Proposition \ref{dela} with $h=2$ and $s=N$ to show that, for $\Delta N/2 \lqs m/2 \lqs k/4$,
\begin{align}
    \left| \spr{2/k}{m/2} \cdot T_{L_3}(m/2,2/k) \right| & \ll e^{WN/2} \label{iin3}\\
    \left| T_{L_3}(m/2,2/k) \right| & \ll 1  \label{iin4}
\end{align}
with absolute implied constants, $L_3:=\lfloor \pi e \Delta \cdot N/2 \rfloor$ and $N \gqs 2 \cdot R_\Delta$. Inequality \eqref{iin3} is valid with $\spr{2/k}{m/2}$ replaced by $\spr{(k-1)/2k}{m}$ using \eqref{uu1}:
\begin{equation}
    \left| \spr{(k-1)/2k}{m} \cdot  T_{L_3}(m/2,2/k) \right| \ll e^{WN/2} \label{iin6}.
\end{equation}

\item{\bf The term $T_{L_2}(m/2,1/k)$.} To prove
\begin{align}
    \left| \spr{(k-1)/2k}{m} \cdot  T_{L_2}(m/2,1/k) \right| & \ll e^{WN/2} \label{ijn3}\\
    \left| T_{L_2}(m/2,1/k) \right| & \ll 1  \label{ijn4}
\end{align}
for $\Delta N/2 \lqs m/2 \lqs k/4$, choose $L_2=L_3$ and note that  \eqref{iin3} and \eqref{iin4} are valid with $2/k$ replaced by $1/k$ using Corollary \ref{delacor}.

\item{\bf The term $T_{L_1}(m,1/k)$.} Use Proposition \ref{dela} with $h=1$ and $s=N$ to show that, also for $\Delta N \lqs m \lqs k/2$,
\begin{align}
    \left| \spr{1/k}{m} \cdot  T_{L_1}(m,1/k) \right| & \ll e^{WN} \label{iin9}\\
    \left| T_{L_1}(m,1/k) \right| & \ll 1  \label{iin10}
\end{align}
with absolute implied constants, $L_1:=\lfloor \pi e \Delta \cdot N \rfloor$  and $N \gqs  R_\Delta$. Taking square roots of both sides of \eqref{iin9} and using \eqref{uu2} shows
\begin{equation}
    \left| \spr{(k-1)/2k}{m}  \cdot T_{L_1}(m,1/k) \right|  \ll e^{WN/2} \label{iin12}.
\end{equation}

\item{\bf The term $T_{L_4}(m,1/(2k))$.} To prove
\begin{align}
    \left| \spr{(k-1)/2k}{m}  \cdot T_{L_4}(m,1/(2k)) \right| & \ll e^{WN/2} \label{ijn7}\\
    \left| T_{L_4}(m,1/(2k)) \right| & \ll 1  \label{ijn8}
\end{align}
for $\Delta N \lqs m \lqs k/2$, choose $L_4=L_1$ and note that  \eqref{iin9} and \eqref{iin10} are valid with $1/k$ replaced by $1/(2k)$ using Corollary \ref{delacor}.

\end{itemize}

 The inequalities \eqref{iin4} - \eqref{ijn8} establish \eqref{eon1}, \eqref{eon2} and the arguments of Proposition \ref{hard} now go through, applied to \eqref{eva3}. This allows us to remove the $\exp(T_{L_i})$ factors in \eqref{eva3} at the expense of adding an $O(e^{WN/2})$ error. Write $L$ for $L_1$, $L_4$ and $L^*$ for $L_2$, $L_3$. The interval $\Delta N \lqs m \lqs k/2$ corresponds to
\begin{equation*}
     1 + \frac{\Delta}{1-\Delta} \lqs z \lqs 3/2.
\end{equation*}
This completes the proof of \eqref{ott2c2}.
\end{proof}

It  simplifies things to work with the conjugate of \eqref{nodd}:
\begin{equation}\label{nodd2}
    \overline{Q_{\left(\frac{k-1}{2}\right) k\sigma}(N)}= \frac{1}{k^2}
\exp\left( N\left[ \frac{-\pi i}{4}\left(z+1+\frac{2}z \right)\right]\right)
\exp\left( \frac{-\pi i}{4}\left(z+3\right)\right)
\exp\left( \frac {\pi i \sigma z}N\right)
\spr{(k-1)/2k}{N-k}.
\end{equation}
From \eqref{dli} we have
\begin{equation*}
    \cl(2\pi z)= -i\li(e^{2\pi i z}) +i\pi^2(z^2-3z+13/6) \qquad(1<z<2)
\end{equation*}
so that
\begin{equation*}
    \frac {\cl(2\pi  z)}{4\pi z} - \frac{\pi i}{4}\left( z+1+\frac 2{z}\right) = -\pi i +\frac 1{4\pi i z}\Bigl[\li(e^{2\pi i z})-\li(1) \Bigr].
\end{equation*}
Set
\begin{align}
    r_\mathcal D(z)  := & \frac 1{4\pi i z}\Bigl[\li(e^{2\pi i z})-\li(1) \Bigr] \label{rcz}\\
    q_\mathcal D(z)  := & \left( \frac{z }{2\sin(\pi(z -1)/2)}\right)^{1/2} \exp\left( -\frac{\pi i}{4}\left(z+3\right)\right)  \notag\\
    v_\mathcal D(z;N, \sigma) := & \frac{\pi i \sigma z}{N} + \sum_{\ell=1}^{L-1} \frac{g_{\ell}(z)-g^*_{\ell}(z)}{N^{2\ell-1}}
    +\sum_{\ell=1}^{L^*-1} \frac{2 g^*_{\ell}(z) - g_{\ell}(z)}{(N/2)^{2\ell-1}} \label{vcnz}
\end{align}
for $L=\lfloor \alpha \cdot N\rfloor$ and $L^*=\lfloor \alpha \cdot N/2\rfloor$.
With $N$ odd, define
\begin{equation} \label{d2sum}
    \mathcal D_2(N,\sigma)  :=  \frac {-2}{N^{1/2}} \Re \sum_{k  \text{ odd}\ : \ z \in  \left(1+\delta, \ 3/2 - \delta'\right)}
     \frac{1}{k^2} \exp \bigl(N \cdot r_\mathcal D(z) \bigr) q_\mathcal D(z) \exp \bigl(v_\mathcal D(z;N,\sigma)\bigr) .
\end{equation}
It follows from \eqref{nodd2} and Theorem \ref{cplxc2} that for $\sigma \in \R$ and an absolute implied constant
\begin{equation} \label{cosdo}
\mathcal D_1(N,\sigma) = \mathcal D_2(N,\sigma) +O(e^{WN/2}) \qquad (N \text{ odd}).
\end{equation}

\subsection{Expressing $\mathcal D_2(N,\sigma)$ as an integral for $N$ odd}
Similarly to Proposition \ref{addx} we have
\begin{prop} \label{addz}
Suppose $1/2 \lqs \Re(z) \lqs 3/2$ and $|z-1| \gqs \varepsilon >0$ and assume
$
    \max\{1+\frac{1}{\varepsilon}, \ 16 \} \ < \frac{\pi e}{\alpha}
$.
Then
\begin{equation} \label{wez}
    \sum_{\ell=d}^{L-1} \frac{g_{\ell}(z)-g^*_{\ell}(z)}{N^{2\ell-1}}
    +\sum_{\ell=d}^{L^*-1} \frac{2 g^*_{\ell}(z) - g_{\ell}(z)}{(N/2)^{2\ell-1}} \ll \frac{1}{N^{2d-1}}e^{-\pi|y|/2}
\end{equation}
for $d \gqs 2$ where $L=\lfloor \alpha \cdot N\rfloor$, $L^*=\lfloor \alpha \cdot N/2\rfloor$ and the implied constant depends only on $\varepsilon$, $\alpha$ and $d$.
\end{prop}

Fixing the choice of constants in \eqref{fix} and with $\varepsilon = 0.0061$
and
\begin{equation} \label{gazid}
    g_{\mathcal D, \ell}(z):=g_\ell(z) - g_\ell^*(z)+2^{2\ell-1}\bigl(2g_\ell^*(z) - g_\ell(z)\bigr)
\end{equation}
 we obtain:
\begin{cor} \label{acdcz}
With $\delta, \delta' \in [0.0061, 0.01]$ and  $z \in \C$ such that $1+\delta < \Re(z) < 3/2-\delta'$  we have
\begin{equation*}
    v_\mathcal D(z;N, \sigma) =  \frac{\pi i \sigma z}{N} + \sum_{\ell=1}^{d-1} \frac{g_{\mathcal D, \ell}(z)}{N^{2\ell-1}}
     + O\left( \frac 1{N^{2d-1}}\right)
\end{equation*}
 for $2 \lqs d\lqs L^*=\lfloor 0.006 \pi e \cdot N/2 \rfloor$ and an implied constant depending only on  $d$.
\end{cor}


Similarly to Theorem \ref{rzjb} we have

\begin{theorem} \label{rest}
The functions $r_\mathcal D(z)$, $q_\mathcal D(z)$ and $v_\mathcal D(z;N, \sigma)$ are holomorphic for $1< \Re(z)<3/2$. In this strip
\begin{alignat}{2}
    \Re\left(r_\mathcal D \left(z \right) +\frac{\pi i j}{z}\right) & \lqs \frac{1}{4\pi |z|^2} \left(x \cl(2\pi x) +\pi^2 |y|\left[\frac 13 +4j \right]  \right) \qquad & & (y\gqs 0) \label{efc1}\\
   \Re\left(r_\mathcal D \left(z \right) +\frac{\pi i j}{z}\right) & \lqs \frac{1}{4\pi |z|^2} \left(x \cl(2\pi x) +\pi^2 |y|\left[\frac 13 -4(j-1/2) \right] \right) \qquad & & (y\lqs 0). \label{efc2}
\end{alignat}
for $j\in \R$.
Also, in the box with $1+\delta \lqs \Re(z) \lqs 3/2-\delta'$ and $-1\lqs \Im(z) \lqs 1$,
\begin{equation} \label{efc3}
    q_\mathcal D(z), \quad \exp \bigl(v_\mathcal D (z;N,\sigma)\bigr) \ll 1
\end{equation}
for an  implied constant depending only on $\sigma \in \R$.
\end{theorem}

Let $C$ be the positively oriented rectangle with horizontal sides $C^+$, $C^-$ having imaginary parts $1/N^2$, $-1/N^2$ and vertical sides $C_{L}$, $C_{R}$ having real parts $1+\delta$ and $3/2-\delta'$ respectively, as used in \cite[Sect. 4.4]{OS1}.
Recalling \eqref{tantan}, \eqref{tantan2} and arguing as in Proposition \ref{cpcmb},  we find
\begin{align*}
  \mathcal D_2(N,\sigma) & =  \frac {(-2)}{N^{1/2}}  \frac{(-1)}{2N}\Re \int_C
     \exp\bigl(N \cdot r_\mathcal D(z) \bigr)  \frac{q_\mathcal D(z)}{2i\tan(\pi(N/z-1)/2)} \exp \bigl(v_\mathcal D(z;N,\sigma)\bigr) \, dz \\
  & =  \frac{1}{ N^{3/2}} \Re \left[ \sideset{}{'}{\sum}_{j \lqs 0} (-1)^j \int_{C^+}
    \exp \bigl(N [ r_\mathcal D(z)+\pi i j/z]\bigr) q_\mathcal D(z) \exp \bigl(v_\mathcal D(z;N,\sigma)\bigr) \, dz\right. \\
  &  \qquad \left. - \sideset{}{'}{\sum}_{j \gqs 0} (-1)^j \int_{C^-}
    \exp \bigl(N [ r_\mathcal D(z)+\pi i j/z]\bigr) q_\mathcal D(z) \exp \bigl(v_\mathcal D(z;N,\sigma)\bigr) \, dz\right] + O(e^{WN/2}).
\end{align*}
With Theorem \ref{rest}, and reasoning as in Proposition \ref{lgst}, we see that the two $j=0$ terms above dominate and $\mathcal D_2(N,\sigma) = \mathcal D_3(N,\sigma) +O(e^{WN/2})$ for
\begin{equation} \label{save}
    \mathcal D_3(N,\sigma) := \frac{-1}{ N^{3/2}} \Re \int_{1.01}^{1.49}
    \exp \bigl(-N \cdot p(z)/2 \bigr)  q_\mathcal D(z) \exp \bigl(v_\mathcal D(z;N,\sigma)\bigr) \, dz  \qquad (N \text{ odd})
\end{equation}
since $r_\mathcal D(z) = -p(z)/2$.

\subsection{$\mathcal D_1(N-1,\sigma)$ for $N$ odd}
Assume $N$ is odd. If $v$ is even then $v-k$ is odd and $(v-k)(v-k+1) \equiv v-k+1 \bmod 8$. Hence, with $v=N-1$, the conjugate of \eqref{rgg} becomes
\begin{equation}\label{neven}
    \overline{Q_{\left(\frac{k-1}{2}\right) k\sigma}(N-1)}= \frac{1}{k^2}
\exp\left( N\left[ \frac{-\pi i}{4}\left(z-1+\frac{4}{z} \right)\right]\right)
\exp\left( \frac{\pi i}{4}\left(z-3\right)\right)
\exp\left( \frac {\pi i \sigma z}N\right)
\spr{(k-1)/2k}{N-1-k}.
\end{equation}
For $m$ even, \eqref{sub} implies
\begin{align*}
    \spr{(k-1)/(2k)}{m-1} & = 2(-1)^{m/2+1} \sin(\pi m/(2k)) \prod_{\substack{1 \lqs j \lqs m \\ j \text{ even}}} \frac{(-1)^{j/2+1}}{2 \sin(\pi j/2k)}
\prod_{\substack{1 \lqs j \lqs m \\ j \text{ odd}}} \frac{(-1)^{(j-1)/2}}{2 \cos(\pi j/2k)}\\
& = 2(-1)^{m/2+1}  \sin(\pi m/(2k)) \cdot \spr{(k-1)/(2k)}{m} .
\end{align*}
It follows that for $N$ odd we have
\begin{equation} \label{tin}
    \spr{(k-1)/(2k)}{N-1-k}=2(-1)^{(N-k)/2+1}\sin(\pi (N/k-1)/2) \cdot  \spr{(k-1)/(2k)}{N-k}
\end{equation}
and can use our results from the last subsection.
Recall $r_\mathcal D(z)$ and $v_\mathcal D(z;N,\sigma)$ from \eqref{rcz}, \eqref{vcnz} and set
\begin{equation*}
    q^*_\mathcal D(z)  :=  2\sin(\pi(z -1)/2) \left( \frac{z }{2\sin(\pi(z -1)/2)}\right)^{1/2} \exp\left( \frac{\pi i}{4}\left(z-1\right)\right).
\end{equation*}
With $N$ odd, define
\begin{equation*}
    \mathcal D_2(N-1,\sigma)  :=  \frac {-2}{N^{1/2}} \Re \sum_{k \text{ odd} \ : \ z \in  \left(1+\delta, \ 3/2 - \delta'\right) }
     \frac{(-1)^{(k+1)/2}}{k^2} \exp \left(N  \left[r_\mathcal D(z) -\frac{\pi i}{2z}\right] \right) q^*_\mathcal D(z) \exp \bigl(v_\mathcal D(z;N,\sigma)\bigr) .
\end{equation*}
It follows from \eqref{neven}, \eqref{tin} and Theorem \ref{cplxc2} that
$$
\mathcal D_1(N-1,\sigma) = \mathcal D_2(N-1,\sigma) +O(e^{WN/2})  \qquad (N \text{ odd}).
$$
The next result is mostly a restatement of Theorem \ref{rest}.

\begin{theorem}
The functions $r_\mathcal D(z)-\frac{\pi i}{2z}$, $q^*_\mathcal D(z)$ and $v_\mathcal D(z;N,\sigma)$ are holomorphic for $1< \Re(z)<3/2$. In this strip
\begin{alignat}{2}
    \Re\left(r_\mathcal D \left(z \right) -\frac{\pi i}{2z} +\frac{\pi i (j-1/2)}{z}\right) & \lqs \frac{1}{4\pi |z|^2} \left(x \cl(2\pi x) +\pi^2 |y|\left[\frac 13 +4(j-1) \right]  \right) \qquad & & (y\gqs 0) \label{efc11}\\
   \Re\left(r_\mathcal D \left(z \right) -\frac{\pi i}{2z} +\frac{\pi i (j+1/2)}{z}\right) & \lqs \frac{1}{4\pi |z|^2} \left(x \cl(2\pi x) +\pi^2 |y|\left[\frac 13 -4(j-1/2) \right] \right) \qquad & & (y\lqs 0) \label{efc21}
\end{alignat}
for $j\in \R$.
Also, in the box with $1+\delta \lqs \Re(z) \lqs 3/2-\delta'$ and $-1\lqs \Im(z) \lqs 1$,
\begin{equation} \label{efc31}
    q^*_\mathcal D(z), \quad \exp \bigl(v_\mathcal D (z;N,\sigma)\bigr) \ll 1
\end{equation}
for an  implied constant depending only on $\sigma \in \R$.
\end{theorem}

With the rectangle $C$ from the last subsection and recalling \eqref{coscos}, \eqref{coscos2}
\begin{multline*}
  \mathcal D_2(N-1,\sigma)  =  \frac {(-2)}{N^{1/2}} \frac{(-1)}{2N} \Re
    \int_C \exp\left(N \left[r_\mathcal D(z) -\frac{\pi i}{2z}\right] \right) \frac{q^*_\mathcal D(z)}{2i\cos(\pi N/(2z))} \exp \bigl(v_\mathcal D (z;N,\sigma)\bigr) \, dz \\
   =   \frac{-1}{N^{3/2}} \Re \left[ i\sum_{j \lqs 0} (-1)^j \int_{C^+}
    \exp\left(N \left[r_\mathcal D \left(z \right) -\frac{\pi i}{2z} +\frac{\pi i (j-1/2)}{z}\right]\right) q^*_\mathcal D(z) \exp \bigl(v_\mathcal D (z;N,\sigma)\bigr) \, dz\right. \\
    \qquad \left. +i\sum_{j \gqs 0} (-1)^j \int_{C^-}
    \exp\left(N \left[r_\mathcal D \left(z \right) -\frac{\pi i}{2z} +\frac{\pi i (j+1/2)}{z}\right]\right) q^*_\mathcal D(z) \exp \bigl(v_\mathcal D (z;N,\sigma)\bigr) \, dz\right] +O(e^{WN/2}).
\end{multline*}
With \eqref{efc11}, \eqref{efc21} we see the $j=0$ term on $C^-$  dominates so that $\mathcal D_2(N-1,\sigma) = \mathcal D_3(N-1,\sigma) +O(e^{WN/2})$ for
\begin{equation} \label{save1}
    \mathcal D_3(N-1,\sigma) := \frac{-1}{ N^{3/2}} \Re \int_{1.01}^{1.49}
     \exp \bigl(-N \cdot p(z)/2 \bigr) i q^*_\mathcal D(z) \exp \bigl(v_\mathcal D (z;N,\sigma)\bigr) \, dz  \qquad (N \text{ odd}).
\end{equation}
Thus, with the definitions \eqref{save} and \eqref{save1} we have shown that for all $N$
\begin{equation}\label{cosdk}
    \mathcal D_1(N,\sigma) = \mathcal D_3(N,\sigma) +O(e^{WN/2}).
\end{equation}

\subsection{The asymptotic behavior of $\mathcal D_1(N,\sigma)$}

Recall \eqref{gazid}
and for $j \in \Z_{\gqs 0}$ put
\begin{equation*}
    u_{\mathcal D, \sigma, j}(z):=\sum_{m_1+3m_2+5m_3+ \dots =j}\frac{(\pi i \sigma z +g_{\mathcal D, 1}(z))^{m_1}}{m_1!}\frac{g_{\mathcal D, 2}(z)^{m_2}}{m_2!} \cdots \frac{g_{\mathcal D, j}(z)^{m_j}}{m_j!},
\end{equation*}
with $u_{\mathcal D, \sigma, 0}=1$. The proof of the next proposition is similar to Proposition \ref{gasc}'s and uses Corollary \ref{acdcz}.

\begin{prop} \label{gasd}
For $1.01 \lqs \Re(z) \lqs 1.49$ and $|\Im(z)|\lqs 1$, say, there is a holomorphic function $\zeta_{\mathcal D, d}(z;N,\sigma)$ of $z$ so that
\begin{equation*}
    \exp\bigl(v_{\mathcal D}(z;N,\sigma)\bigr) = \sum_{j=0}^{d-1} \frac{u_{\mathcal D, \sigma,j}(z)}{N^j} + \zeta_{\mathcal D, d}(z;N,\sigma) \quad \text{for} \quad \zeta_{\mathcal D, d}(z;N,\sigma) = O\left(\frac{1}{N^d} \right)
\end{equation*}
with an implied constant depending only on $\sigma$ and  $d$ where $1 \lqs d \lqs 2L^*-1$ and $L^*=\lfloor 0.006 \pi e \cdot N/2 \rfloor$.
\end{prop}

We restate Theorem \ref{c2sed} here. Recall that $z_0=1+\log(1-w_0)/(2\pi i)$ where $w_0$  is the dilogarithm zero $w(0,-1)$.
\begin{maind}  Let $\overline{N}$ denote $N \bmod 2$. With
\begin{equation}\label{d0c}
    d_{0}\bigl(\overline{N}\bigr) = z_0 \sqrt{2  e^{-\pi i z_0}\bigl(e^{-\pi i z_0}+(-1)^N \bigr)}
\end{equation}
 and  explicit  $d_{1}\bigl(\sigma, \overline{N}\bigr) $, $d_{2}\bigl(\sigma, \overline{N}\bigr), \dots $ depending on $\sigma \in \R$ and $\overline{N}$, we have
\begin{equation} \label{presd}
   \mathcal D_1(N,\sigma) = \Re\left[\frac{w_0^{-N/2}}{N^{2}} \left( d_{0}\bigl(\overline{N}\bigr) +\frac{d_{1}\bigl(\sigma, \overline{N}\bigr)}{N}+ \dots +\frac{d_{m-1}\bigl(\sigma, \overline{N}\bigr)}{N^{m-1}}\right)\right] + O\left(\frac{|w_0|^{-N/2}}{N^{m+2}}\right)
\end{equation}
for an implied constant depending only on  $\sigma$ and $m$.
\end{maind}
\begin{proof}
Let $v=\Im(z_0)/\Re(z_0) \approx 0.216279$ and $c=1+iv$. We replace the path of integration $[1.01,1.49]$ in \eqref{save} and \eqref{save1} with the path $\mathcal P$ through $z_0$ made up of the lines joining $1.01$, $1.01c$, $1.49c$ and $1.49$. This path is used in \cite[Sect. 5.2]{OS1} and it is proved there that $\Re(p(z)-p(z_0))>0$ for $z\in \mathcal P$ except at $z=z_0$.

For $N$ odd, applying the saddle-point method to \eqref{save},  as in Theorem \ref{c2se}, gives
\begin{equation*}
    \mathcal D_3(N,\sigma)  = \Re \left[  e^{-N p(z_0)/2}
    \sum_{t=0}^{d-2} \frac{-2}{N^{t+2}}    \sum_{s=0}^{t} \G(s+1/2) a_{2s}(q_{\mathcal D} \cdot u_{\mathcal D, \sigma, t-s})
    \right]+ O\left( \frac{|w_0|^{-N/2}}{N^{d+1}}\right).
\end{equation*}
Therefore we set
\begin{equation} \label{dtyso}
    d_t\bigl(\sigma, \overline{N}\bigr):=   -2\sum_{s=0}^t \G(s+1/2) a_{2s}(q_{\mathcal D} \cdot u_{\mathcal D, \sigma, t-s}) \qquad (N \text{ odd}).
\end{equation}
Since $\sqrt{w_0}  = e^{p(z_0)/2}$ and  \eqref{cosdk} is true, we obtain \eqref{presd} in the statement of the theorem in this odd case.

For $N$ even, \eqref{save1} implies
\begin{equation*}
    \mathcal D_3(N,\sigma) = \frac{-1}{ (N+1)^{3/2}} \Re \int_{1.01}^{1.49}
     \exp \bigl(-(N+1) \cdot p(z)/2 \bigr) i q^*_\mathcal D(z) \exp \bigl(v_\mathcal D (z;N+1,\sigma)\bigr) \, dz  \qquad (N \text{ even})
\end{equation*}
and
applying the saddle-point method  yields
\begin{equation*}
    \mathcal D_3(N,\sigma)  = \Re \left[  e^{-(N+1) p(z_0)/2}
    \sum_{t=0}^{d-2} \frac{-2}{(N+1)^{t+2}}    \sum_{s=0}^{t} \G(s+1/2) a_{2s}(i q^*_{\mathcal D} \cdot u_{\mathcal D, \sigma, t-s})
    \right]+ O\left( \frac{|w_0|^{-N/2}}{N^{d+1}}\right).
\end{equation*}
Define
\begin{equation} \label{dtyse1}
    d_t^*(\sigma):=   -2 e^{-p(z_0)/2} \sum_{s=0}^t \G(s+1/2) a_{2s}(i q^*_{\mathcal D} \cdot u_{\mathcal D, \sigma, t-s})
\end{equation}
and we want to convert the above series in $1/(N+1)$ to one in $1/N$.
The method to do this is given in the proof of Theorem \ref{c2ses}. Let
\begin{equation} \label{dtyse2}
    d_t\bigl(\sigma, \overline{N}\bigr):= \sum_{j=0}^t (-1)^{t-j}\binom{t+1}{j+1}  d_j^*(\sigma) \qquad (N \text{ even})
\end{equation}
and with \eqref{cosdk} we obtain \eqref{presd} in the statement of the theorem in this even case.

To calculate $d_0\bigl(\sigma, \overline{N}\bigr)$, we begin with $N$ odd and see from  \eqref{dtyso} and \eqref{a2sb} that
\begin{equation*}
    d_0\bigl(\sigma, \overline{N}\bigr) =  -2\sqrt{\pi} a_0(q_{\mathcal D} \cdot 1) = -2\sqrt{\pi} \frac{\omega}{2(\omega^2 p_0/2)^{1/2}} q_0
\end{equation*}
for $q_0=q_\mathcal D(z_0)$, $p_0=p''(z_0)/2$ and the direction $\omega=z_0$. Short computations (see \eqref{pzzz22}) provide
\begin{equation*}
    p_0  =\frac{\pi i}{z_0(1-e^{-2\pi i z_0})}, \qquad q_0^2  =\frac{ -i z_0 e^{-\pi i z_0}}{e^{-\pi i z_0}+1}
\end{equation*}
so that
\begin{equation*}
    d_0\bigl(\sigma, \overline{N}\bigr)^2 =  2z_0^2 e^{-\pi i z_0}\bigl(e^{-\pi i z_0}-1 \bigr)  \qquad (N \text{ odd})
\end{equation*}
and \eqref{d0c} follows in this case. The $N$ even case is similar: from  \eqref{dtyse1}, \eqref{dtyse2} and \eqref{a2sb}
\begin{equation*}
    d_0\bigl(\sigma, \overline{N}\bigr) =  -2 e^{-p(z_0)/2}  \sqrt{\pi} a_0(i q^*_{\mathcal D} \cdot 1) = -2 w_0^{-1/2} \sqrt{\pi} \frac{\omega}{2(\omega^2 p_0/2)^{1/2}} i q^*_0
\end{equation*}
for $q^*_0=q^*_\mathcal D(z_0)$. We see that $\left(q_0^*\right)^2  = i z_0 (e^{\pi i z_0}+1)$ and so
\begin{equation*}
    d_0\bigl(\sigma, \overline{N}\bigr)^2 =  2z_0^2 e^{-\pi i z_0}\bigl(e^{-\pi i z_0}+1 \bigr)  \qquad (N \text{ even})
\end{equation*}
and \eqref{d0c} follows in this case also.
\end{proof}

Table \ref{d1n1} gives an example of the accuracy of \eqref{presd} in Theorem \ref{c2sed}.
\begin{table}[h]
\begin{center}
\begin{tabular}{c|cccc|c}
$N$ & $m=1$ & $m=2$ & $m=3$ & $m=4$ & $\mathcal D_1(N,1)$  \\ \hline
$1000$ & $-1.7713 \times 10^{9}$ &  $-1.7785 \times 10^{9}$ & $-1.7778 \times 10^{9}$ &  $-1.77778 \times 10^{9}$ &  $-1.77778 \times 10^{9}$ \\
$1001$ & $-2.10996 \times 10^{9}$ &  $-2.11483 \times 10^{9}$ & $-2.1142 \times 10^{9}$ &  $-2.11418 \times 10^{9}$ &  $-2.11418 \times 10^{9}$
\end{tabular}
\caption{Theorem \ref{c2sed}'s approximations to $\mathcal D_1(N,1)$.} \label{d1n1}
\end{center}
\end{table}

\section{The sum $\mathcal E_1(N, \sigma)$} \label{sec-e1}

Let $\sigma \in \Z$. In this section  we prove Theorem \ref{thme}, giving the asymptotic expansion as $N \to \infty$ of
\begin{equation*}
    \mathcal E_1(N,\sigma) := \sum_{h/k \in \mathcal E(N)} Q_{hk\sigma}(N) = 2 \Re \sum_{ \frac{N}{3}  <k \lqs \frac{N}{2}}  Q_{1k\sigma}(N).
\end{equation*}

\subsection{Higher-order poles}

Recall from \eqref{qhkl} that
\begin{equation*}
  Q_{hk\sigma}(N) :=  2\pi i \res_{z=h/k} \frac{e^{2\pi i \sigma z}}{(1-e^{2\pi i z})(1-e^{2\pi i2 z}) \cdots (1-e^{2\pi i N z})}
\end{equation*}
and  the expression on the right above has a pole at $z=h/k$ of order $s=\lfloor N/k \rfloor$. We calculated $Q_{hk\sigma}(N)$ in the case of a simple pole ($s=1$ or equivalently $N/2<k\lqs N$) in Proposition \ref{simple}  and require the double pole case ($s=2$ or $N/3<k\lqs N/2$)  in this section. In general, we have
\begin{equation*}
    e^{2\pi i \sigma z} = e^{2\pi i \sigma h/k} \sum_{r=0}^\infty \frac{(2\pi i \sigma)^r}{r!} (z-h/k)^r
\end{equation*}
and
for $m \in \Z_{\neq 0}$ write
$$
\frac{1}{1-e^{2\pi i m z}} = \sum_{r=0}^\infty \frac{\beta_r(m,h/k)}{r!}(z-h/k)^r \times
\begin{cases} (z-h/k)^{-1} & \text{ if \ } k  \mid m\\
 1 & \text{ if \ } k \nmid m .
\end{cases}
$$
Therefore, for any $k$,
\begin{equation}\label{cbet}
    Q_{hk\sigma}(N)=2\pi i \cdot e^{2\pi i \sigma h/k} \sum_{r_0+r_1+ \cdots + r_N = s-1} (2\pi i \sigma)^{r_0} \frac{\beta_{r_1}(1,h/k) \beta_{r_2}(2,h/k)  \cdots \beta_{r_N}(N,h/k)}{r_0! r_1!  \cdots r_N!}
\end{equation}
where
\begin{alignat}{2}
    \beta_r(m,h/k) & =  -(2\pi i m)^{r-1} B_r \qquad & & (k \mid m), \notag \\
   \beta_r(m,h/k) & = \left. \frac{d^r}{dz^r} \frac{1}{1-e^{2\pi i m z}} \right|_{z=h/k} \qquad & & (k \nmid m). \label{bet}
\end{alignat}
Formula \eqref{bet} implies for example,
\begin{eqnarray*}
  \beta_0(m,h/k) &=& \frac{1}{1-e^{2\pi i m h/k}} \\
  \beta_1(m,h/k) &=& -2\pi i m \beta_0(m,h/k)\Bigl(1-\beta_0(m,h/k)\Bigr) \\
  \beta_2(m,h/k) &=& (2\pi i m)^2 \beta_0(m,h/k)\Bigl(1-3\beta_0(m,h/k)+2\beta_0(m,h/k)^2\Bigr)
\end{eqnarray*}
for $k \nmid m$.

\subsection{Second-order poles}
For  $N/3<k\lqs N/2$ (and $s=2$), formula \eqref{cbet} shows that
\begin{equation*}
    Q_{hk\sigma}(N)=2\pi i \cdot e^{2\pi i \sigma h/k} \left[2\pi i \sigma + \frac{\beta_1(1,h/k)}{\beta_0(1,h/k)}+ \cdots + \frac{\beta_1(N,h/k)}{\beta_0(N,h/k)}\right] \prod_{j=1}^N \beta_0(j,h/k)
\end{equation*}
and hence, recalling the root of unity identity after \eqref{i2},
$$
Q_{hk\sigma}(N)=\frac{- e^{2\pi i \sigma h/k}}{2k^4}\left(\frac{N(N+1)-3k-2\sigma}{2}-\sum_{1\lqs m\lqs N,\ k \nmid m} \frac{m}{1-e^{2\pi i m h/k}} \right)\prod_{j=1}^{N-2k}\frac{1}{1-e^{2\pi i h j/k}}.
$$
For the case we need, $h=1$,
\begin{align*}
    \sum_{1\lqs m\lqs N,\ k \nmid m} \frac{m}{e^{2\pi i m/k}-1} & = \sum_{1\lqs m\lqs N,\ k \nmid m} \frac{m \cdot e^{-\pi i m/k}}{e^{\pi i m/k}-e^{-\pi i m/k}} \\
    & = \frac{1}{2i}\sum_{1\lqs m\lqs N,\ k \nmid m} \frac{m (\cos(-\pi  m/k)+i\sin(-\pi  m/k))}{\sin(\pi  m/k)}\\
    & = \frac{-1}{2i}\sum_{1\lqs m\lqs N,\ k \nmid m} i m + \frac{1}{2i}\sum_{1\lqs m\lqs N,\ k \nmid m} m \cot(\pi m/k).
\end{align*}
Therefore
\begin{equation} \label{cyan}
     Q_{1k\sigma}(N) =\frac{1}{2k^2} \phi(N,k,\sigma)
\exp\left( N\left[ \frac{-i\pi}{2}\left(\frac{N}{k}-1 +2\frac k N \right)\right]\right)
\exp\left( \frac{-i\pi}2 \frac{N}{k}\right)
\exp\left( \frac 1N\left[ 2 i\pi \sigma\frac{N}{k}\right]\right)
 \spr{1/k}{N-2k}
\end{equation}
for
\begin{equation} \label{pink}
    \phi(N,k,\sigma) := \frac{1}{4k^2}(N^2+N-4\sigma) + \frac{1}{2\pi ik}\sum_{1\lqs j\lqs N,\ k \nmid j} \frac{\pi j}{k} \cot \left(\frac{\pi j}{k}\right).
\end{equation}
Also note that
\begin{equation}\label{blue}
    \left| Q_{1k\sigma}(N) \right| = \frac{1}{2k^2} \left| \phi(N,k,\sigma)
\cdot  \spr{1/k}{N-2k} \right|.
\end{equation}

\begin{lemma} \label{phib}
For $N/3<k\lqs N/2$ and an  implied constant depending only on $\sigma$
\begin{equation*}
    \phi(N,k,\sigma)=O(N).
\end{equation*}
\end{lemma}
\begin{proof}
Verify that $1/|\sin(\pi j/k)| < 2k/\pi$ for $k \nmid j$ (as in \cite[Sect. 3.3]{OS1}). Therefore
\begin{equation*}
    |\cot(\pi j/k)| < 2k/\pi \qquad (k \nmid j)
\end{equation*}
and the lemma follows.
\end{proof}

Set $z=z(N,k):=N/k$. Applications of Propositions \ref{gprop} and  \ref{hard}, with $m=N-2k$ and  $s=N/2$, prove the following.
\begin{theorem} \label{cplxd2}
Fix $W>0$. Let $\Delta$ be in the range $0.0048 \lqs \Delta \lqs 0.0079$ and set $\alpha = \Delta \pi e$.  Suppose $\delta$ and $\delta'$ satisfy
\begin{equation*}
    \frac{\Delta}{1-\Delta} < \delta  \lqs \frac{1}{e}, \ 0<\delta' \lqs \frac{1}{e} \quad \text{ and  } \quad \delta \log 1/\delta, \ \ \delta' \log 1/\delta' \lqs W.
\end{equation*}
Then for all $N \gqs 2 \cdot R_\Delta$  we have
\begin{equation}\label{ott1d2}
    \spr{1/k}{N-2k} = O\left(e^{WN/2}\right) \quad \text{ for } \quad  z \in \left[2, \ 2+ \delta \right] \cup  \left[ 5/2 - \delta' , \ 3 \right)
\end{equation}
and
\begin{multline} \label{ott2d2}
    \spr{1/k}{N-2k} =  \frac{1}{N^{1/2}}  \exp\left(N\frac {\cl(2\pi  z)}{2\pi z}  \right) \left( \frac{z}{2\sin(\pi (z-2))}\right)^{1/2}
      \\
      \quad \times \exp\left(\sum_{\ell=1}^{L-1} \frac{g_{\ell}(z)}{N^{2\ell-1}} \right) +O\left(e^{WN/2}\right)
        \quad \text{ for } \quad z \in  \left(2+ \delta, \ 5/2 - \delta' \right)
  \end{multline}
with $L=\lfloor \alpha \cdot N/2\rfloor$. The implied constants in \eqref{ott1d2}, \eqref{ott2d2} are absolute.
\end{theorem}

\subsection{Estimating $\phi(N,k,\sigma)$}
With Lemma \ref{phib} and \eqref{ott1d2}, we see that
\begin{equation} \label{rst9}
    \mathcal E_1(N, \sigma) = 2 \Re \sum_{ k \ : \ z \in  \left(2+ \delta, \ 5/2 - \delta' \right)}  Q_{1k\sigma}(N) +O(N e^{WN/2})
\end{equation}
 and so we may restrict our attention to  indices $k$ corresponding to this range. Let
$$
f(x):=x \cot(x),
$$
a smooth function of $x\in \R$ except at $x=\pm \pi, \pm 2\pi \dots$ and with $f(0)=1$. Note the identities
\begin{equation*}
   f(-x)=f(x), \quad f(\pi+x)=f(x)+\pi \cot(x), \quad f(\pi-x)=f(x)-\pi \cot(x)
\end{equation*}
for example. Let $m=N-2k$ as before, so that $0 \lqs m <k$. With \eqref{rst9} we may assume
\begin{equation*}
    \delta k < m <k/2 -\delta' k,
\end{equation*}
and in particular, $m \neq 0$. For $m<k/2$, the sum we need from \eqref{pink} is
\begin{multline} \label{es9}
    \sum_{1\lqs j\lqs N,\ k \nmid j} f\left(\frac{\pi j}k\right)
     = \sum_{m < j <k-m} \left( f\left(\frac{\pi j}k\right) + f\left(\frac{\pi (k+j)}k\right) \right) \\
     +\sum_{1\lqs j \lqs m} \left( f\left(\frac{\pi j}k\right) + f\left(\frac{\pi (k-j)}k\right) +f\left(\frac{\pi (k+j)}k\right)  +f\left(\frac{\pi (2k-j)}k\right) +f\left(\frac{\pi (2k+j)}k\right) \right)  \\
     = 5\sum_{1\lqs j \lqs m}  f\left(\frac{\pi j}k\right) + 2 \sum_{m < j <k-m}  f\left(\frac{\pi j}k\right).
\end{multline}
With $\rho(z):=\log \bigl((\sin z)/z\bigr)$, we have
\begin{equation*}
    f(x)=1+x \rho'(x)
\end{equation*}
and for $d \in \Z_{\gqs 1}$
\begin{align}
    f^{(d)}(x) & = x \cot^{(d)}(x)+d \cot^{(d-1)}(x) \label{frh}\\
     & = x \rho^{(d+1)}(x)+d \rho^{(d)}(x). \label{frh2}
\end{align}
Since $\rho^{(d)}(0)$ equals $0$ for $d$ odd, (and equals $-2^d|B_d|/d$ for $d$ even), we see
\begin{equation} \label{fm0}
    f^{(d)}(0)=0 \qquad (d \text{ odd}).
\end{equation}
Also note the relation
\begin{equation} \label{fm00}
    f^{(d)}(\pi-x)= (-1)^d \left( f^{(d)}(x) - \pi \cot^{(d)}(x) \right).
\end{equation}
Applying Euler-Maclaurin summation to \eqref{es9}, as in \cite[Chap. 2]{Ra} or \cite[p. 285]{Ol}, and simplifying with \eqref{fm0}, \eqref{fm00} produces
\begin{multline} \label{ess}
    \sum_{1\lqs j\lqs N,\ k \nmid j} f\left(\frac{\pi j}k\right) = 5 \int_0^m f\left(\frac{\pi x}k\right) \, dx + 2 \int_m^{k-m} f\left(\frac{\pi x}k\right) \, dx - \frac{5}{2} + \frac{1}{2} f\left(\frac{\pi m}k\right) + \pi \cot\left(\frac{\pi m}k\right) \\
    + \sum_{\ell =1}^{L-1} \frac{B_{2\ell }}{(2\ell )!} \left( \frac{\pi }{k}\right)^{2\ell -1}\left\{ f^{(2\ell -1)}\left(\frac{\pi m}{k}\right) +
    2\pi \cot^{(2\ell -1)}\left(\frac{\pi m}{k}\right) \right\} + \varepsilon_L(m,1/k)
\end{multline}
for
\begin{equation}\label{vek}
\varepsilon_L(m,1/k) :=\left( \frac{\pi }{k}\right)^{2L} \left[5\int_0^m + 2\int_m^{k-m}\right]\frac{B_{2L}-B_{2L}(x-\lfloor x \rfloor)}{(2L)!} f^{(2L)}\left(\frac{\pi x }{k}\right) \, dx.
\end{equation}
With the evaluation
\begin{equation*}
    \int_0^t x \cot(x) \, dx = \frac 12 \cl(2t)-t \clp(2t)
\end{equation*}
we find
\begin{equation*}
    5 \int_0^m f\left(\frac{\pi x}k\right) \, dx + 2 \int_m^{k-m} f\left(\frac{\pi x}k\right) \, dx = \frac k{2\pi} \cl(2\pi m/k)- N \clp(2\pi m/k).
\end{equation*}
Using \eqref{frh} also, \eqref{ess} becomes
\begin{multline} \label{ess2}
    \sum_{1\lqs j\lqs N,\ k \nmid j} f\left(\frac{\pi j}k\right) = \frac k{2\pi} \cl(2\pi m/k)- N \clp(2\pi m/k)
     - \frac{5}{2} + \frac{\pi N}{2 k}  \cot\left(\frac{\pi m}k\right) \\
    + \sum_{\ell =1}^{L-1} \frac{B_{2\ell }}{(2\ell )!} \left( \frac{\pi }{k}\right)^{2\ell -1}\left\{ \frac{\pi N}{k} \cot^{(2\ell -1)}\left(\frac{\pi m}{k}\right) +
    (2 \ell -1) \cot^{(2\ell -2)}\left(\frac{\pi m}{k}\right) \right\} + \varepsilon_L(m,1/k).
\end{multline}
Define
\begin{equation*}
    \tilde{g}_\ell(z) := \frac{B_{2\ell }}{(2\ell )!} \left( \pi z\right)^{2\ell -1}\left\{ \pi z \cot^{(2\ell -1)}\left(\pi z\right) +
    (2 \ell -1) \cot^{(2\ell -2)}\left(\pi z\right) \right\}.
\end{equation*}
With \eqref{pink} and \eqref{ess2} we have demonstrated that
\begin{multline} \label{sunny}
    \phi(N,k,\sigma) = \left[\frac {\cl(2\pi z)}{4\pi^2 i} - \frac {z \clp(2\pi z)}{2\pi i} + \frac{z^2}{4}\right]
    +\frac{1}{N} \left[ \frac{z^2 \cot(\pi z)}{4i} + \frac{z^2}{4} - \frac{5z}{4\pi i}\right]\\
    - \frac{\sigma z^2}{N^2} + \frac{z}{2\pi i} \sum_{\ell =1}^{L-1} \frac{\tilde{g}_\ell(z)}{N^{2\ell}} + \frac{ \varepsilon_L(m,1/k) }{2\pi i k}
\end{multline}
which we write as
\begin{equation*}
    \phi(N,k,\sigma) = \sum_{\ell=0}^{2L-1} \frac{\phi_{\sigma, \ell}(z)}{N^{\ell}} + \frac{ \varepsilon_L(m,1/k) }{2\pi i k},
\end{equation*}
though only $\phi_{\sigma, 2}(z)$ depends on $\sigma$.

\begin{prop} \label{vrc}
For $1 \lqs m \lqs k/2$ we have
\begin{equation*}
    \frac{\left| \varepsilon_L(m,1/k) \right|}{2\pi k} \lqs 2\pi^2 (2L-1) \left( \frac{2L-1}{2\pi e m}\right)^{2L-1}.
\end{equation*}
\end{prop}
\begin{proof}
The arguments here are similar to those in \cite[Sect. 3]{OS1}. Use the inequalities
\begin{equation*}
    |B_{2n}-B_{2n}(x-\lfloor x \rfloor)|  \lqs  2|B_{2n}|, \qquad \frac{|B_{2n}|}{(2n)!}  \lqs  \frac{\pi^2}{3(2\pi)^{2n}}
\end{equation*}
from \cite[Thm 1.1, p. 283]{Ol} and \cite[(9.6)]{Ra} to see that
\begin{equation} \label{nyt}
    \left| \varepsilon_L(m,1/k) \right| \lqs \frac{2 \pi^2}{3(2\pi)^{2L}}\left( \frac{\pi }{k}\right)^{2L} \left[5\int_0^m + 2\int_m^{k-m}\right]\left| f^{(2L)}\left(\frac{\pi x }{k}\right) \right| \, dx.
\end{equation}
By \cite[(11.1)]{Ra}
\begin{equation} \label{simo2}
    -\rho'(w) = \sum_{r=1}^\infty \frac{2^{2r} |B_{2r}|}{(2r)!} w^{2r-1} \qquad (|w|<\pi)
\end{equation}
so that
\begin{equation*}
    \rho^{(d)}(x) \lqs 0 \quad \text{ for all }\quad x\in [0,\pi), \ d \in \Z_{\gqs 0}.
\end{equation*}
Hence  \eqref{frh2} implies
\begin{equation*}
    f^{(d)}(x) \lqs 0 \quad \text{ for all }\quad x\in [0,\pi), \ d \in \Z_{\gqs 1}
\end{equation*}
and $\left| f^{(2L)}\left(\pi x /k\right) \right| = - f^{(2L)}\left(\pi x /k\right)$ in \eqref{nyt}. On integrating and applying \eqref{fm0}, \eqref{fm00} we obtain
\begin{align*}
    \left| \varepsilon_L(m,1/k) \right| & \lqs -\frac{\pi}{3}\left( \frac{1 }{2k}\right)^{2L-1} \left( f^{(2L-1)}\left(\frac{\pi m }{k}\right) + 2\pi \cot^{(2L-1)}\left(\frac{\pi m }{k}\right)  \right) \\
    & = \frac{\pi}{3}\left( \frac{1 }{2k}\right)^{2L-1} \left(
    2\pi (2L-1)! \left(\frac{k}{\pi m }\right)^{2L}
    - \frac{\pi N }{k}
    \rho^{(2L)}\left(\frac{\pi m }{k}\right) - (2L-1) \rho^{(2L-1)}\left(\frac{\pi m }{k}\right)  \right)
\end{align*}
with the last line coming from \eqref{frh2} and the further identity
\begin{equation*}
    \cot^{(d)}(x) = \frac{(-1)^d d!}{x^{d+1}}+\rho^{(d+1)}(x).
\end{equation*}
Use
\begin{equation*}
    \left| \rho^{(d+1)}(x) \right| \lqs \frac{2\pi d!}{3}  \left( \frac 2 \pi \right)^d \qquad (|x| \lqs \pi/2, \ d \in \Z_{\gqs 0})
\end{equation*}
from \eqref{simo2}, and $(n-1)!  <  3 \left( n / e\right)^n$ from Stirling's formula, to complete the proof.
\end{proof}

\subsection{Approximating $\mathcal E_1(N, \sigma)$}

With \eqref{dli} for $m=2$ we find, for $2<z<3$,
\begin{equation*}
    \frac {\cl(2\pi  z)}{2\pi z} - \frac{i \pi}{2}(z-1+2/z) =-2\pi i +\frac 1{2\pi i z} \Bigl[ \li(e^{2\pi i z}) -\li(1) -4\pi^2 \Bigr].
\end{equation*}
Put
\begin{align}
    r_{\mathcal E}(z)  & :=  \frac 1{2\pi i z} \Bigl[ \li(e^{2\pi i z}) -\li(1) -4\pi^2 \Bigr] \label{rdx}\\
    q_{\mathcal E}(z;N, \sigma) & :=   \left(\frac{z}{2\sin(\pi z)} \right)^{1/2}\exp\left( \frac{-\pi i z}2 \right) \times  \sum_{\ell=0}^{2L-1} \frac{\phi_{\sigma, \ell}(z)}{N^{\ell}}  \label{qdx}\\
 v_\mathcal E(z;N, \sigma) & :=  \frac{2\pi i \sigma z}{N} + \sum_{\ell=1}^{L-1} \frac{g_{\ell}(z)}{N^{2\ell-1}}  \label{vdx}
\end{align}
for $L := \lfloor \alpha \cdot N/2\rfloor$ in \eqref{qdx} and \eqref{vdx}. Also set
\begin{equation*}
    \mathcal E_2(N, \sigma)  := \frac{1}{N^{1/2}} \Re \sum_{k \ : \ z \in  \left(2+ \delta, \ 5/2 - \delta' \right)}
    \frac{1}{k^2} \exp \bigl(N \cdot r_\mathcal E(z) \bigr) q_{\mathcal E}(z;N, \sigma) \exp \bigl(v_\mathcal E(z;N, \sigma)\bigr) .
\end{equation*}

The terms summed for $\mathcal E_2(N, \sigma)$ above differ from the terms in $\mathcal E_1(N, \sigma)$  only in the removal of the error terms from the approximations of $\spr{1/k}{N-2k}$ and $\phi(N,k,\sigma)$. The next proposition lets us control what happens on removing the error term for $\phi(N,k,\sigma)$.

\begin{prop} \label{delad}
Suppose   $\Delta$ and $W$ satisfy $0.0048 \lqs \Delta \lqs 0.0079$ and $\Delta \log 1/\Delta \lqs W$.
For the integers $k$, $s$ and $m$ we require
\begin{equation*}
    1<k \lqs s, \quad  R_\Delta \lqs s, \quad \Delta s \lqs m \lqs k/2.
\end{equation*}
 Then for $L:=\lfloor \pi e \Delta \cdot s \rfloor$ we have
\begin{align*}
 \spr{1/k}{m} \frac{ \varepsilon_L(m,1/k) }{2\pi i k}  & =O(s e^{sW}) \\
    \frac{ \varepsilon_L(m,1/k) }{2\pi i k}  & =O(s) .
\end{align*}
\end{prop}
\begin{proof}
We may copy the proof of Proposition \ref{dela} in \cite[Sect. 3.4]{OS1}. The bound used for $T_L(m,h/k)$ in that result is $\left( \frac{2L-1}{2\pi e m}\right)^{2L-1}$. The corresponding bound for $\varepsilon_L(m,1/k)/(2\pi i k)$ in Proposition \ref{vrc} is bigger by a factor $2L-1 \ll s$.
\end{proof}

Choosing $s=N/2$ and $m=N-2k$ in Proposition \ref{delad} shows
\begin{align}
 \spr{1/k}{N-2k} \frac{ \varepsilon_L(m,1/k) }{2\pi i k}  & = O(N e^{WN/2}) \label{slmad}\\
    \frac{ \varepsilon_L(m,1/k) }{2\pi i k}  & = O(N) \label{slmbd}
\end{align}
for $N \gqs 2 \cdot R_\Delta$, $L= \lfloor \alpha \cdot N/2\rfloor$ and $2+\Delta/(1-\Delta/2) \lqs z \lqs 5/2$.

\begin{prop} \label{e1e2n}
For an implied constant depending only on $\sigma$
$$
\mathcal E_1(N, \sigma) = \mathcal E_2(N, \sigma) +O(N e^{WN/2}).
$$
\end{prop}
\begin{proof}
Starting with \eqref{rst9}, write
\begin{equation*}
 \sum_{ k \ : \ z \in  \left(2+ \delta, \ 5/2 - \delta' \right)}  Q_{1k\sigma}(N)  =  \sum_{ k \ : \ z \in  \left(2+ \delta, \ 5/2 - \delta' \right)}  \frac{Q_{1k\sigma}(N)}{\phi(N,k,\sigma)}\left( \sum_{\ell=0}^{2L-1} \frac{\phi_{\sigma, \ell}(z)}{N^{\ell}} + \frac{ \varepsilon_L(m,1/k) }{2\pi i k}\right)
\end{equation*}
where
\begin{equation*}
    \frac{Q_{1k\sigma}(N)}{\phi(N,k,\sigma)} = \frac{1}{2k^2}
\exp\left( N  \frac{-i\pi (z-1 +2/z )}{2}-\frac{\pi i z}2  + \frac {2 \pi i \sigma z}N\right)
 \spr{1/k}{N-2k}
\end{equation*}
by \eqref{cyan}.
We have
\begin{multline*}
     \sum_{ k \ : \ z \in  \left(2+ \delta, \ 5/2 - \delta' \right)}  \frac{Q_{1k\sigma}(N)}{\phi(N,k,\sigma)} \frac{ \varepsilon_L(m,1/k) }{2\pi i k}\\
    \ll \sum_{ k \ : \ z \in  \left(2+ \delta, \ 5/2 - \delta' \right)}  \frac{1}{k^2} \left| \spr{1/k}{N-2k} \frac{ \varepsilon_L(m,1/k) }{2\pi i k} \right|
    \ll N e^{WN/2}
\end{multline*}
using \eqref{slmad} and that
\begin{equation*}
    \frac{\Delta}{1-\Delta/2} < \frac{\Delta}{1-\Delta} < \delta
\end{equation*}
so  the bound \eqref{slmad} is valid for $z \in  (2+ \delta, \ 5/2 - \delta')$.
Therefore,
\begin{equation} \label{dlksirj}
    \mathcal E_1(N, \sigma) = 2 \Re \sum_{ k \ : \ z \in  \left(2+ \delta, \ 5/2 - \delta' \right)}  \frac{Q_{1k\sigma}(N)}{\phi(N,k,\sigma)}\left( \sum_{\ell=0}^{2L-1} \frac{\phi_{\sigma, \ell}(z)}{N^{\ell}}\right) +O(N e^{WN/2}).
\end{equation}
Next note that
\begin{equation}
   \left| \sum_{\ell=0}^{2L-1} \frac{\phi_{\sigma, \ell}(z)}{N^{\ell}} \right|  \lqs \left| \phi(N,k,\sigma)\right| + \left|\frac{ \varepsilon_L(m,1/k) }{2\pi i k} \right| \ll N \label{yestag}
\end{equation}
by Lemma \ref{phib} and \eqref{slmbd}. With \eqref{yestag} we see that replacing $\spr{1/k}{N-2k}$ in \eqref{dlksirj} by the main term on the right of  \eqref{ott2d2} changes $\mathcal E_1(N, \sigma)$ by at most $O(N e^{WN/2})$, as required.
\end{proof}

Comparing \eqref{rdx}-\eqref{vdx} and \eqref{rbw}-\eqref{vbw} gives the relations
\begin{equation*}
    r_{\mathcal E}(z)   =  r_{\mathcal C}(z), \qquad
    q_{\mathcal E}(z;N, \sigma)  =  q_{\mathcal C}(z)  \sum_{\ell=0}^{2L-1} \frac{\phi_{\sigma, \ell}(z)}{N^{\ell}}, \qquad
 v_\mathcal E(z;N, \sigma)  =  v_\mathcal C(z;N, \sigma)
\end{equation*}
so that we may reuse our work from Section \ref{c1ns}. We fix the choice of constants as in \eqref{fix}.

\begin{lemma} \label{phot}
The function $q_{\mathcal E}(z;N, \sigma)$ is holomorphic for $2< \Re(z)<5/2$. In the box with $2+\delta \lqs \Re(z) \lqs 5/2-\delta'$ and $-1\lqs \Im(z) \lqs 1$,
\begin{equation} \label{efex3}
    q_{\mathcal E}(z;N, \sigma) \ll 1
\end{equation}
for an  implied constant depending only on $\sigma \in \R$.
\end{lemma}
\begin{proof}
The first issue is that $\phi_{\sigma, 0}(z)$ has only been defined in \eqref{sunny} for $z \in \R$. Use \eqref{dli} and its derivative with $m=2$ to show
\begin{equation}\label{p0z}
    \phi_{\sigma, 0}(z)=\frac 1{4\pi^2} \left[\li(1)- \li(e^{2\pi i z}) +6\pi^2 -2\pi i z \log(1-e^{2\pi i z}) \right]
\end{equation}
giving the analytic continuation of $\phi_{\sigma, 0}(z)$ to all $z$ with $2< \Re(z)<5/2$. It follows, as in Theorem \ref{rzjb}, that $q_{\mathcal E}(z;N, \sigma)$ is holomorphic in $z$ as required. The bound \eqref{efex3} follows from
\begin{equation*} 
    \frac{\tilde g_{\ell}(z)}{N^{2\ell-1}} \ll N \bigl(F_{N,\varepsilon}(2\ell-1)+ F_{N,\varepsilon}(2\ell)\bigr) e^{-\pi|y|},
 \end{equation*}
 with $F_{N,\varepsilon}$ defined in \eqref{hoot}, as in Proposition \ref{addx} and Corollary \ref{acdcx}.
\end{proof}

With the rectangle $C_1$ from Figure \ref{c1c2} we find
\begin{equation*}
    \mathcal E_2(N, \sigma) = \frac{-1}{N^{3/2}} \Re
      \int_{C_1}
     \exp \bigl(N \cdot r_\mathcal C(z)\bigr) \frac{q_{\mathcal E}(z;N, \sigma)}{2i\tan \bigl(\pi N/z \bigr)} \exp \bigl(v_\mathcal C(z;N,\sigma)\bigr)\, dz
\end{equation*}
where
$$
\frac{1}{2i\tan(\pi N/z)} = \begin{cases}
1/2+\sum_{j \lqs -1}  e^{2\pi i j N/z} & \text{ \ if \ } \Im z >0 \\
-1/2-\sum_{j \gqs 1} e^{2\pi i j N/z} & \text{ \ if \ } \Im z <0.
\end{cases}
$$
The arguments of Propositions \ref{cpcmb}, \ref{lgst} and \ref{lgst2} now go through almost unchanged:
\begin{multline*}
    \mathcal E_2(N, \sigma)  = \frac{-1}{N^{3/2}} \Re \left[\sideset{}{'}{\sum}_{j \lqs 0}\int_{C_1^+}
     \exp \bigl(N [ r_\mathcal C(z) +2\pi i j/z]  \bigr)  q_{\mathcal E}(z;N, \sigma) \exp \bigl(v_\mathcal C(z;N,\sigma)\bigr) \, dz \right.\\
     -\left. \sideset{}{'}{\sum}_{j \gqs 0}\int_{C_1^-}
    \exp \bigl(N [ r_\mathcal C(z) +2\pi i j/z]  \bigr)  q_{\mathcal E}(z;N, \sigma) \exp \bigl(v_\mathcal C(z;N,\sigma)\bigr) \, dz \right] + O(e^{WN/2}),
\end{multline*}
the term with $j = -1$ is the largest and
\begin{equation}\label{cos2e}
    \mathcal E_2(N,\sigma)  = \mathcal E_3(N,\sigma) + O(e^{WN/2})
\end{equation}
for $W=0.05$, an implied constant depending only on $\sigma$, and
\begin{equation}\label{newee}
    \mathcal E_3(N,\sigma)  := \frac{1}{N^{3/2}}  \Re \int_{2.01}^{2.49}
     \exp \bigl(-N \cdot p(z) \bigr) q_{\mathcal E}(z;N, \sigma) \exp \bigl(v_\mathcal C(z;N,\sigma)\bigr) \, dz.
\end{equation}

\subsection{The asymptotic behavior of $\mathcal E_1(N, \sigma)$}

Arguing as in Lemma \ref{phot} shows the next result.

\begin{prop} \label{gase}
For $2.01 \lqs \Re(z) \lqs 2.49$ and $|\Im(z)|\lqs 1$, say, there is a holomorphic function $\xi_r(z;N,\sigma)$ of $z$ so that
\begin{equation*}
    q_{\mathcal E}(z;N, \sigma) = q_\mathcal C(z) \sum_{k=0}^{r-1} \frac{\phi_{\sigma, k}(z)}{N^{k}} + \xi_r(z;N,\sigma) \quad \text{for} \quad \xi_r(z;N,\sigma) = O\left(\frac{1}{N^r} \right)
\end{equation*}
with an implied constant depending only on $\sigma$ and  $r$ where $1 \lqs r \lqs 2L-1$ and $L=\lfloor 0.006 \pi e \cdot N/2 \rfloor$.
\end{prop}

We restate Theorem \ref{thme}:

\begin{maine}  With $e_{0}=-3z_1 e^{-\pi i z_1}/2$ and  explicit  $e_{1}(\sigma)$, $e_{2}(\sigma), \dots $ depending on $\sigma \in \R$ we have
\begin{equation} \label{prese}
   \mathcal E_1(N,\sigma) = \Re\left[\frac{w(0,-2)^{-N}}{N^{2}} \left( e_{0}+\frac{e_{1}(\sigma)}{N}+ \dots +\frac{e_{m-1}(\sigma)}{N^{m-1}}\right)\right] + O\left(\frac{|w(0,-2)|^{-N}}{N^{m+2}}\right)
\end{equation}
for an implied constant depending only on  $\sigma$ and $m$.
\end{maine}
\begin{proof}
With Propositions \ref{gasc} and  \ref{gase}, write
\begin{multline*}
    q_{\mathcal E}(z;N, \sigma) \exp\bigl(v_{\mathcal C}(z;N,\sigma)\bigr) = q_\mathcal C(z) \left( \sum_{k=0}^{r-1} \frac{\phi_{\sigma, k}(z)}{N^{k}}  \right)
    \left( \sum_{j=0}^{d-1} \frac{u_{\sigma, j}(z)}{N^j}\right) \\
    + q_{\mathcal E}(z;N, \sigma) \zeta_{d}(z;N,\sigma)
    + \xi_r(z;N,\sigma)  \exp\bigl(v_{\mathcal C}(z;N,\sigma)\bigr) - \xi_r(z;N,\sigma) \zeta_{d}(z;N,\sigma).
\end{multline*}
Then putting this into \eqref{newee} and moving the line of integration to $\mathcal Q$ (see Figure \ref{pthb}) gives
\begin{multline} \label{bcafe}
    \mathcal E_3(N, \sigma) =   \frac{1}{N^{3/2}}\Re \int_{\mathcal Q}
     \exp \bigl(-N \cdot p(z) \bigr) q_\mathcal C(z) \left( \sum_{k=0}^{r-1} \frac{\phi_{\sigma, k}(z)}{N^{k}}  \right)
    \left( \sum_{j=0}^{d-1} \frac{u_{\sigma, j}(z)}{N^j}\right) \, dz \\
     + O\left( \frac{|w(0,-2)|^{-N}}{N^{3/2}}\left(\frac 1{N^d} + \frac 1{N^r} +\frac 1{N^{d+r}}  \right)\right).
\end{multline}
The integral in \eqref{bcafe} is
\begin{equation*}
   \sum_{k=0}^{r-1} \sum_{j=0}^{d-1} \frac{1}{N^{3/2 +k +j}}  \int_{\mathcal Q}
      \exp \bigl(-N \cdot p(z) \bigr) q_\mathcal C(z) \phi_{\sigma, k}(z) u_{\sigma, j}(z) \, dz
\end{equation*}
 and applying the saddle-point method, Theorem \ref{sdle}, gives
\begin{equation*}
   \sum_{k=0}^{r-1} \sum_{j=0}^{d-1} \frac{2 e^{-N \cdot p(z_1)}}{N^{3/2 +k +j}}\left(\sum_{s=0}^{S-1} \Gamma(s+1/2) \frac{a_{2s}( q_\mathcal C \cdot \phi_{\sigma, k} \cdot u_{\sigma, j}  )}{N^{s+1/2}} + O\left(\frac 1{N^{S+1/2}}\right)\right).
\end{equation*}
 Letting $S=r=d$ we obtain, as in the proof of Theorem \ref{c2se},
\begin{equation} \label{hante}
    \mathcal E_3(N, \sigma) = \Re \left[  e^{-N \cdot p(z_1)}
    \sum_{t=0}^{d-2} \frac{2}{N^{t+2}}    \sum_{s=0}^{t} \sum_{k=0}^{t-s} \G(s+1/2) a_{2s}( q_\mathcal C \cdot \phi_{\sigma, k} \cdot   u_{\sigma, t-s-k})
    \right]+ O\left( \frac{|w(0,-2)|^{-N}}{N^{d+1}}\right).
\end{equation}
Hence, recalling Proposition \ref{e1e2n}, \eqref{cos2e} and with
\begin{equation} \label{ctyse}
    e_t(\sigma):=   2\sum_{s=0}^{t} \sum_{k=0}^{t-s} \G(s+1/2) a_{2s}( q_\mathcal C \cdot \phi_{\sigma, k} \cdot   u_{\sigma, t-s-k}),
\end{equation}
we obtain \eqref{prese} in the statement of the theorem.

Computing $e_0(\sigma)$ with \eqref{ctyse} gives
\begin{equation*}
    e_0(\sigma) = 2 \sqrt{\pi} a_0( q_\mathcal C \cdot \phi_{\sigma, 0} \cdot   1) = 2 \sqrt{\pi} \frac{\omega}{2(\omega^2 p_0)^{1/2}} q_{\mathcal C}(z_1) \phi_{\sigma, 0}(z_1).
\end{equation*}
With the identity
\begin{equation*}\label{pz0}
    2\pi i z^2 p'(z)=\li \left( e^{2\pi i z}\right)-\li(1) +2\pi i z \log \left(1-e^{2\pi i z}\right)
\end{equation*}
from \cite[Sect. 2.3]{OS1} we find that
\begin{equation*}
    \phi_{\sigma, 0}(z_1) = \frac{6\pi^2 - 2\pi i z_1^2 p'(z_1)}{4\pi^2} = \frac 32.
\end{equation*}
Combine this with the calculations in \eqref{oad2}
to get $e_0(\sigma)^2 = 9z_1^2 e^{-2\pi i z_1}/4$ and the formula for $e_0=e_0(\sigma)$ in the statement of the theorem follows.
\end{proof}

For example, a comparison of both sides of \eqref{prese} in Theorem \ref{thme} with $\sigma=1$ and some different values of $m$ and $N$  is shown in Table \ref{e1n1}.
\begin{table}[h]
\begin{center}
\begin{tabular}{c|cccc|c}
$N$ & $m=1$ & $m=2$ & $m=3$ & $m=4$ & $\mathcal E_1(N, 1)$  \\ \hline
$800$ & $879.611$ &  $905.272$ & $909.048$ &  $909.358$ &  $909.337$ \\
$1000$ & $-789369.$ &  $-784383.$ & $-784458.$ &  $-784480.$ &  $-784480.$
\end{tabular}
\caption{Theorem \ref{thme}'s approximations to $\mathcal E_1(N,1)$.} \label{e1n1}
\end{center}
\end{table}

\begin{proof}[\bf Proof of Theorem \ref{mainb}]
Recall the sets $\mathcal B(K,N)$, $\mathcal C(N)$, $\mathcal D(N)$ and  $\mathcal E(N)$ from \eqref{subbe}, \eqref{cnsub}, \eqref{dnsub} and \eqref{ensub} respectively. Then
\begin{equation*}
    \farey_N- \bigl(\farey_{100} \cup \mathcal A(N) \bigr) = \mathcal B(101,N) \cup \mathcal C(N) \cup \mathcal D(N) \cup \mathcal E(N).
\end{equation*}
Summing $Q_{hk\sigma}(N)$ for $h/k \in \mathcal B(101,N)$  is $O(e^{WN})$ for any $W > \cl(\pi/3)/(6\pi) \approx 0.0538$ by Theorem \ref{ewu}. Since
\begin{equation*}
    -\log|w(1,-3)|\approx 0.0356795, \qquad -\log|w(0,-1)|/2 \approx 0.0340381, \qquad -\log|w(0,-2)|\approx 0.0256706
\end{equation*}
we see from Theorems \ref{thmc}, \ref{c2sed} and \ref{thme} that the sums of $Q_{hk\sigma}(N)$ for $h/k \in \mathcal C(N)$, $\mathcal D(N)$ and  $\mathcal E(N)$   are $O(e^{0.0357N})$, $O(e^{0.0341N})$ and $O(e^{0.0257N})$ respectively. This completes the proof.
\end{proof}

As a final remark, comparing Tables \ref{e1n1} and \ref{c2n1} we notice that $\mathcal E_1(N, 1)$ is almost exactly 3 times the size of $\mathcal C_2(N, 1)$ and that their asymptotic expansions  also seem to match. This is true for other values of $\sigma$ too. From Theorems \ref{c2se} and \ref{thme} we have
\begin{equation}\label{??}
     3 \cdot c_{t}(\sigma) = e_t(\sigma)
\end{equation}
for the first expansion coefficients at $t=0$. Numerically, \eqref{??} seems to be true for all $t$, as we mentioned before in \eqref{simsum3}.

{\small
\bibliography{raddata}

\begin{thebibliography}{Mun08}

\bibitem[And03]{An}
George~E. Andrews.
\newblock Partitions: at the interface of {$q$}-series and modular forms.
\newblock {\em Ramanujan J.}, 7(1-3):385--400, 2003.
\newblock Rankin memorial issues.

\bibitem[DG02]{DG}
Morley Davidson and Stephen~M. Gagola, Jr.
\newblock On {R}ademacher's conjecture and a recurrence relation of {E}uler.
\newblock {\em Quaest. Math.}, 25(3):317--325, 2002.

\bibitem[DG14]{DrGe}
Michael Drmota and Stefan Gerhold.
\newblock Disproof of a conjecture by {R}ademacher on partial fractions.
\newblock {\em Proc. Amer. Math. Soc. Ser. B}, 1:121--134, 2014.

\bibitem[IK04]{IwKo}
Henryk Iwaniec and Emmanuel Kowalski.
\newblock {\em Analytic number theory}, volume~53 of {\em American Mathematical
  Society Colloquium Publications}.
\newblock American Mathematical Society, Providence, RI, 2004.

\bibitem[Max03]{max}
Leonard~C. Maximon.
\newblock The dilogarithm function for complex argument.
\newblock {\em R. Soc. Lond. Proc. Ser. A Math. Phys. Eng. Sci.},
  459(2039):2807--2819, 2003.

\bibitem[Mun08]{Mu}
Augustine~O. Munagi.
\newblock The {R}ademacher conjecture and {$q$}-partial fractions.
\newblock {\em Ramanujan J.}, 15(3):339--347, 2008.

\bibitem[Olv74]{Ol}
F.~W.~J. Olver.
\newblock {\em Asymptotics and special functions}.
\newblock Academic Press [A subsidiary of Harcourt Brace Jovanovich,
  Publishers], New York-London, 1974.
\newblock Computer Science and Applied Mathematics.

\bibitem[O'Sa]{OS1}
Cormac O'Sullivan.
\newblock Asymptotics for the partial fractions of the restricted partition
  generating function {I}.
\newblock Available on the arXiv.

\bibitem[O'Sb]{OS3}
Cormac O'Sullivan.
\newblock Zeros of the dilogarithm.
\newblock Available on the arXiv. To appear in Mathematics of Computation.

\bibitem[O'S15]{OS}
Cormac O'Sullivan.
\newblock On the partial fraction decomposition of the restricted partition
  generating function.
\newblock {\em Forum Math.}, 27(2):735--766, 2015.

\bibitem[Rad37]{Ra2}
Hans Rademacher.
\newblock A convergent series for the partition function p(n).
\newblock {\em Proc. Natl. Acad. Sci. USA}, 23(2):78--84, 1937.

\bibitem[Rad73]{Ra}
Hans Rademacher.
\newblock {\em Topics in analytic number theory}.
\newblock Springer-Verlag, New York, 1973.
\newblock Edited by E. Grosswald, J. Lehner and M. Newman, Die Grundlehren der
  mathematischen Wissenschaften, Band 169.

\bibitem[Rud76]{Ru}
Walter Rudin.
\newblock {\em Principles of mathematical analysis}.
\newblock McGraw-Hill Book Co., New York-Auckland-D\"usseldorf, third edition,
  1976.
\newblock International Series in Pure and Applied Mathematics.

\bibitem[SZ13]{SZ}
Andrew~V. Sills and Doron Zeilberger.
\newblock Rademacher's infinite partial fraction conjecture is (almost
  certainly) false.
\newblock {\em J. Difference Equ. Appl.}, 19(4):680--689, 2013.

\bibitem[Woj06]{Woj}
John Wojdylo.
\newblock On the coefficients that arise from {L}aplace's method.
\newblock {\em J. Comput. Appl. Math.}, 196(1):241--266, 2006.

\bibitem[Zag07]{Zag07}
Don Zagier.
\newblock The dilogarithm function.
\newblock In {\em Frontiers in number theory, physics, and geometry. {II}},
  pages 3--65. Springer, Berlin, 2007.

\end{thebibliography}
}

\textsc{Dept. of Mathematics, The CUNY Graduate Center , New York, NY 10016-4309, U.S.A.}

{\em E-mail address:} \texttt{cosullivan@gc.cuny.edu}

{\em Web page:} \texttt{http://fsw01.bcc.cuny.edu/cormac.osullivan}

\end{document}